\documentclass{amsart}
\usepackage{amsfonts}

\newtheorem{theorem}{Theorem}
\theoremstyle{plain}

\newtheorem{definition}{Definition}

\newtheorem{lemma}{Lemma}

\newtheorem{remark}{Remark}

\numberwithin{equation}{section}
\input{tcilatex}

\begin{document}
\title[Everywhere h\"{o}lder continuity]{Not convex densities and everywhere
h\"{o}lder continuity }
\author{Tiziano Granucci}
\address{ISIS\ Salvemini duca D'Aosta, Firenze, via Giusti 22, 50100, Italy}
\email{tizianogranucci@libero.it, tiziano.granucci@polotecnicofi.com}
\urladdr{https://www.tizianogranucci.com/}
\thanks{Declarations: Data sharing not applicable to this article as no
datasets were generated or analysed during the current study. The author has
no conicts of interest to declare that are relevant to the content of this
article.\\
I would like to thank all my family and friends for the support given to me
over the years: Elisa Cirri, Caterina Granucci, Delia Granucci, Irene
Granucci, Laura and Fiorenza Granucci, Massimo Masi and Monia Randolfi. Also
I would like to thank my professors Luigi Barletti, Giorgio Busoni, Elvira
Mascolo, Giorgio Talenti and Vincenso Vepri }
\date{Sectember 27, 2021}
\dedicatory{Dedicated to the memory of Fiorella Pettini.}
\thanks{This paper is in final form and no version of it will be submitted
for publication elsewhere.}

\begin{abstract}
In this paper we study the everywhere h\"{o}der continuity of the minima of
the following class of vectorial integral functionals 
\begin{equation*}
\int\limits_{\Omega }\sum\limits_{\alpha =1}^{n}f_{\alpha }\left(
x,u^{\alpha },\nabla u^{\alpha }\right) +G\left( x,u,\nabla u\right) \,dx
\end{equation*}%
The energy densities satisfy suitable structure assumptions and may have
neither radial nor quasi-diagonal structure. The regularity of minimizers is
obtained by proving that each component stays in a suitable De Giorgi class
and, from this, we conclude about the h\"{o}lder continuity.
\end{abstract}

\maketitle

\section{Introduction}

In this paper we study the regularity of the local minima of the following
integral functional%
\begin{equation}
J\left( u,\Omega \right) =\int\limits_{\Omega }\sum\limits_{\alpha
=1}^{n}f_{\alpha }\left( x,u^{\alpha }\left( x\right) ,\nabla u^{\alpha
}\left( x\right) \right) +G\left( x,u\left( x\right) ,\nabla u\left(
x\right) \right) \,dx  \label{1.1}
\end{equation}%
where $\Omega $ is a open subset of $%
\mathbb{R}
^{N}$ and $u\in W^{1,p}\left( \Omega ,%
\mathbb{R}
^{n}\right) $ with $n\geq 2$, $N\geq 1$ and $1<p<N$.

Moreover the following hypotheses hold

\begin{description}
\item[H.1.1] For every $\alpha =1,...,n$ the function$\ f_{\alpha }:\Omega
\times 
\mathbb{R}
\times 
\mathbb{R}
^{N}\rightarrow 
\mathbb{R}
$ is a Caratheodory function and the following growth conditions hold%
\begin{equation}
\left\vert \xi ^{\alpha }\right\vert ^{p}-b_{\alpha }\left( x\right)
\left\vert s\right\vert ^{\gamma _{\alpha }}-a_{\alpha }\left( x\right) \leq
f_{\alpha }\left( x,s,\xi ^{\alpha }\right) \leq L_{\alpha }\left(
\left\vert \xi ^{\alpha }\right\vert ^{p}+b_{\alpha }\left( x\right)
\left\vert s\right\vert ^{\gamma _{\alpha }}+a_{\alpha }\left( x\right)
\right)
\end{equation}%
for almost every $x\in \Omega $, for every $s\in 
\mathbb{R}
$ and for every $\xi ^{\alpha }\in 
\mathbb{R}
^{N}$ where $L_{\alpha }>1$, $1<p\leq \gamma _{\alpha }<p^{\ast }=\frac{Np}{%
N-p}$, $b_{\alpha }$ and $a_{\alpha }$ are two not-negative function, $%
b_{\alpha }\in L_{loc}^{\sigma _{\alpha }}\left( \Omega \right) $ and $%
a_{\alpha }\in L_{loc}^{\kappa }\left( \Omega \right) $ with $\sigma
_{\alpha }=\frac{p^{\ast }}{p^{\ast }-\gamma _{\alpha }-\epsilon p^{\ast }}$%
, $\kappa =\frac{N}{p-\epsilon N}$ and $0<\epsilon <\frac{p}{N}$.

\item[H.1.2] $G:\Omega \times 
\mathbb{R}
^{n}\times 
\mathbb{R}
^{Nn}\rightarrow 
\mathbb{R}
$ is a Caratheodory function and the following growth conditions hold%
\begin{equation}
\left\vert G\left( x,u,\xi \right) \right\vert \leq C\left( \left\vert \xi
\right\vert ^{q}+\left\vert u\right\vert ^{q}+a\left( x\right) \right)
\end{equation}%
for almost every $x\in \Omega $, for every $u\in 
\mathbb{R}
^{n}$ and for every $\xi \in 
\mathbb{R}
^{Nn}$ where $C>1$, $1\leq q<\frac{p^{2}}{N}<p$, $a$ is a not-negative
function and $a\in L_{loc}^{\kappa }\left( \Omega \right) $ with $\kappa =%
\frac{N}{p-\epsilon N}$ and $0<\epsilon <\frac{p}{N}$.

\item[H.1.3] The function $G\left( x,u,\cdot \right) $ is rank one convex
then%
\begin{equation*}
G\left( x,u,\lambda \xi ^{1}+\left( 1-\lambda \right) \xi ^{2}\right) \leq
\lambda G\left( x,u,\xi ^{1}\right) +\left( 1-\lambda \right) G\left(
x,u,\xi ^{2}\right)
\end{equation*}%
for a. e. $x\in \Omega $, for every $u\in 
\mathbb{R}
^{n}$, for every $\lambda \in \left[ 0,1\right] $, and for every $\xi
^{1},\xi ^{2}\in 
\mathbb{R}
^{Nn}$ with $rank\left\{ \xi ^{1}-\xi ^{2}\right\} \leq 1$.

\item[H.1.4] The function $G\left( x,\cdot ,\xi \right) $ is h\"{o}lder
continuous and for almost every $x\in \Omega $ and for every $\xi \in 
\mathbb{R}
^{Nn}$ it follows 
\begin{equation}
\left\vert G\left( x,u,\xi \right) -G\left( x,v,\xi \right) \right\vert \leq
c\left( x\right) \left\vert \xi \right\vert ^{\delta }\left\vert
u-v\right\vert ^{\beta }  \label{H.1.4}
\end{equation}%
for a. e. $x\in \Omega $, for every $u,v\in 
\mathbb{R}
^{n}$ and for every $\xi \in 
\mathbb{R}
^{Nn}$ with $0<\delta <q$ and $0<\beta <\min \left\{ \frac{p^{\ast }\left(
q-\delta \right) }{q},\,p-\delta ,1\right\} $, $c$ is a not-negative
function and $c\left( x\right) \in L_{loc}^{\sigma }\left( \Omega \right) $
with $\sigma >\frac{p^{\ast }q}{p^{\ast }\left( q-\delta \right) -\beta q}$
and 
\begin{equation}
\frac{\delta }{p}+\frac{\beta }{p^{\ast }}+\frac{1}{\sigma }<\frac{p}{N}
\label{H.1.4.5}
\end{equation}

\item[H.1.4 (bis)] H.1.4 holds with $c\left( x\right) \in L_{loc}^{\infty
}\left( \Omega \right) $ and 
\begin{equation}
\frac{\delta }{p}+\frac{\beta }{p^{\ast }}<\frac{p}{N}  \label{H.1.4.6}
\end{equation}
\end{description}

\bigskip

\begin{remark}
\label{rem1} Using hypotheses H.1.2 and H.1.4 we get%
\begin{equation*}  \label{rem1.1}
\begin{tabular}{ll}
$\left\vert G\left( x,u,\xi \right) \right\vert $ & $\leq \left\vert G\left(
x,u,\xi \right) -G\left( x,0,\xi \right) \right\vert +\left\vert G\left(
x,0,\xi \right) \right\vert $ \\ 
& $\leq c\left( x\right) \left\vert \xi \right\vert ^{\delta }\left\vert
u\right\vert ^{\beta }+C\left( \left\vert \xi \right\vert ^{q}+a\left(
x\right) \right) $%
\end{tabular}%
\end{equation*}%
for a. e. $x\in \Omega $, for every $u\in 
\mathbb{R}
^{n}$\ and for every $\xi \in 
\mathbb{R}
^{Nn}$. Then the following inequality holds 
\begin{equation}  \label{rem1.2}
\left\vert G\left( x,u,\xi \right) \right\vert \leq c\left( x\right)
\left\vert \xi \right\vert ^{\delta }\left\vert u\right\vert ^{\beta
}+C\left( \left\vert \xi \right\vert ^{q}+a\left( x\right) \right)
\end{equation}%
for a. e. $x\in \Omega $, for every $u\in 
\mathbb{R}
^{n}$\ and for every $\xi \in 
\mathbb{R}
^{Nn}$.
\end{remark}

Our principal result is the following Theorem.

\begin{theorem}
\label{TH1} If $u\in W_{loc}^{1,p}\left( \Omega ,%
\mathbb{R}
^{n}\right) $, with $N\geq 2$, $n\geq 1$ and $1<p<N$, is a local minimizer
of the functional $J\left( u,\Omega \right) $ and the hypotheses H.1.1,
H.1.2, H.1.3 and H.1.4 (or H.1.4 (bis)) hold then every componets $u^{\alpha
}$ of the vectorial function $u$ are a locally h\"{o}lder continuous
functions.
\end{theorem}

Theorem 1 derives from the following two results.\bigskip\ 

\begin{theorem}
\label{TH2} If $u\in W_{loc}^{1,p}\left( \Omega ,%
\mathbb{R}
^{n}\right) $ with $1<p<N$ is a local minimizer of the functional $J\left(
u,\Omega \right) $ and the hypotheses H.1.1, H.1.2, H.1.3 and H.1.4 (or
H.1.4 (bis)) hold then two positive constants $C_{C,1}$, $C_{C,2}$ and a
radius $R_{0}>0$ exists such that for every $0<\varrho <R<R_{0}$ and for
every $k\in 
\mathbb{R}
$\ it follows%
\begin{equation}
\begin{tabular}{l}
$\int\limits_{A_{k,\varrho }^{\alpha }}\left\vert \nabla u^{\alpha
}\right\vert ^{p}\,dx\leq \frac{C_{C,1}}{\left( R-\varrho \right) ^{p}}%
\int\limits_{A_{k,R}^{\alpha }}\left( u^{\alpha }-k\right)
^{p}\,dx+C_{C,2}\left( 1+\left\vert k\right\vert ^{p}R^{-\epsilon _{\alpha
}N}\right) \left\vert A_{k,R}^{\alpha }\right\vert ^{1-\frac{p}{N}+\epsilon
} $%
\end{tabular}%
\end{equation}%
and%
\begin{equation}
\int\limits_{B_{k,\varrho }^{\alpha }}\left\vert \nabla u^{\alpha
}\right\vert ^{p}\,dx\leq \frac{C_{C,1}}{\left( R-\varrho \right) ^{p}}%
\int\limits_{B_{k,R}^{\alpha }}\left( u^{\alpha }-k\right)
^{p}\,dx+C_{C,2}\left( 1+\left\vert k\right\vert ^{p}R^{-\epsilon _{\alpha
}N}\right) \left\vert B_{k,R}^{\alpha }\right\vert ^{1-\frac{p}{N}+\epsilon }
\end{equation}%
where $A_{k,s}^{\alpha }=\left\{ u^{\alpha }>k\right\} \cap B_{s}\left(
x_{0}\right) $ and $B_{k,s}^{\alpha }=\left\{ u^{\alpha }<k\right\} \cap
B_{s}\left( x_{0}\right) $.\bigskip
\end{theorem}

\begin{theorem}
\label{TH3} If $u\in W_{loc}^{1,p}\left( \Omega ,%
\mathbb{R}
^{n}\right) $, with $N\geq 2$, $n\geq 1$ and $1<p<N$, is a local minimizer
of the functional $J\left( u,\Omega \right) $ and the hypotheses H.1.1,
H.1.2, H.1.3 and H.1.4 (or H.1.4 (bis)) hold then every componets $u^{\alpha
}$ of the vector function $u$ are a locally boundedness functions.
\end{theorem}

\bigskip

Furthermore, in this article, we will also consider the polyconvex case.

\begin{description}
\item[H.2.3] $G:\Omega \times \ 
\mathbb{R}
^{N\times n}\rightarrow R$ is a Carath\'{e}odory function defined as%
\begin{equation*}
G(x,u,\xi )=\sum\limits_{\alpha =1}^{n}G_{\alpha }\left( x,u,\left(
adj_{m-1}\xi \right) ^{\alpha }\right)
\end{equation*}%
here $G_{\alpha }:\Omega \times \ 
\mathbb{R}
^{\frac{N!}{N!\left( N-n\right) !}}\rightarrow 
\mathbb{R}
$ is a Carath\'{e}odory convex function satisfying the following growth
conditions%
\begin{equation*}
0\leq G_{\alpha }\left( x,0,\left( adj_{n-1}\xi \right) ^{\alpha }\right)
\leq k_{4}|\left( adj_{n-1}\xi \right) ^{\alpha }|^{r}+a\left( x\right)
\end{equation*}%
for every $\xi \in 
\mathbb{R}
^{N\times n}$, for almost every $x\in \Omega $, here $k_{3}$ is a real
positive constant, $1\leq r<\frac{p}{n-1}$ and $a\in L_{loc}^{\sigma }\left(
\Omega \right) $ is a non negative function, $N\geq n\geq 2$.

\item[H.2.4] The function $G_{\alpha }\left( x,\cdot ,\xi \right) $ is h\"{o}%
lder continuous and for almost every $x\in \Omega $ and for every $\xi \in 
\mathbb{R}
^{Nn}$ it follows 
\begin{equation}
\left\vert G_{\alpha }\left( x,u,\left( adj_{n-1}\xi \right) ^{\alpha
}\right) -G_{\alpha }\left( x,v,\left( adj_{n-1}\xi \right) ^{\alpha
}\right) \right\vert \leq c\left( x\right) \left\vert \left( adj_{n-1}\xi
\right) ^{\alpha }\right\vert ^{\delta }\left\vert u-v\right\vert ^{\beta }
\end{equation}%
for a. e. $x\in \Omega $, for every $u,v\in 
\mathbb{R}
^{n}$ and for every $\xi \in 
\mathbb{R}
^{Nn}$ with $0<\delta <\frac{r}{n-1}$ and $0<\beta <\min \left\{ \frac{%
p^{\ast }\left( r-\left( n-1\right) \delta \right) }{r},1\right\} $, $c$ is
a not-negative function and $c\left( x\right) \in L_{loc}^{\sigma }\left(
\Omega \right) $ with $\sigma >\frac{p^{\ast }r}{p^{\ast }\left( r-\left(
n-1\right) \delta \right) -\beta r}$ and 
\begin{equation}
\frac{(n-1)\delta }{p}+\frac{\beta }{p^{\ast }}+\frac{1}{\sigma }<\frac{p}{N}
\end{equation}%
for $\alpha =1,...,n$.

\item[H.2.4 (bis)] H.1.4 holds with $c\left( x\right) \in L_{loc}^{\infty
}\left( \Omega \right) $ and 
\begin{equation}
\frac{(n-1)\delta }{p}+\frac{\beta }{p^{\ast }}<\frac{p}{N}  \label{H.2.4.6}
\end{equation}%
for $\alpha =1,...,n$.
\end{description}

\bigskip

\begin{remark}
Using hypotheses H.1.2 and H.1.4 we get%
\begin{equation*}
\begin{tabular}{ll}
$\left\vert G_{\alpha }\left( x,u,\left( adj_{n-1}\xi \right) ^{\alpha
}\right) \right\vert $ & $\leq \left\vert G_{\alpha }\left( x,u,\left(
adj_{n-1}\xi \right) ^{\alpha }\right) -G_{\alpha }\left( x,0,\left(
adj_{n-1}\xi \right) ^{\alpha }\right) \right\vert +\left\vert G_{\alpha
}\left( x,0,\left( adj_{n-1}\xi \right) ^{\alpha }\right) \right\vert $ \\ 
& $\leq c\left( x\right) \left\vert \left( adj_{n-1}\xi \right) ^{\alpha
}\right\vert ^{\delta }\left\vert u\right\vert ^{\beta }+C\left( \left\vert
\left( adj_{n-1}\xi \right) ^{\alpha }\right\vert ^{r}+a\left( x\right)
\right) $%
\end{tabular}%
\end{equation*}%
for a. e. $x\in \Omega $, for every $u\in 
\mathbb{R}
^{n}$\ and for every $\xi \in 
\mathbb{R}
^{Nn}$. Then the following inequality holds 
\begin{equation}
\left\vert G_{\alpha }\left( x,u,\left( adj_{n-1}\xi \right) ^{\alpha
}\right) \right\vert \leq c\left( x\right) \left\vert \left( adj_{n-1}\xi
\right) ^{\alpha }\right\vert ^{\delta }\left\vert u\right\vert ^{\beta
}+C\left( \left\vert \left( adj_{n-1}\xi \right) ^{\alpha }\right\vert
^{r}+a\left( x\right) \right)  \label{rem2.2}
\end{equation}%
for a. e. $x\in \Omega $, for every $u\in 
\mathbb{R}
^{n}$\ and for every $\xi \in 
\mathbb{R}
^{Nn}$.
\end{remark}

\bigskip Our second regular result is the following Theorem

\begin{theorem}
If $u\in W_{loc}^{1,p}\left( \Omega ,%
\mathbb{R}
^{n}\right) $, with $N\geq 2$, $n\geq 1$ and $1<p<N$, is a local minimizer
of the functional $J\left( u,\Omega \right) $ and the hypotheses H.1.1,
H.1.2, H.2.3 and H.2.4 (or H.2.4 (bis)) hold then every componets $u^{\alpha
}$ of the vectorial function $u$ are a locally h\"{o}lder continuous
functions.
\end{theorem}

Theorem 4 derives from the following two theorems.\bigskip\ 

\begin{theorem}
If $u\in W_{loc}^{1,p}\left( \Omega ,%
\mathbb{R}
^{n}\right) $ with $1<p<N$ is a local minimizer of the functional $J\left(
u,\Omega \right) $ and the hypotheses H.1.1, H.1.2, H.2.3 and H.2.4 (or
H.2.4 (bis)) hold then two positive constants $C_{C,1}$, $C_{C,2}$ and a
radius $R_{0}>0$ exists such that for every $0<\varrho <R<R_{0}$ and for
every $k\in 
\mathbb{R}
$\ it follows%
\begin{equation}
\begin{tabular}{l}
$\int\limits_{A_{k,\varrho }^{\alpha }}\left\vert \nabla u^{\alpha
}\right\vert ^{p}\,dx\leq \frac{C_{C,1}}{\left( R-\varrho \right) ^{p}}%
\int\limits_{A_{k,R}^{\alpha }}\left( u^{\alpha }-k\right)
^{p}\,dx+C_{C,2}\left( 1+\left\vert k\right\vert ^{p}R^{-\epsilon _{\alpha
}N}\right) \left\vert A_{k,R}^{\alpha }\right\vert ^{1-\frac{p}{N}+\epsilon
} $%
\end{tabular}%
\end{equation}%
and%
\begin{equation}
\int\limits_{B_{k,\varrho }^{\alpha }}\left\vert \nabla u^{\alpha
}\right\vert ^{p}\,dx\leq \frac{C_{C,1}}{\left( R-\varrho \right) ^{p}}%
\int\limits_{B_{k,R}^{\alpha }}\left( u^{\alpha }-k\right)
^{p}\,dx+C_{C,2}\left( 1+\left\vert k\right\vert ^{p}R^{-\epsilon _{\alpha
}N}\right) \left\vert B_{k,R}^{\alpha }\right\vert ^{1-\frac{p}{N}+\epsilon }
\end{equation}%
where $A_{k,s}^{\alpha }=\left\{ u^{\alpha }>k\right\} \cap B_{s}\left(
x_{0}\right) $, $B_{k,s}^{\alpha }=\left\{ u^{\alpha }<k\right\} \cap
B_{s}\left( x_{0}\right) $. and $\alpha =1,...,m$.
\end{theorem}

\begin{theorem}
If $u\in W_{loc}^{1,p}\left( \Omega ,%
\mathbb{R}
^{n}\right) $, with $N\geq 2$, $n\geq 1$ and $1<p<N$, is a local minimizer
of the functional $J\left( u,\Omega \right) $ and the hypotheses H.1.1,
H.1.2, H.2.3 and H.2.4 (or H.2.4 (bis)) hold then every componets $u^{\alpha
}$ of the vector function $u$ are a locally boundedness functions.
\end{theorem}

The results stated in the previous theorems [1-6] are not trivial, in
particular we recall the following fundamental counterexample by De Giorgi
[14]; the function 
\begin{equation*}
u(x)=x\left\vert x\right\vert ^{-\delta }
\end{equation*}%
with $\delta =\frac{N}{2}\left( 1-\frac{1}{\sqrt{1+\left( 2N-2\right) ^{2}}}%
\right) $ it is a local minimum of the functional%
\begin{equation}
F\left( u,B_{1}\right) =\int\limits_{B_{1}}\sum\limits_{\alpha ,\beta
=1}^{N}\sum\limits_{i,j=1}^{N}A_{\alpha ,\beta }^{i,j}(x)\partial
_{i}u^{\alpha }\partial _{j}u^{\beta }\,dx
\end{equation}%
with%
\begin{equation*}
A_{\alpha ,\beta }^{i,j}(x)=\delta _{\alpha ,\beta }\delta _{i,j}+\left[
\left( N-2\right) \delta _{\alpha ,i}+N\frac{x_{\alpha }x_{i}}{\left\vert
x\right\vert ^{2}}\right] \left[ \left( N-2\right) \delta _{\beta ,j}+N\frac{%
x_{\beta }x_{j}}{\left\vert x\right\vert ^{2}}\right]
\end{equation*}%
$N\geq 3$ and $B_{1}=\left\{ x\in 
\mathbb{R}
^{N}:\left\vert x\right\vert \leq 1\right\} $. It is easily observed that $%
u(x)=x\left\vert x\right\vert ^{-\delta }$ it is not a bounded function.
Recently, in [8-10, 25, 26] new classes of vectorial problems, with regular
weak solutions, have been introduced. In particular in [8], Cupini, Focardi,
Leonetti and Mascolo introduced the following class of vectorial functionals%
\begin{equation}
\int\limits_{\Omega }f\left( x,\nabla u\right) \,dx  \label{1.3}
\end{equation}%
where $\Omega \subset 
\mathbb{R}
^{N},u:\Omega \rightarrow 
\mathbb{R}
^{n},N>1,n\geq 1$ and%
\begin{equation*}
f\left( x,\nabla u\right) =\sum\limits_{\alpha =1}^{n}F_{\alpha }\left(
x,\nabla u^{\alpha }\right) +G\left( x,\nabla u\right)
\end{equation*}%
where $F_{\alpha }\times \ R^{N\times n}\rightarrow R$ is a Carath\'{e}odory
function satisfying the following standard growth condition%
\begin{equation*}
k_{1}\left\vert \xi ^{\alpha }\right\vert ^{p}-a\left( x\right) \leq
F_{\alpha }\left( x,\xi ^{\alpha }\right) \leq k_{2}\left\vert \xi ^{\alpha
}\right\vert ^{p}+a\left( x\right)
\end{equation*}%
for every $\xi ^{\alpha }\in 
\mathbb{R}
^{N}$ and for almost every $x\in \Omega $ ,where $k_{1}$ and $k_{2}$ are two
real positive constants, $p>1$ and $a\in L_{loc}^{\sigma }\left( \Omega
\right) $ is a non negative function. In [8],Cupini, Focardi, Leonetti and
Mascolo analyze two different types of hypotheses on the $G$ function. They
started by assuming that $G:\Omega \times \ R^{N\times n}\rightarrow R$ is a
Carath\'{e}odory rank one convex function satisfying the following growth
condition%
\begin{equation*}
|G(x,\xi )|\leq k_{3}|\xi |^{q}+b(x)
\end{equation*}%
for every $\xi \in 
\mathbb{R}
^{N\times n}$, for almost every $x\in \Omega $, here $k_{3}$ is a real
positive constant, $1\leq q<p$ and $b\in L_{loc}^{\sigma }\left( \Omega
\right) $ is a nonnegative function. Moreover Cupini, Focardi, Leonetti and
Mascolo in [8] study the case where $n\geq m\geq 3$, and $G:\Omega \times \ 
\mathbb{R}
^{N\times n}\rightarrow R$ is a Carath\'{e}odory function defined as%
\begin{equation*}
G(x,\xi )=\sum\limits_{\alpha =1}^{n}G_{\alpha }\left( x,\left( adj_{n-1}\xi
\right) ^{\alpha }\right)
\end{equation*}%
here $G_{\alpha }:\Omega \times \ 
\mathbb{R}
^{\frac{N!}{N!\left( N-n\right) !}}\rightarrow 
\mathbb{R}
$ is a Carath\'{e}odory convex function satisfying the following growth
conditions%
\begin{equation*}
0\leq G_{\alpha }\left( x,\left( adj_{n-1}\xi \right) ^{\alpha }\right) \leq
k_{4}|\left( adj_{n-1}\xi \right) ^{\alpha }|^{r}+b(x)
\end{equation*}%
for every $\xi \in 
\mathbb{R}
^{N\times n}$, for almost every $x\in \Omega $, here $k_{3}$ is a real
positive constant, $1\leq r<p$ and $b\in L_{loc}^{\sigma }\left( \Omega
\right) $ is a non negative function. In both cases, by imposing appropriate
hypotheses on the parameters $q$ and $r$ , Cupini, Focardi, Leonetti and
Mascolo proved that the localminimizers of the vectorial functional (\ref%
{1.3}) are locally h\"{o}lder continuous functions. In [25] the author with
M. Randolfi proved a regularity result for the minima of vector functionals
with anisotropic growths of the following type%
\begin{equation}
\int\limits_{\Omega }\sum\limits_{\alpha =1}^{n}F_{\alpha }\left( x,\nabla
u^{\alpha }\right) +G\left( x,\nabla u\right) \,dx  \label{1.4}
\end{equation}%
with%
\begin{equation}
\sum\limits_{\alpha =1}^{n}\Phi _{i,\alpha }\left( \left\vert \xi
_{i}^{\alpha }\right\vert \right) \leq F_{\alpha }\left( x,\xi ^{\alpha
}\right) \leq L\left[ \bar{B}_{\alpha }^{\beta _{\alpha }}\left( \left\vert
\xi ^{\alpha }\right\vert \right) +a\left( x\right) \right]  \label{1.5}
\end{equation}%
where $\Phi _{i,\alpha }$ are N functions belonging to the class $\triangle
_{2}^{m_{\alpha }}\cap \nabla _{2}^{r_{\alpha }}$, $\bar{B}_{\alpha }$\ is
the Sobolev function associated with $\Phi _{i,\alpha }$'s, $\beta _{\alpha
}\in \left( 0,1\right] $ and $a\in L_{loc}^{\sigma }\left( \Omega \right) $
is a non negative function; oppure with%
\begin{equation}
\sum\limits_{\alpha =1}^{n}\Phi _{i,\alpha }\left( \left\vert \xi
_{i}^{\alpha }\right\vert \right) -a\left( x\right) \leq F_{\alpha }\left(
x,\xi ^{\alpha }\right) \leq L_{1}\left[ \sum\limits_{\alpha =1}^{n}\Phi
_{i,\alpha }\left( \left\vert \xi _{i}^{\alpha }\right\vert \right) +a\left(
x\right) \right]  \label{1.6}
\end{equation}%
where $\Phi _{i,\alpha }$ are N-functions belonging to the class $\triangle
_{2}^{m_{\alpha }}\cap \nabla _{2}^{r_{\alpha }}$ and $a\in L_{loc}^{\sigma
}\left( \Omega \right) $ is a non negative function, moreover, appropriate
hypotheses are made on the density $G$, for more details we refer to [25].
In particular, using the techniques presented in [25,....], the author with
M. Randolfi have shown that the minima of the functional (\ref{1.4}) are
locally bounded functions in the case (\ref{1.5}) and locally Holder
continuous in the case (\ref{1.6}), we refer to [25] for more details. In
recent years a large number of articles have been written dealing with the
regularity of weak solutions of vactorial problems under standard growths,
general growths and anisotropic growths, refer to [1, 2, 4-12, 15-22, 23,
25-32, 35-37].

Our results can therefore be framed within a vast area of {}{}research
called everywhere regulairy that was born with the fundamental work
Uhlenbeck [37] and which has developed up to now with results that generally
have concerned functionals of the form 
\begin{equation*}
\int\limits_{B_{R}\left( x_{0}\right) }G\left( \left\vert \nabla u\left(
x\right) \right\vert \right) \,dx
\end{equation*}%
The literature in this area is very wide, in the bibliography we report only
some articles. Generally results are obtained on the boudedness of the
gradient of the minima using particular conditions of regularity and
ellipticity using very different techniques, we mention [1,2, 4-7, 15,
17-21, 27-29, 35-37] only because they are those consulted by the author
during the preparation of the article but we refer to the bibliographies of
the aforementioned articles for more bibliographical references. Our results
differ from those previously mentioned as it deals with some functionals of
the form%
\begin{equation*}
\int\limits_{B_{R}\left( x_{0}\right) }F\left( x,u,\nabla u\right) \,dx
\end{equation*}%
in this case, as far as the limited knowledge of the author is concerned,
there are few results, in particular we refer to [8-12] and [25]. Theorem 1
therefore generalizes the results presented in [8], [25] and is, in a
certain sense, complementary to the results presented by the author in [26].

\section{\protect\bigskip Preliminary results}

Before giving the proofs of Theorem 1 and Theorem 2, for completeness we
introduce a list of results that we will use during the proof.

\subsection{Lemmata}

\begin{lemma}[Young Inequality]
Let $\varepsilon >0$, $a,b>0$ and $1<p,q<+\infty $ with $\frac{1}{p}+\frac{1%
}{q}=1$\ then it follows%
\begin{equation}
ab\leq \varepsilon \frac{a^{p}}{p}+\frac{b^{q}}{\varepsilon ^{\frac{q}{p}}q}
\end{equation}
\end{lemma}

\begin{lemma}[H\"{o}lder Inequality]
Assume $1\leq p,q\leq +\infty $ with $\frac{1}{p}+\frac{1}{q}=1$\ then if $%
u\in L^{p}\left( \Omega \right) $\ and $v\in L^{p}\left( \Omega \right) $\
it follows%
\begin{equation}
\int\limits_{\Omega }\left\vert uv\right\vert \,dx\leq \left(
\int\limits_{\Omega }\left\vert u\right\vert ^{p}\,dx\right) ^{\frac{1}{p}%
}\left( \int\limits_{\Omega }\left\vert v\right\vert ^{q}\,dx\right) ^{\frac{%
1}{q}}
\end{equation}
\end{lemma}

\begin{lemma}
\label{lem3} Let $Z\left( t\right) $ be a nonnegative and bounded function
on the set $\left[ \varrho ,R\right] $; if for every $\varrho \leq t<s\leq R$
we get%
\begin{equation}
Z\left( t\right) \leq \theta Z\left( s\right) +\frac{A}{\left( s-t\right)
^{\lambda }}+\frac{B}{\left( s-t\right) ^{\mu }}+C
\end{equation}%
where $A,B,C\geq 0$, $\lambda >\mu >0$ and $0\leq \theta <1$ then it follows%
\begin{equation}
Z\left( \varrho \right) \leq C\left( \theta ,\lambda \right) \left( \frac{A}{%
\left( R-\varrho \right) ^{\lambda }}+\frac{B}{\left( R-\varrho \right)
^{\mu }}+C\right)
\end{equation}%
where $C\left( \theta ,\lambda \right) >0$ is a real constant depending only
on $\theta $ and $\lambda $.
\end{lemma}

Definitions

\begin{definition}
\bigskip A function $f:%
\mathbb{R}
^{nm}\rightarrow 
\mathbb{R}
\cup \left\{ +\infty \right\} $ is said to be rank one convex if%
\begin{equation*}
f\left( \lambda A+\left( 1-\lambda \right) B\right) \leq \lambda f\left(
A\right) +\left( 1-\lambda \right) f\left( B\right)
\end{equation*}%
for every $\lambda \in \left[ 0,1\right] $, $A$,$B\in 
\mathbb{R}
^{nm}$ with $rank\left\{ A-B\right\} \leq 1$.
\end{definition}

\begin{definition}
A Borel measurable function and locally integrable function $f:%
\mathbb{R}
^{nm}\rightarrow 
\mathbb{R}
$ is said to be quasiconvex if%
\begin{equation*}
f\left( A\right) \leq \frac{1}{\left\vert D\right\vert }\int\limits_{D}f%
\left( A+\nabla \varphi \right) \,dx
\end{equation*}%
for every bounded domain $D\subset 
\mathbb{R}
^{n}$, for every $A\in 
\mathbb{R}
^{nm}$ and for every $\varphi \in W_{0}^{1,\infty }\left( D;%
\mathbb{R}
^{nm}\right) $.
\end{definition}

\begin{definition}
A function $f:%
\mathbb{R}
^{nm}\rightarrow 
\mathbb{R}
\cup \left\{ +\infty \right\} $ is said to be polyconvex if there exists a
function $g:%
\mathbb{R}
^{nm}\rightarrow 
\mathbb{R}
\cup \left\{ +\infty \right\} $ convex such that%
\begin{equation*}
f\left( A\right) =g\left( T\left( A\right) \right)
\end{equation*}%
where $T:%
\mathbb{R}
^{nm}\rightarrow 
\mathbb{R}
^{\tau \left( n,m\right) }$ is such that%
\begin{equation*}
T(A)=(A,adj_{2}\left( A\right) ,...,adj_{n\wedge m}\left( A\right) )
\end{equation*}%
where $adj_{s}\left( A\right) $ stands for the matrix of all $s\times s$
minors of tha matrix $A\in 
\mathbb{R}
^{nm}$, 2$\leq s\leq n\wedge m=\min \left\{ n,m\right\} $ and%
\begin{equation*}
\tau \left( n,m\right) =\sum\limits_{s=1}^{n\wedge m}\sigma \left( s\right)
\end{equation*}%
where $\sigma \left( s\right) =\frac{n!m!}{\left( s!\right) ^{2}\left(
m-s\right) !\left( n-s\right) !}$.
\end{definition}

In particular we recall the following theorem.

\begin{theorem}

\begin{enumerate}
\item Let $f:%
\mathbb{R}
^{nm}\rightarrow 
\mathbb{R}
$ then%
\begin{equation*}
f\text{ convex }\Longrightarrow f\text{ polyconvex}\Longrightarrow f\text{
quasiconvex}\Longrightarrow f\text{ rank one convex. }
\end{equation*}

\item If $m=1$ or $n=1$ then all these notions are equivalent.

\item If $f\in C^{2}\left( 
\mathbb{R}
^{nm}\right) $ then rank 0ne convexity is equivalent to Legendre-Hadamard
condition%
\begin{equation*}
\sum\limits_{i,j=1}^{m}\sum\limits_{\alpha ,\beta =1}^{n}\frac{\partial ^{2}f%
}{\partial A_{\alpha }^{i}\partial A_{\beta }^{j}}\left( A\right) \lambda
^{i}\lambda ^{j}\mu _{\alpha }\mu _{\beta }\geq 0
\end{equation*}%
for every $\lambda \in 
\mathbb{R}
^{m}$, $\mu \in 
\mathbb{R}
^{n}$, $A=\left( A_{\alpha }^{i}\right) _{1\leq i\leq m,1\leq \alpha \leq
n}\in 
\mathbb{R}
^{nm}$.

\item If $f:%
\mathbb{R}
^{nm}\rightarrow 
\mathbb{R}
$ is convex, polyconvex, quasiconvex or rank one convex then $f$ is locally
Lipschitz.
\end{enumerate}
\end{theorem}

\subsection{Sobolev Spaces}

\begin{theorem}[Sobolev Inequality]
Let $\Omega $ be a open subset of $%
\mathbb{R}
^{N}$ if $u\in W_{0}^{1,p}\left( \Omega \right) $ with $1\leq p<N$ there
exists a real positive constant $C_{SN}$, depending only on $p$ and $N$,
such that%
\begin{equation}
\left\Vert u\right\Vert _{L^{p^{\ast }}\left( \Omega \right) }\leq
C_{SN}\left\Vert \nabla u\right\Vert _{L^{p}\left( \Omega \right) }
\end{equation}%
where $p^{\ast }=\frac{Np}{N-p}$.
\end{theorem}

\begin{theorem}
(Rellich-Sobolev Immersion Theorem) Let $\Omega $ be a open bounded subset
of $%
\mathbb{R}
^{N}$ with lipschitz boundary then if $u\in W^{1,p}\left( \Omega \right) $
with $1\leq p<N$ there exists a real positive constant $C_{IS}$, depending
only on $p$ and $N$, such that%
\begin{equation}
\left\Vert u\right\Vert _{L^{p^{\ast }}\left( \Omega \right) }\leq
C_{IS}\left\Vert u\right\Vert _{W^{1,p}\left( \Omega \right) }
\end{equation}%
where $p^{\ast }=\frac{Np}{N-p}$.
\end{theorem}

For completeness we remember that if $\Omega $ is a open subset of $%
\mathbb{R}
^{N}$ and $u$\ is a Lebesgue measurable function then $L^{p}\left( \Omega
\right) $ is the set of the class of the Lebesgue measurable function such
that $\int\limits_{\Omega }\left\vert u\right\vert ^{p}\,dx<+\infty $ and $%
W^{1,p}\left( \Omega \right) $\ is the set of the function $u\in L^{p}\left(
\Omega \right) $ such that its waek derivate $\partial _{i}u\in L^{p}\left(
\Omega \right) $. The spaces $L^{p}\left( \Omega \right) $ and $%
W^{1,p}\left( \Omega \right) $ are Banach spaces with the respective norms 
\begin{equation}
\left\Vert u\right\Vert _{L^{p}\left( \Omega \right) }=\left(
\int\limits_{\Omega }\left\vert u\right\vert ^{p}\,dx\right) ^{\frac{1}{p}}
\end{equation}%
and%
\begin{equation}
\left\Vert u\right\Vert _{W^{1,p}\left( \Omega \right) }=\left\Vert
u\right\Vert _{L^{p}\left( \Omega \right) }+\sum\limits_{i=1}^{N}\left\Vert
\partial _{i}u\right\Vert _{L^{p}\left( \Omega \right) }
\end{equation}%
We say that the function $u:\Omega \subset 
\mathbb{R}
^{N}\rightarrow 
\mathbb{R}
^{n}$ belong \ \ \ in $W^{1,p}\left( \Omega ,%
\mathbb{R}
^{n}\right) $ if $u^{\alpha }\in W^{1,p}\left( \Omega \right) $ for every $%
\alpha =1,...,n$, where $u^{\alpha }$ is the $\alpha $ component of the
vector-valued function $u$; we end by remembering that $W^{1,p}\left( \Omega
,%
\mathbb{R}
^{n}\right) $ is a Banach space with the norm%
\begin{equation}
\left\Vert u\right\Vert _{W^{1,p}\left( \Omega ,%
\mathbb{R}
^{n}\right) }=\sum\limits_{\alpha =1}^{n}\left\Vert u^{\alpha }\right\Vert
_{W^{1,p}\left( \Omega \right) }
\end{equation}

\begin{definition}
Let $\Omega \subset 
\mathbb{R}
^{N}$ be a bounded open set and $v:\Omega \rightarrow 
\mathbb{R}
$, we say that $v\in W_{loc}^{1,p}\left( \Omega \right) $ belong to the De
Giorgi class $DG^{+}\left( \Omega ,p,\lambda ,\lambda _{\ast },\chi
,\varepsilon ,R_{0},k_{0}\right) $ with $p>1$, $\lambda >0$, $\lambda _{\ast
}>0$, $\chi >0$, $\varepsilon >0$, $R_{0}>0$ and $k_{0}\geq 0$ if%
\begin{equation}
\int\limits_{A_{k,\varrho }}\left\vert \nabla v\right\vert ^{p}\,dx\leq 
\frac{\lambda }{\left( R-\varrho \right) ^{p}}\int\limits_{A_{k,R}}\left(
v-k\right) ^{p}\,dx+\lambda _{\ast }\left( \chi ^{p}+k^{p}R^{-N\varepsilon
}\right) \left\vert A_{k,R}\right\vert ^{1-\frac{p}{N}+\varepsilon }
\end{equation}%
for all $k\geq k_{0}\geq 0$ and for all pair of balls $B_{\varrho }\left(
x_{0}\right) \subset B_{R}\left( x_{0}\right) \subset \subset \Omega $ with $%
0<\varrho <R<R_{0}$ and $A_{k,s}=B_{s}\left( x_{0}\right) \cap \left\{
v>k\right\} $ with $s>0$.
\end{definition}

\begin{definition}
Let $\Omega \subset 
\mathbb{R}
^{N}$ be a bounded open set and $v:\Omega \rightarrow 
\mathbb{R}
$, we say that $v\in W_{loc}^{1,p}\left( \Omega \right) $ belong to the De
Giorgi class $DG^{-}\left( \Omega ,p,\lambda ,\lambda _{\ast },\chi
,\varepsilon ,R_{0},k_{0}\right) $ with $p>1$, $\lambda >0$, $\lambda _{\ast
}>0$, $\chi >0$ and $k_{0}\geq 0$ if%
\begin{equation}
\int\limits_{B_{k,\varrho }}\left\vert \nabla v\right\vert ^{p}\,dx\leq 
\frac{\lambda }{\left( R-\varrho \right) ^{p}}\int\limits_{B_{k,R}}\left(
k-v\right) ^{p}\,dx+\lambda _{\ast }\left( \chi ^{p}+\left\vert k\right\vert
^{p}R^{-N\varepsilon }\right) \left\vert B_{k,R}\right\vert ^{1-\frac{p}{N}%
+\varepsilon }
\end{equation}%
for all $k\leq -k_{0}\leq 0$ and for all pair of balls $B_{\varrho }\left(
x_{0}\right) \subset B_{R}\left( x_{0}\right) \subset \subset \Omega $ with $%
0<\varrho <R<R_{0}$ and $B_{k,s}=B_{s}\left( x_{0}\right) \cap \left\{
v<k\right\} $ with $s>0$.
\end{definition}

\begin{definition}
We set $DG\left( \Omega ,p,\lambda ,\lambda _{\ast },\chi ,\varepsilon
,R_{0},k_{0}\right) =DG^{+}\left( \Omega ,p,\lambda ,\lambda _{\ast },\chi
,\varepsilon ,R_{0},k_{0}\right) \cap DG^{-}\left( \Omega ,p,\lambda
,\lambda _{\ast },\chi ,\varepsilon ,R_{0},k_{0}\right) $.
\end{definition}

\begin{theorem}
\label{TH6} Let $v\in DG\left( \Omega ,p,\lambda ,\lambda _{\ast },\chi
,\varepsilon ,R_{0},k_{0}\right) $ and $\tau \in (0,1)$, then there exists a
constant $C>1$ depending only upon the data and not-dependent on $v$ and $%
x_{0}\in \Omega $ such that for every pair of balls $B_{\tau \varrho }\left(
x_{0}\right) \subset B_{\varrho }\left( x_{0}\right) \subset \subset \Omega $
with $0<\varrho <R_{0}$ 
\begin{equation}
\left\Vert v\right\Vert _{L^{\infty }\left( B_{\tau \varrho }\left(
x_{0}\right) \right) }\leq \max \left\{ \lambda _{\ast }\varrho ^{\frac{%
N\varepsilon }{p}};\frac{C}{\left( 1-\tau \right) ^{\frac{N}{p}}}\left[ 
\frac{1}{\left\vert B_{\varrho }\left( x_{0}\right) \right\vert }%
\int\limits_{B_{\varrho }\left( x_{0}\right) }\left\vert v\right\vert
^{p}\,dx\right] ^{\frac{1}{p}}\right\} .
\end{equation}
\end{theorem}

For more details on De Giorgi's classes and for the proof of the Theorem 6
refer to [13],[22], [24], [31] and [32].

\section{Proof of Theorm 2}

\subsection{The case of H.1.4}

Let us consider $y\in \Omega $ then we fix $R_{0}=\frac{1}{4}\min \left\{ 
\frac{1}{\sqrt[N]{\varpi _{N}}},dist\left( \partial \Omega ,y\right)
\right\} $, where\ $\varpi _{N}=\left\vert B_{1}\left( 0\right) \right\vert $%
, and we define $\Sigma =\left\{ x\in \Omega :\left\vert x-y\right\vert \leq
R_{0}\right\} $. We fix $x_{0}\in \Sigma $, $R_{1}=\frac{1}{4}dist\left(
\partial \Sigma ,x_{0}\right) $, $R_{0}<\min\limits_{\alpha =1,...,n}\left\{
R_{1},A_{\alpha },B_{\alpha }\right\} $, where

$A_{\alpha }=\frac{1}{\left[ 2\cdot 4^{\gamma _{\alpha }}L_{\alpha
}C_{N,p,\gamma _{\alpha }}D^{\frac{p^{\ast }-\gamma _{\alpha }}{p^{\ast }}%
}\left\Vert b_{\alpha }\right\Vert _{L^{\sigma _{\alpha }}\left( \Sigma
\right) }\left\Vert u^{\alpha }\right\Vert _{W^{1,p}\left( \Sigma \right)
}^{\gamma _{\alpha }-p}\varpi _{N}^{\epsilon _{\alpha }}\right] ^{\frac{1}{%
\epsilon _{\alpha }N}}}$, $B_{\alpha }=\frac{1}{\left[ 2D^{p^{\ast }-\gamma
_{\alpha }}\left\Vert b_{\alpha }\right\Vert _{L^{\sigma _{\alpha }}\left(
\Sigma \right) }\left\Vert u^{\alpha }\right\Vert _{W^{1,p}\left( \Sigma
\right) }^{\gamma _{\alpha }-p}\right] ^{\frac{1}{\epsilon _{\alpha }N}}}$, $%
C_{N,p,\gamma _{\alpha }}$ and $D^{\frac{p^{\ast }-\gamma _{\alpha }}{%
p^{\ast }}}$ are universal positive constants.

We fix $0<\varrho \leq t<s\leq R<R_{0}$, $B_{z}\left( x_{0}\right) =\left\{
x:\left\vert x-x_{0}\right\vert <z\right\} $, $k\in 
\mathbb{R}
$\ and we choose $\eta \in C_{c}^{\infty }\left( B_{s}\left( x_{0}\right)
\right) $\ such that $\eta =1$\ on $B_{t}\left( x_{0}\right) $, $0\leq \eta
\leq 1$\ on $B_{s}\left( x_{0}\right) $\ and $\left\vert \nabla \eta
\right\vert \leq \frac{2}{s-t}$\ on $B_{s}\left( x_{0}\right) $. Let us
define%
\begin{equation*}
\varphi =-\eta ^{p}w
\end{equation*}%
where $w\in W^{1,p}\left( \Sigma ,%
\mathbb{R}
^{n}\right) $ with%
\begin{equation*}
w^{1}=\max \left( u^{1}-k,0\right) ,w^{\alpha }=0,\alpha =2,...,n
\end{equation*}%
Let us observe that $\varphi =0$ $\mathcal{L}^{N}$-a.e. in $\Omega
\backslash \left( \left\{ \eta >0\right\} \cap \left\{ u^{1}>k\right\}
\right) $ thus%
\begin{equation}
\nabla u+\nabla \varphi =\nabla u
\end{equation}%
$\mathcal{L}^{N}$-a.e. in $\Omega \backslash \left( \left\{ \eta >0\right\}
\cap \left\{ u^{1}>k\right\} \right) $. Let us define 
\begin{equation}
A=\left( 
\begin{tabular}{l}
$p\eta ^{-1}\nabla \eta \left( k-u^{1}\right) $ \\ 
$\nabla u^{2}$ \\ 
$\vdots $ \\ 
$\nabla u^{n}$%
\end{tabular}%
\right)
\end{equation}%
since%
\begin{equation}
\nabla w=\left( 
\begin{tabular}{l}
$\nabla u^{1}$ \\ 
$0$ \\ 
$\vdots $ \\ 
$0$%
\end{tabular}%
\right)
\end{equation}%
$\mathcal{L}^{N}$-a.e. in $\Omega \backslash \left( \left\{ \eta >0\right\}
\cap \left\{ u^{1}>k\right\} \right) $ then we deduce that%
\begin{equation}
\nabla u+\nabla \varphi =\left( 1-\eta ^{p}\right) \nabla u+\eta ^{p}A
\end{equation}%
$\mathcal{L}^{N}$-a.e. in $\Omega \backslash \left( \left\{ \eta >0\right\}
\cap \left\{ u^{1}>k\right\} \right) $. Since $u$ is a local minimizer of
the functional (1.1) then we get%
\begin{equation}
J\left( u,\Sigma \right) \leq J\left( u+\varphi ,\Sigma \right)
\end{equation}%
it is%
\begin{equation}
\begin{tabular}{l}
$\int\limits_{\Sigma }\sum\limits_{\alpha =1}^{n}f_{\alpha }\left(
x,u^{\alpha },\nabla u^{\alpha }\right) +G\left( x,u,\nabla u\right) \,dx$
\\ 
$\leq \int\limits_{\Sigma }\sum\limits_{\alpha =1}^{n}f_{\alpha }\left(
x,u^{\alpha }+\varphi ^{\alpha },\nabla u^{\alpha }+\nabla \varphi ^{\alpha
}\right) +G\left( x,u+\varphi ,\nabla u+\nabla \varphi \right) \,dx$%
\end{tabular}%
\end{equation}%
and 
\begin{equation}
\begin{tabular}{l}
$\int\limits_{\Sigma }\sum\limits_{\alpha =2}^{n}f_{\alpha }\left(
x,u^{\alpha },\nabla u^{\alpha }\right) \,dx+\int\limits_{\Sigma
}f_{1}\left( x,u^{1},\nabla u^{1}\right) +G\left( x,u,\nabla u\right) \,dx$
\\ 
$\leq \int\limits_{\Sigma }\sum\limits_{\alpha =2}^{n}f_{\alpha }\left(
x,u^{\alpha },\nabla u^{\alpha }\right) \,dx+\int\limits_{\Sigma
}f_{1}\left( x,u^{1}+\varphi ^{1},\nabla u^{1}+\nabla \varphi ^{1}\right)
+G\left( x,u+\varphi ,\nabla u+\nabla \varphi \right) \,dx$%
\end{tabular}%
\end{equation}%
From (3.7) we deduce%
\begin{equation}
\begin{tabular}{l}
$\int\limits_{\Sigma }f_{1}\left( x,u^{1},\nabla u^{1}\right) +G\left(
x,u,\nabla u\right) \,dx$ \\ 
$\leq \int\limits_{\Sigma }f_{1}\left( x,u^{1}+\varphi ^{1},\nabla
u^{1}+\nabla \varphi ^{1}\right) +G\left( x,u+\varphi ,\nabla u+\nabla
\varphi \right) \,dx$ \\ 
$=\int\limits_{B_{r}\left( x_{0}\right) }f_{1}\left( x,u^{1}+\varphi
^{1},\nabla u^{1}+\nabla \varphi ^{1}\right) \,dx+\int\limits_{\Sigma
-B_{r}\left( x_{0}\right) \backslash }f_{1}\left( x,u^{1},\nabla
u^{1}\right) \,dx$ \\ 
$+\int\limits_{B_{s}\left( x_{0}\right) }G\left( x,u+\varphi ,\nabla
u+\nabla \varphi \right) \,dx+\int\limits_{\Sigma -B_{s}\left( x_{0}\right)
\backslash }G\left( x,u,\nabla u\right) \,dx$%
\end{tabular}%
\end{equation}%
and%
\begin{equation}
\begin{tabular}{l}
$\int\limits_{B_{s}\left( x_{0}\right) }f_{1}\left( x,u^{1},\nabla
u^{1}\right) +G\left( x,u,\nabla u\right) \,dx$ \\ 
$\leq \int\limits_{B_{s}\left( x_{0}\right) }f_{1}\left( x,u^{1}+\varphi
^{1},\nabla u^{1}+\nabla \varphi ^{1}\right) \,dx+$ \\ 
$+\int\limits_{B_{s}\left( x_{0}\right) }G\left( x,u+\varphi ,\nabla
u+\nabla \varphi \right) \,dx$%
\end{tabular}%
\end{equation}%
Let us define $E_{k,s}^{1}=\left\{ \eta >0\right\} \cap \left\{
u^{1}>k\right\} \cap B_{s}\left( x_{0}\right) \subset B_{s}\left(
x_{0}\right) $\ then%
\begin{equation}
\begin{tabular}{l}
$\int\limits_{E_{k,s}^{1}}f_{1}\left( x,u^{1},\nabla u^{1}\right) +G\left(
x,u,\nabla u\right) \,dx+\int\limits_{B_{s}\left( x_{0}\right)
-E_{k,s}^{1}}f_{1}\left( x,u^{1},\nabla u^{1}\right) +G\left( x,u,\nabla
u\right) \,dx$ \\ 
$\leq \int\limits_{E_{k,s}^{1}}f_{1}\left( x,u^{1}+\varphi ^{1},\nabla
u^{1}+\nabla \varphi ^{1}\right) \,dx+\int\limits_{B_{s}\left( x_{0}\right)
-E_{k,s}^{1}}f_{1}\left( x,u^{1},\nabla u^{1}\right) \,dx$ \\ 
$+\int\limits_{E_{k,s}^{1}}G\left( x,u+\varphi ,\nabla u+\nabla \varphi
\right) \,dx$ \\ 
$+\int\limits_{B_{s}\left( x_{0}\right) -E_{k,s}^{1}}G\left( x,u,\nabla
u\right) \,dx$%
\end{tabular}%
\end{equation}%
and%
\begin{equation}
\begin{tabular}{l}
$\int\limits_{E_{k,s}^{1}}f_{1}\left( x,u^{1},\nabla u^{1}\right) +G\left(
x,u,\nabla u\right) \,dx$ \\ 
$\leq \int\limits_{E_{k,s}^{1}}f_{1}\left( x,u^{1}+\varphi ^{1},\nabla
u^{1}+\nabla \varphi ^{1}\right) \,dx+$ \\ 
$+\int\limits_{E_{k,s}^{1}}G\left( x,u+\varphi ,\nabla u+\nabla \varphi
\right) \,dx$%
\end{tabular}
\label{45}
\end{equation}%
Since $\nabla u-A$ is a rank-one matrix, then by hypothesis H.1.3 it follows 
\begin{equation}
\begin{tabular}{l}
$\int\limits_{E_{k,s}^{1}}G\left( x,u+\varphi ,\nabla u+\nabla \varphi
\right) \,dx$ \\ 
$=\int\limits_{E_{k,s}^{1}}G\left( x,u+\varphi ,\left( 1-\eta ^{p}\right)
\nabla u+\eta ^{p}A\right) \,dx$ \\ 
$\leq \int\limits_{E_{k,s}^{1}}\left( 1-\eta ^{p}\right) G\left( x,u+\varphi
,\nabla u\right) +\eta ^{p}G\left( x,u+\varphi ,A\right) \,dx$ \\ 
$=\int\limits_{E_{k,s}^{1}}G\left( x,u+\varphi ,\nabla u\right) +\eta ^{p}%
\left[ G\left( x,u+\varphi ,A\right) -G\left( x,u+\varphi ,\nabla u\right) %
\right] \,dx$%
\end{tabular}
\label{46}
\end{equation}%
using (\ref{45}) and (\ref{46}) we get%
\begin{equation}
\begin{tabular}{l}
$\int\limits_{E_{k,s}^{1}}f_{1}\left( x,u^{1},\nabla u^{1}\right) \,dx$ \\ 
$\leq \int\limits_{E_{k,s}^{1}}f_{1}\left( x,u^{1}+\varphi ^{1},\nabla
u^{1}+\nabla \varphi ^{1}\right) \,dx+$ \\ 
$+\int\limits_{E_{k,s}^{1}}\eta ^{p}\left[ G\left( x,u+\varphi ,A\right)
-G\left( x,u+\varphi ,\nabla u\right) \right] \,dx$ \\ 
$+\int\limits_{E_{k,s}^{1}}G\left( x,u+\varphi ,\nabla u\right) -G\left(
x,u,\nabla u\right) \,dx$%
\end{tabular}%
\end{equation}%
Using hypothesis H.1.1 we have%
\begin{equation}
\begin{tabular}{l}
$\int\limits_{E_{k,s}^{1}}\left\vert \nabla u^{1}\right\vert
^{p}-b_{1}\left( x\right) \left\vert u^{1}\right\vert ^{\gamma
_{1}}-a_{1}\left( x\right) \,dx$ \\ 
$\leq L_{1}\int\limits_{E_{k,s}^{1}}\left\vert \nabla u^{1}+\nabla \varphi
^{1}\right\vert ^{p}+b_{1}\left( x\right) \left\vert u^{1}+\varphi
^{1}\right\vert ^{\gamma _{1}}+a_{1}\left( x\right) \,dx$ \\ 
$+\int\limits_{E_{k,s}^{1}}\eta ^{p}\left\vert G\left( x,u+\varphi ,A\right)
-G\left( x,u+\varphi ,\nabla u\right) \right\vert \,dx$ \\ 
$+\int\limits_{E_{k,s}^{1}}\left\vert G\left( x,u+\varphi ,\nabla u\right)
-G\left( x,u,\nabla u\right) \right\vert \,dx$%
\end{tabular}%
\end{equation}%
then%
\begin{equation}
\begin{tabular}{l}
$\int\limits_{E_{k,s}^{1}}\left\vert \nabla u^{1}\right\vert ^{p}\,dx$ \\ 
$\leq \int\limits_{E_{k,s}^{1}}b_{1}\left( x\right) \left\vert
u^{1}\right\vert ^{\gamma _{1}}+a_{1}\left( x\right) \,dx$ \\ 
$+L_{1}\int\limits_{E_{k,s}^{1}}\left\vert \nabla u^{1}+\nabla \varphi
^{1}\right\vert ^{p}+b_{1}\left( x\right) \left\vert u^{1}+\varphi
^{1}\right\vert ^{\gamma _{1}}+a_{1}\left( x\right) \,dx$ \\ 
$+\int\limits_{E_{k,s}^{1}}\eta ^{p}\left\vert G\left( x,u+\varphi ,A\right)
-G\left( x,u+\varphi ,\nabla u\right) \right\vert \,dx$ \\ 
$+\int\limits_{E_{k,s}^{1}}\left\vert G\left( x,u+\varphi ,\nabla u\right)
-G\left( x,u,\nabla u\right) \right\vert \,dx$%
\end{tabular}
\label{49}
\end{equation}%
For every $x\in E_{k,s}^{1}$ it follows 
\begin{equation}
u^{1}=\left( 1-\eta ^{p}\right) u^{1}+\eta ^{p}\left( w^{1}+k\right)
\label{50}
\end{equation}%
and%
\begin{equation}
u^{1}+\varphi ^{1}=u^{1}\left( 1-\eta ^{p}\right) +k\eta ^{p}  \label{51}
\end{equation}%
Using (\ref{50}) and (\ref{51}) for every $x\in E_{k,s}^{1}$ we get%
\begin{equation}
\nabla u^{1}+\nabla \varphi ^{1}=\nabla u^{1}\left( 1-\eta ^{p}\right)
-pw^{1}\eta ^{p-1}\nabla \eta  \label{52}
\end{equation}%
and%
\begin{equation}
\left\vert \nabla u^{1}+\nabla \varphi ^{1}\right\vert ^{p}\leq
2^{p-1}\left( 1-\eta ^{p}\right) \left\vert \nabla u^{1}\right\vert
^{p}+2^{p-1}p^{p}\eta ^{p-1}\left( u^{1}-k\right) ^{p}\left\vert \nabla \eta
\right\vert ^{p}  \label{53}
\end{equation}%
Recalling the properties of the function $\eta $ and the inequalities (\ref%
{52}) and (\ref{53}) it follows that%
\begin{equation}
\begin{tabular}{l}
$\int\limits_{E_{k,s}^{1}}\left\vert \nabla u^{1}+\nabla \varphi
^{1}\right\vert ^{p}\,dx$ \\ 
$\leq 2^{p-1}\int\limits_{E_{k,s}^{1}}\left( 1-\eta ^{p}\right) \left\vert
\nabla u^{1}\right\vert ^{p}\,dx$ \\ 
$+2^{2p-1}p^{p}\int\limits_{E_{k,s}^{1}-E_{k,t}^{1}}\frac{\left(
u^{1}-k\right) ^{p}}{\left( s-t\right) ^{p}}\,dx$%
\end{tabular}
\label{54}
\end{equation}%
Moreover, since on the set $E_{k,s}^{1}$ we have%
\begin{equation}
\begin{tabular}{l}
$\left\vert u^{1}\right\vert ^{\gamma _{1}}=\left( \left( 1-\eta ^{p}\right)
u^{1}+\eta ^{p}\left( w^{1}+k\right) \right) ^{\gamma _{1}}$ \\ 
$\leq 2^{\gamma _{1}-1}\left( 1-\eta ^{p}\right) ^{\gamma _{1}}\left(
u^{1}\right) ^{\gamma _{1}}+4^{\gamma _{1}-1}\eta ^{p\gamma _{1}}\left(
w^{1}\right) ^{\gamma _{1}}+4^{\gamma _{1}-1}\eta ^{p\gamma _{1}}\left\vert
k\right\vert ^{\gamma _{1}}$%
\end{tabular}%
\end{equation}%
and 
\begin{equation}
\begin{tabular}{l}
$\left( u^{1}\left( 1-\eta ^{p}\right) +k\eta ^{p}\right) ^{\gamma _{1}}$ \\ 
$\leq 2^{\gamma _{1}-1}\left( 1-\eta ^{p}\right) ^{\gamma _{1}}\left(
u^{1}\right) ^{\gamma _{1}}+2^{\gamma _{1}-1}\eta ^{p\gamma _{1}}\left\vert
k\right\vert ^{\gamma _{1}}$%
\end{tabular}%
\end{equation}%
then we get%
\begin{equation}
\begin{tabular}{l}
$\int\limits_{E_{k,s}^{1}}b_{1}\left( x\right) \left\vert u^{1}\right\vert
^{\gamma _{1}}+L_{1}b_{1}\left( x\right) \left\vert u^{1}+\varphi
^{1}\right\vert ^{\gamma _{1}}\,dx$ \\ 
$\leq L_{1}\int\limits_{E_{k,s}^{1}}b_{1}\left( x\right) \left( \left\vert
u^{1}\right\vert ^{\gamma _{1}}+\left\vert u^{1}+\varphi ^{1}\right\vert
^{\gamma _{1}}\right) \,dx$ \\ 
$\leq 2\,4^{\gamma _{1}-1}L_{1}\int\limits_{E_{k,s}^{1}}b_{1}\left( x\right)
\left( \left( 1-\eta ^{p}\right) ^{\gamma _{1}}\left( u^{1}\right) ^{\gamma
_{1}}+\eta ^{p\gamma _{1}}\left( w^{1}\right) ^{\gamma _{1}}+\eta ^{p\gamma
_{1}}\left\vert k\right\vert ^{\gamma _{1}}\right) \,dx$%
\end{tabular}
\label{57}
\end{equation}%
Using (\ref{49}), (\ref{54}) and (\ref{57}) it follows%
\begin{equation}
\begin{tabular}{l}
$\int\limits_{E_{k,s}^{1}}\left\vert \nabla u^{1}\right\vert ^{p}\,dx$ \\ 
$\leq 2^{p-1}L_{1}\int\limits_{E_{k,s}^{1}}\left( 1-\eta ^{p}\right)
\left\vert \nabla u^{1}\right\vert
^{p}\,dx+2^{2p-1}p^{p}L_{1}\int\limits_{E_{k,s}^{1}-E_{k,t}^{1}}\frac{\left(
u^{1}-k\right) ^{p}}{\left( s-t\right) ^{p}}\,dx+2L_{1}\int%
\limits_{E_{k,s}^{1}}a_{1}\left( x\right) \,dx$ \\ 
$+2\,4^{\gamma _{1}-1}L_{1}\int\limits_{E_{k,s}^{1}}b_{1}\left( x\right)
\left( \left( 1-\eta ^{p}\right) ^{\gamma _{1}}\left( u^{1}\right) ^{\gamma
_{1}}+\eta ^{p\gamma _{1}}\left( w^{1}\right) ^{\gamma _{1}}+\eta ^{p\gamma
_{1}}\left\vert k\right\vert ^{\gamma _{1}}\right) \,dx$ \\ 
$+\int\limits_{E_{k,s}^{1}}\eta ^{p}\left\vert G\left( x,u+\varphi ,A\right)
-G\left( x,u+\varphi ,\nabla u\right) \right\vert \,dx$ \\ 
$+\int\limits_{E_{k,s}^{1}}\left\vert G\left( x,u+\varphi ,\nabla u\right)
-G\left( x,u,\nabla u\right) \right\vert \,dx$%
\end{tabular}%
\end{equation}%
then, adding the term $\int\limits_{E_{k,s}^{1}}b_{1}\left( x\right) \left(
u^{1}\right) ^{\gamma _{1}}\,dx$\ to both sides of the inequality we obtain%
\begin{equation}
\begin{tabular}{l}
$\int\limits_{E_{k,s}^{1}}\left\vert \nabla u^{1}\right\vert
^{p}+b_{1}\left( x\right) \left( u^{1}\right) ^{\gamma _{1}}\,dx$ \\ 
$\leq c\left( p,\gamma \right) \int\limits_{E_{k,s}^{1}}\left( 1-\eta
^{p}\right) \left( \left\vert \nabla u^{1}\right\vert ^{p}+b_{1}\left(
x\right) \left( u^{1}\right) ^{\gamma _{1}}\right)
\,dx+\int\limits_{E_{k,s}^{1}}b_{1}\left( x\right) \left( u^{1}\right)
^{\gamma _{1}}\,dx$ \\ 
$+2^{2p-1}p^{p}\int\limits_{E_{k,s}^{1}-E_{k,t}^{1}}\frac{\left(
u^{1}-k\right) ^{p}}{\left( s-t\right) ^{p}}\,dx+2L_{1}\int%
\limits_{E_{k,s}^{1}}a_{1}\left( x\right) \,dx$ \\ 
$+2\,4^{\gamma _{1}-1}L_{1}\int\limits_{E_{k,s}^{1}}b_{1}\left( x\right)
\left( \eta ^{p\gamma _{1}}\left( w^{1}\right) ^{\gamma _{1}}+\eta ^{p\gamma
_{1}}\left\vert k\right\vert ^{\gamma _{1}}\right) \,dx$ \\ 
$+\int\limits_{E_{k,s}^{1}}\eta ^{p}\left\vert G\left( x,u+\varphi ,A\right)
-G\left( x,u+\varphi ,\nabla u\right) \right\vert \,dx$ \\ 
$+\int\limits_{E_{k,s}^{1}}\left\vert G\left( x,u+\varphi ,\nabla u\right)
-G\left( x,u,\nabla u\right) \right\vert \,dx$%
\end{tabular}%
\end{equation}%
where $c\left( p,\gamma \right) =L_{1}(2^{p-1},2\,4^{\gamma _{1}-1})$, \
remembering that%
\begin{equation}
\int\limits_{E_{k,s}^{1}}b_{1}\left( x\right) \left( u^{1}\right) ^{\gamma
_{1}}\,dx\leq 2\,4^{\gamma _{1}-1}L_{1}\int\limits_{E_{k,s}^{1}}b_{1}\left(
x\right) \left( \left( 1-\eta ^{p}\right) ^{\gamma _{1}}\left( u^{1}\right)
^{\gamma _{1}}+\eta ^{p\gamma _{1}}\left( w^{1}\right) ^{\gamma _{1}}+\eta
^{p\gamma _{1}}\left\vert k\right\vert ^{\gamma _{1}}\right) \,dx
\end{equation}%
it follows%
\begin{equation}
\begin{tabular}{l}
$\int\limits_{E_{k,s}^{1}}\left\vert \nabla u^{1}\right\vert
^{p}+b_{1}\left( x\right) \left( u^{1}\right) ^{\gamma _{1}}\,dx$ \\ 
$\leq 2c\left( p,\gamma \right) \int\limits_{E_{k,s}^{1}}\left( 1-\eta
^{p}\right) \left( \left\vert \nabla u^{1}\right\vert ^{p}+b_{1}\left(
x\right) \left( u^{1}\right) ^{\gamma _{1}}\right) \,dx$ \\ 
$+2^{2p-1}p^{p}\int\limits_{E_{k,s}^{1}-E_{k,t}^{1}}\frac{\left(
u^{1}-k\right) ^{p}}{\left( s-t\right) ^{p}}\,dx+2L_{1}\int%
\limits_{E_{k,s}^{1}}a_{1}\left( x\right) \,dx$ \\ 
$+4^{\gamma _{1}}L_{1}\int\limits_{E_{k,s}^{1}}b_{1}\left( x\right) \eta
^{p\gamma _{1}}\left( w^{1}\right) ^{\gamma _{1}}\,dx+4^{\gamma
_{1}}L_{1}\int\limits_{E_{k,s}^{1}}b_{1}\left( x\right) \eta ^{p\gamma
_{1}}\left\vert k\right\vert ^{\gamma _{1}}\,dx$ \\ 
$+\int\limits_{E_{k,s}^{1}}\eta ^{p}\left\vert G\left( x,u+\varphi ,A\right)
-G\left( x,u+\varphi ,\nabla u\right) \right\vert \,dx$ \\ 
$+\int\limits_{E_{k,s}^{1}}\left\vert G\left( x,u+\varphi ,\nabla u\right)
-G\left( x,u,\nabla u\right) \right\vert \,dx$%
\end{tabular}
\label{61}
\end{equation}%
Now we have to estimate the term $\int\limits_{E_{k,s}^{1}}b_{1}\left(
x\right) \eta ^{p\gamma _{1}}\left( w^{1}\right) ^{\gamma _{1}}\,dx$, using
the H\"{o}lder'Inequality we have%
\begin{equation}
\begin{tabular}{l}
$\int\limits_{E_{k,s}^{1}}b_{1}\left( x\right) \eta ^{\left( p-1\right)
\gamma _{1}}\left( \eta w^{1}\right) ^{\gamma _{1}}\,dx$ \\ 
$\leq \left[ \int\limits_{E_{k,s}^{1}}\left( b_{1}\left( x\right) \eta
^{\left( p-1\right) \gamma _{1}}\right) ^{\frac{p^{\ast }}{p^{\ast }-\gamma
_{1}}}\,dx\right] ^{\frac{p^{\ast }-\gamma _{1}}{p^{\ast }}}\left[
\int\limits_{E_{k,s}^{1}}\left( \eta w^{1}\right) ^{p^{\ast }}\,dx\right] ^{%
\frac{\gamma _{1}}{p^{\ast }}}$ \\ 
$\leq \left[ \left( \int\limits_{E_{k,s}^{1}}\left( b_{1}\left( x\right)
\right) ^{\sigma _{1}}\,dx\right) ^{\frac{p^{\ast }}{\left( p^{\ast }-\gamma
_{1}\right) \sigma _{1}}}\left( \int\limits_{E_{k,s}^{1}}\eta ^{\frac{%
p^{\ast }\left( p-1\right) \gamma _{1}\sigma _{1}}{\left( p^{\ast }-\gamma
_{1}\right) \sigma _{1}-p^{\ast }}}\,dx\right) ^{\frac{\left( p^{\ast
}-\gamma _{1}\right) \sigma _{1}-p^{\ast }}{\left( p^{\ast }-\gamma
_{1}\right) \sigma _{1}}}\right] ^{\frac{p^{\ast }-\gamma _{1}}{p^{\ast }}}%
\left[ \int\limits_{E_{k,s}^{1}}\left( \eta w^{1}\right) ^{p^{\ast }}\,dx%
\right] ^{\frac{\gamma _{1}}{p^{\ast }}}$ \\ 
$\leq \left[ \left\Vert b_{1}\right\Vert _{L^{\sigma _{1}}\left( B_{s}\left(
x_{0}\right) \right) }^{\frac{p^{\ast }}{\left( p^{\ast }-\gamma _{1}\right) 
}}\cdot \left\vert E_{k,s}^{1}\right\vert ^{\frac{\left( p^{\ast }-\gamma
_{1}\right) \sigma _{1}-p^{\ast }}{\left( p^{\ast }-\gamma _{1}\right)
\sigma _{1}}}\right] ^{\frac{p^{\ast }-\gamma _{1}}{p^{\ast }}}\left[
\int\limits_{B_{s}\left( x_{0}\right) }\left( \eta w^{1}\right) ^{p^{\ast
}}\,dx\right] ^{\frac{\gamma _{1}}{p^{\ast }}}$%
\end{tabular}%
\end{equation}%
Noting that $\frac{\gamma _{1}}{p^{\ast }}=\frac{\gamma _{1}-p}{p^{\ast }}+%
\frac{p}{p^{\ast }}$, it follows%
\begin{equation}
\begin{tabular}{l}
$\int\limits_{E_{k,s}^{1}}b_{1}\left( x\right) \eta ^{\left( p-1\right)
\gamma _{1}}\left( \eta w^{1}\right) ^{\gamma _{1}}\,dx$ \\ 
$\leq \left[ \left\Vert b_{1}\right\Vert _{L^{\sigma _{1}}\left( B_{s}\left(
x_{0}\right) \right) }^{\frac{p^{\ast }}{\left( p^{\ast }-\gamma _{1}\right) 
}}\cdot \left\vert E_{k,s}^{1}\right\vert ^{\frac{\left( p^{\ast }-\gamma
_{1}\right) \sigma _{1}-p^{\ast }}{\left( p^{\ast }-\gamma _{1}\right)
\sigma _{1}}}\right] ^{\frac{p^{\ast }-\gamma _{1}}{p^{\ast }}}\left[
\int\limits_{B_{s}\left( x_{0}\right) }\left( \eta w^{1}\right) ^{p^{\ast
}}\,dx\right] ^{\frac{\gamma _{1}-1}{p^{\ast }}}\left[ \int\limits_{B_{s}%
\left( x_{0}\right) }\left( \eta w^{1}\right) ^{p^{\ast }}\,dx\right] ^{%
\frac{1}{p^{\ast }}}$%
\end{tabular}%
\end{equation}%
and 
\begin{equation}
\begin{tabular}{l}
$\int\limits_{E_{k,s}^{1}}b_{1}\left( x\right) \eta ^{\left( p-1\right)
\gamma _{1}}\left( \eta w^{1}\right) ^{\gamma _{1}}\,dx$ \\ 
$\leq \left[ \left\Vert b_{1}\right\Vert _{L^{\sigma _{1}}\left( B_{s}\left(
x_{0}\right) \right) }^{\frac{p^{\ast }}{\left( p^{\ast }-\gamma _{1}\right) 
}}\cdot \left\vert E_{k,s}^{1}\right\vert ^{\frac{\left( p^{\ast }-\gamma
_{1}\right) \sigma _{1}-p^{\ast }}{\left( p^{\ast }-\gamma _{1}\right)
\sigma _{1}}}\left( \int\limits_{B_{s}\left( x_{0}\right) }\left(
u^{1}\right) ^{p^{\ast }}\,dx\right) ^{\frac{\gamma _{1}-p}{p^{\ast }-\gamma
_{1}}}\right] ^{\frac{p^{\ast }-\gamma _{1}}{p^{\ast }}}\left[
\int\limits_{B_{s}\left( x_{0}\right) }\left( \eta w^{1}\right) ^{p^{\ast
}}\,dx\right] ^{\frac{p}{p^{\ast }}}$%
\end{tabular}
\label{4.30}
\end{equation}%
Applying Sobolev's Immersion Theorem there exists a constant $C_{IS}>0$,
dependent only on $p$ and $N$, such that 
\begin{equation}
\left[ \int\limits_{B_{s}\left( x_{0}\right) }\left( u^{1}\right) ^{p^{\ast
}}\,dx\right] ^{\frac{1}{p^{\ast }}}\leq C_{IS}\left\Vert u^{1}\right\Vert
_{W^{1,p}\left( B_{s}\left( x_{0}\right) \right) }  \label{4.31}
\end{equation}%
then, by (\ref{4.30}) and (\ref{4.31}) it follows%
\begin{equation}
\begin{tabular}{l}
$\int\limits_{E_{k,s}^{1}}b_{1}\left( x\right) \eta ^{\left( p-1\right)
\gamma _{1}}\left( \eta w^{1}\right) ^{\gamma _{1}}\,dx$ \\ 
$\leq \left[ D\left\Vert b_{1}\right\Vert _{L^{\sigma _{1}}\left(
B_{s}\left( x_{0}\right) \right) }^{\frac{p^{\ast }}{\left( p^{\ast }-\gamma
_{1}\right) }}\cdot \left\vert E_{k,s}^{1}\right\vert ^{\frac{\left( p^{\ast
}-\gamma _{1}\right) \sigma _{1}-p^{\ast }}{\left( p^{\ast }-\gamma
_{1}\right) \sigma _{1}}}\left( \left\Vert u^{1}\right\Vert _{W^{1,p}\left(
B_{s}\left( x_{0}\right) \right) }\right) ^{\frac{\left( \gamma
_{1}-p\right) p^{\ast }}{p^{\ast }-\gamma _{1}}}\right] ^{\frac{p^{\ast
}-\gamma _{1}}{p^{\ast }}}\left[ \int\limits_{B_{s}\left( x_{0}\right)
}\left( \eta w^{1}\right) ^{p^{\ast }}\,dx\right] ^{\frac{p}{p^{\ast }}}$%
\end{tabular}
\label{66}
\end{equation}%
where $D=C_{IS}^{\frac{\gamma _{1}-p}{p^{\ast }-\gamma _{1}}}$.

Applying Sobolev's inequality we have%
\begin{equation}
\left[ \int\limits_{B_{s}\left( x_{0}\right) }\left( \eta w^{1}\right)
^{p^{\ast }}\,dx\right] ^{\frac{p}{p^{\ast }}}\leq C_{N,p,\gamma _{1}}\left[
\int\limits_{B_{s}\left( x_{0}\right) }\eta ^{p}\left\vert \nabla
w^{1}\right\vert ^{p}+\frac{2^{p}\left( w^{1}\right) ^{p}}{\left( s-t\right)
^{p}}\,dx\right]  \label{67}
\end{equation}%
then, using (\ref{66}) and (\ref{67}), it follows%
\begin{equation}
\begin{tabular}{l}
$\int\limits_{E_{k,s}^{1}}b_{1}\left( x\right) \eta ^{\left( p-1\right)
\gamma _{1}}\left( \eta w^{1}\right) ^{\gamma _{1}}\,dx$ \\ 
$\leq C_{N,p,\gamma _{1}}\left[ D\left\Vert b_{1}\right\Vert _{L^{\sigma
_{1}}\left( B_{s}\left( x_{0}\right) \right) }^{\frac{p^{\ast }}{\left(
p^{\ast }-\gamma _{1}\right) }}\cdot \left\vert E_{k,s}^{1}\right\vert ^{%
\frac{\left( p^{\ast }-\gamma _{1}\right) \sigma _{1}-p^{\ast }}{\left(
p^{\ast }-\gamma _{1}\right) \sigma _{1}}}\left( \left\Vert u^{1}\right\Vert
_{W^{1,p}\left( B_{s}\left( x_{0}\right) \right) }\right) ^{\frac{\left(
\gamma _{1}-p\right) p^{\ast }}{p^{\ast }-\gamma _{1}}}\right] ^{\frac{%
p^{\ast }-\gamma _{1}}{p^{\ast }}}$ \\ 
$\cdot \left[ \int\limits_{E_{k,s}^{1}}\eta ^{p}\left\vert \nabla
w^{1}\right\vert ^{p}+\frac{2^{p}\left( u^{1}-k\right) ^{p}}{\left(
s-t\right) ^{p}}\,dx\right] $%
\end{tabular}%
\end{equation}%
Since $\frac{\left( p^{\ast }-\gamma _{1}\right) \sigma _{1}-p^{\ast }}{%
p^{\ast }\sigma _{1}}=\varepsilon $ and $\left\vert E_{k,s}^{1}\right\vert
\leq \varpi _{N}R_{0}^{N}$ we get 
\begin{equation}
\begin{tabular}{l}
$\int\limits_{E_{k,s}^{1}}b_{1}\left( x\right) \eta ^{\left( p-1\right)
\gamma _{1}}\left( \eta w^{1}\right) ^{\gamma _{1}}\,dx$ \\ 
$\leq C_{N,p,\gamma _{1}}D^{\frac{p^{\ast }-\gamma _{1}}{p^{\ast }}%
}\left\Vert b_{1}\right\Vert _{L^{\sigma _{1}}\left( \Sigma \right)
}\left\Vert u^{1}\right\Vert _{W^{1,p}\left( \Sigma \right) }^{\left( \gamma
_{1}-p\right) }\varpi _{N}^{\varepsilon }R_{0}^{\varepsilon N}$ \\ 
$\cdot \left[ \int\limits_{E_{k,s}^{1}}\eta ^{p}\left\vert \nabla
w^{1}\right\vert ^{p}+\frac{2^{p}\left( u^{1}-k\right) ^{p}}{\left(
s-t\right) ^{p}}\,dx\right] $%
\end{tabular}
\label{69}
\end{equation}%
now, using (\ref{61}) and (\ref{69}) it follows%
\begin{equation}
\begin{tabular}{l}
$\int\limits_{E_{k,s}^{1}}\left\vert \nabla u^{1}\right\vert
^{p}+b_{1}\left( x\right) \left( u^{1}\right) ^{\gamma _{1}}\,dx$ \\ 
$\leq 2c\left( p,\gamma \right) \int\limits_{E_{k,s}^{1}}\left( 1-\eta
^{p}\right) \left( \left\vert \nabla u^{1}\right\vert ^{p}+b_{1}\left(
x\right) \left( u^{1}\right) ^{\gamma _{1}}\right) \,dx$ \\ 
$+2^{2p-1}p^{p}\int\limits_{E_{k,s}^{1}-E_{k,t}^{1}}\frac{\left(
u^{1}-k\right) ^{p}}{\left( s-t\right) ^{p}}\,dx+2L_{1}\int%
\limits_{E_{k,s}^{1}}a_{1}\left( x\right) \,dx+4^{\gamma
_{1}}L_{1}\varepsilon \int\limits_{E_{k,s}^{1}}\eta ^{p}\left\vert \nabla
u^{1}\right\vert ^{p}+\frac{2^{p}\left( u^{1}-k\right) ^{p}}{\left(
s-t\right) ^{p}}\,dx$ \\ 
$+4^{\gamma _{1}}L_{1}C_{N,p,\gamma _{1}}D^{\frac{p^{\ast }-\gamma _{1}}{%
p^{\ast }}}\left\Vert b_{1}\right\Vert _{L^{\sigma _{1}}\left( \Sigma
\right) }\left\Vert u^{1}\right\Vert _{W^{1,p}\left( \Sigma \right)
}^{\left( \gamma _{1}-p\right) }\varpi _{N}^{\varepsilon }R_{0}^{\varepsilon
N}\cdot \left[ \int\limits_{E_{k,s}^{1}}\eta ^{p}\left\vert \nabla
w^{1}\right\vert ^{p}+\frac{2^{p}\left( u^{1}-k\right) ^{p}}{\left(
s-t\right) ^{p}}\,dx\right] $ \\ 
$+4^{\gamma _{1}}L_{1}\int\limits_{E_{k,s}^{1}}b_{1}\left( x\right) \eta
^{p\gamma _{1}}\left\vert k\right\vert ^{\gamma
_{1}}\,dx+\int\limits_{E_{k,s}^{1}}\eta ^{p}\left\vert G\left( x,u+\varphi
,A\right) -G\left( x,u+\varphi ,\nabla u\right) \right\vert \,dx$ \\ 
$+\int\limits_{E_{k,s}^{1}}\left\vert G\left( x,u+\varphi ,\nabla u\right)
-G\left( x,u,\nabla u\right) \right\vert \,dx$%
\end{tabular}%
\end{equation}%
Recalling that we have fixed $R_{0}$ such that $C_{N,p,\gamma _{1}}D^{\frac{%
p^{\ast }-\gamma _{1}}{p^{\ast }}}\left\Vert b_{1}\right\Vert _{L^{\sigma
_{1}}\left( \Sigma \right) }\left\Vert u^{1}\right\Vert _{W^{1,p}\left(
\Sigma \right) }^{\left( \gamma _{1}-p\right) }\varpi _{N}^{\varepsilon
}R_{0}^{\varepsilon N}\leq \frac{1}{2\cdot 4^{\gamma _{1}}L_{1}}$\ then we
have%
\begin{equation}
\begin{tabular}{l}
$\frac{1}{2}\int\limits_{E_{k,s}^{1}}\left\vert \nabla u^{1}\right\vert
^{p}+b_{1}\left( x\right) \left( u^{1}\right) ^{\gamma _{1}}\,dx$ \\ 
$\leq 2c\left( p,\gamma \right) \int\limits_{E_{k,s}^{1}}\left( 1-\eta
^{p}\right) \left( \left\vert \nabla u^{1}\right\vert ^{p}+b_{1}\left(
x\right) \left( u^{1}\right) ^{\gamma _{1}}\right) \,dx$ \\ 
$+\left( 2^{2p-1}p^{p}+2^{p-1}\right) \int\limits_{E_{k,s}^{1}}\frac{\left(
u^{1}-k\right) ^{p}}{\left( s-t\right) ^{p}}\,dx+2L_{1}\int%
\limits_{E_{k,s}^{1}}a_{1}\left( x\right) \,dx$ \\ 
$+4^{\gamma _{1}}L_{1}\int\limits_{E_{k,s}^{1}}b_{1}\left( x\right) \eta
^{p\gamma _{1}}\left\vert k\right\vert ^{\gamma _{1}}\,dx$ \\ 
$+\int\limits_{E_{k,s}^{1}}\eta ^{p}\left\vert G\left( x,u+\varphi ,A\right)
-G\left( x,u+\varphi ,\nabla u\right) \right\vert \,dx$ \\ 
$+\int\limits_{E_{k,s}^{1}}\left\vert G\left( x,u+\varphi ,\nabla u\right)
-G\left( x,u,\nabla u\right) \right\vert \,dx$%
\end{tabular}
\label{71}
\end{equation}%
Now, let's consider the integral $\int\limits_{E_{k,s}^{1}}b_{1}\left(
x\right) \eta ^{p\gamma _{1}}\left\vert k\right\vert ^{\gamma _{1}}\,dx$,
then using the H\"{o}lder inequality we have%
\begin{equation}
\begin{tabular}{l}
$\int\limits_{E_{k,s}^{1}}b_{1}\left( x\right) \eta ^{p\gamma
_{1}}\left\vert k\right\vert ^{\gamma _{1}}\,dx$ \\ 
$\leq \left\vert k\right\vert ^{\gamma _{1}}\left\Vert b_{1}\right\Vert
_{L^{\sigma _{1}}\left( \Sigma \right) }\left\vert E_{k,s}^{1}\right\vert
^{1-\frac{1}{\sigma _{1}}}$%
\end{tabular}%
\end{equation}%
and%
\begin{equation}
\begin{tabular}{l}
$\int\limits_{E_{k,s}^{1}}b_{1}\left( x\right) \eta ^{p\gamma
_{1}}\left\vert k\right\vert ^{\gamma _{1}}\,dx$ \\ 
$\leq k^{\gamma _{1}}\left\Vert b_{1}\right\Vert _{L^{\sigma _{1}}\left(
\Sigma \right) }\left\vert A_{k,s}^{1}\right\vert ^{1-\frac{1}{\sigma _{1}}}$%
\end{tabular}%
\end{equation}%
Since $1-\frac{1}{\sigma _{1}}=\frac{\gamma _{1}}{p^{\ast }}+\varepsilon $
and $\frac{p}{p^{\ast }}=1-\frac{p}{N}$ then%
\begin{eqnarray*}
\left\vert k\right\vert ^{\gamma _{1}}\left\Vert b_{1}\right\Vert
_{L^{\sigma _{1}}\left( \Sigma \right) }\left\vert A_{k,s}^{1}\right\vert
^{1-\frac{1}{\sigma _{1}}} &=&\left\vert k\right\vert ^{\gamma
_{1}-p}\left\Vert b_{1}\right\Vert _{L^{\sigma _{1}}\left( \Sigma \right)
}\left\vert k\right\vert ^{p}\left\vert A_{k,s}^{1}\right\vert ^{1-\frac{1}{%
\sigma _{1}}} \\
&=&\left\Vert b_{1}\right\Vert _{L^{\sigma _{1}}\left( \Sigma \right)
}\left( \left\vert k\right\vert ^{p^{\ast }}\left\vert
A_{k,s}^{1}\right\vert \right) ^{\frac{\gamma _{1}-p}{p^{\ast }}}\left\vert
k\right\vert ^{p}\left\vert A_{k,s}^{1}\right\vert ^{1-\frac{p}{N}%
+\varepsilon }
\end{eqnarray*}%
\ Moreover it follows 
\begin{equation}
\begin{tabular}{l}
$\int\limits_{E_{k,s}^{1}}b_{1}\left( x\right) \eta ^{p\gamma
_{1}}\left\vert k\right\vert ^{\gamma _{1}}\,dx$ \\ 
$\leq \left\Vert b_{1}\right\Vert _{L^{\sigma _{1}}\left( \Sigma \right) }%
\left[ \int\limits_{B_{s}\left( x_{0}\right) }\left\vert u^{1}\right\vert
^{p^{\ast }}\,dx\right] ^{\frac{\gamma _{1}-p}{p^{\ast }}}\left\vert
k\right\vert ^{p}\left\vert A_{k,s}^{1}\right\vert ^{1-\frac{p}{N}%
+\varepsilon }$%
\end{tabular}%
\end{equation}%
using the Sobolev Embedding Theorem we get%
\begin{equation}
\begin{tabular}{l}
$\int\limits_{E_{k,s}^{1}}b_{1}\left( x\right) \eta ^{p\gamma _{1}}k^{\gamma
_{1}}\,dx$ \\ 
$\leq \left\Vert b_{1}\right\Vert _{L^{\sigma _{1}}\left( \Sigma \right) }%
\left[ C_{IS}\left\Vert u^{1}\right\Vert _{W^{1,p}\left( B_{s}\left(
x_{0}\right) \right) }\right] ^{\gamma _{1}-p}\left\vert k\right\vert
^{p}\left\vert A_{k,s}^{1}\right\vert ^{1-\frac{p}{N}+\varepsilon }$%
\end{tabular}%
\end{equation}%
Since%
\begin{eqnarray*}
&&\left\Vert b_{1}\right\Vert _{L^{\sigma _{1}}\left( \Sigma \right) }\left[
C_{IS}\left\Vert u^{1}\right\Vert _{W^{1,p}\left( B_{s}\left( x_{0}\right)
\right) }\right] ^{\gamma _{1}-p}\left\vert k\right\vert ^{p}\left\vert
A_{k,s}^{1}\right\vert ^{1-\frac{p}{N}+\varepsilon } \\
&\leq &\left\Vert b_{1}\right\Vert _{L^{\sigma _{1}}\left( \Sigma \right) } 
\left[ C_{IS}\left\Vert u^{1}\right\Vert _{W^{1,p}\left( \Sigma \right) }%
\right] ^{\gamma _{1}-p}R_{0}^{\varepsilon N}\frac{\left\vert k\right\vert
^{p}}{R^{\varepsilon N}}\left\vert A_{k,s}^{1}\right\vert ^{1-\frac{p}{N}%
+\varepsilon }
\end{eqnarray*}%
\ and $\left\Vert b_{1}\right\Vert _{L^{\sigma _{1}}\left( \Sigma \right) }%
\left[ C_{IS}\left\Vert u^{1}\right\Vert _{W^{1,p}\left( \Sigma \right) }%
\right] ^{\gamma _{1}-p}R_{0}^{\varepsilon N}\leq \frac{1}{2}$ then\ it
follows%
\begin{equation}
\begin{tabular}{l}
$\int\limits_{E_{k,s}^{1}}b_{1}\left( x\right) \eta ^{p\gamma
_{1}}\left\vert k\right\vert ^{\gamma _{1}}\,dx$ \\ 
$\leq \frac{1}{2}\frac{\left\vert k\right\vert ^{p}}{R^{\varepsilon N}}%
\left\vert A_{k,s}^{1}\right\vert ^{1-\frac{p}{N}+\varepsilon }$%
\end{tabular}
\label{76}
\end{equation}%
Using (\ref{71}) and (\ref{76}) we have%
\begin{equation}
\begin{tabular}{l}
$\frac{1}{2}\int\limits_{E_{k,s}^{1}}\left\vert \nabla u^{1}\right\vert
^{p}+b_{1}\left( x\right) \left( u^{1}\right) ^{\gamma _{1}}\,dx$ \\ 
$\leq 2c\left( p,\gamma \right) \int\limits_{E_{k,s}^{1}}\left( 1-\eta
^{p}\right) \left( \left\vert \nabla u^{1}\right\vert ^{p}+b_{1}\left(
x\right) \left( u^{1}\right) ^{\gamma _{1}}\right) \,dx$ \\ 
$+\left( 2^{2p-1}p^{p}+2^{p-1}\right) \int\limits_{E_{k,s}^{1}}\frac{\left(
u^{1}-k\right) ^{p}}{\left( s-t\right) ^{p}}\,dx+2L_{1}\int%
\limits_{E_{k,s}^{1}}a_{1}\left( x\right) \,dx$ \\ 
$+\left\vert k\right\vert ^{p}R^{-\varepsilon N}\left\vert
A_{k,s}^{1}\right\vert ^{1-\frac{p}{N}+\varepsilon }$ \\ 
$+\int\limits_{E_{k,s}^{1}}\eta ^{p}\left\vert G\left( x,u+\varphi ,A\right)
-G\left( x,u+\varphi ,\nabla u\right) \right\vert \,dx$ \\ 
$+\int\limits_{E_{k,s}^{1}}\left\vert G\left( x,u+\varphi ,\nabla u\right)
-G\left( x,u,\nabla u\right) \right\vert \,dx$%
\end{tabular}
\label{77}
\end{equation}

\bigskip Now let's estimate the following term%
\begin{equation*}
\int\limits_{E_{k,s}^{1}}\left\vert G\left( x,u+\varphi ,\nabla u\right)
-G\left( x,u,\nabla u\right) \right\vert \,dx
\end{equation*}%
using hypothesis H.1.4 and H\"{o}lder's inequality we obtain%
\begin{equation}
\begin{tabular}{l}
$\int\limits_{E_{k,s}^{1}}\left\vert G\left( x,u+\varphi ,\nabla u\right)
-G\left( x,u,\nabla u\right) \right\vert \,dx$ \\ 
$\leq \int\limits_{E_{k,s}^{1}}c\left( x\right) \left\vert \nabla
u\right\vert ^{\delta }\left\vert \varphi \right\vert ^{\beta }\,dx$ \\ 
$\leq \left[ \int\limits_{E_{k,s}^{1}}\left\vert \nabla u\right\vert ^{q}\,dx%
\right] ^{\frac{\delta }{q}}\left[ \int\limits_{E_{k,s}^{1}}\left( c\left(
x\right) \right) ^{\frac{q}{q-\delta }}\left\vert \varphi \right\vert ^{%
\frac{\beta q}{q-\delta }}\,dx\right] ^{\frac{q-\delta }{q}}$%
\end{tabular}%
\end{equation}%
and using H\"{o}lder's inequality again it follows that%
\begin{equation}
\begin{tabular}{l}
$\int\limits_{E_{k,s}^{1}}\left\vert G\left( x,u+\varphi ,\nabla u\right)
-G\left( x,u,\nabla u\right) \right\vert \,dx$ \\ 
$\leq \left[ \left\vert E_{k,s}^{1}\right\vert ^{\frac{p-q}{p}}\left(
\int\limits_{E_{k,s}^{1}}\left\vert \nabla u\right\vert ^{p}\,dx\right) ^{%
\frac{q}{p}}\right] ^{\frac{\delta }{q}}$ \\ 
$\left[ \left( \int\limits_{E_{k,s}^{1}}\left( c\left( x\right) \right) ^{%
\frac{qp^{\ast }}{p^{\ast }\left( q-\delta \right) -q\beta }}\,dx\right) ^{%
\frac{p^{\ast }\left( q-\delta \right) -q\beta }{p^{\ast }\left( q-\delta
\right) }}\left( \int\limits_{E_{k,s}^{1}}\left\vert \varphi \right\vert
^{p^{\ast }}\,dx\right) ^{\frac{\beta q}{p^{\ast }\left( q-\delta \right) }}%
\right] ^{\frac{q-\delta }{q}}$ \\ 
$\leq \left\vert E_{k,s}^{1}\right\vert ^{\frac{\left( p-q\right) \delta }{pq%
}}\left( \int\limits_{E_{k,s}^{1}}\left\vert \nabla u\right\vert
^{p}\,dx\right) ^{\frac{\delta }{p}}$ \\ 
$\left[ \left[ \left\vert E_{k,s}^{1}\right\vert ^{1-\frac{qp^{\ast }}{%
\sigma \left[ p^{\ast }\left( q-\delta \right) -q\beta \right] }}\left(
\left\Vert c\right\Vert _{L^{\sigma }\left( E_{k,s}^{1}\right) }\right) ^{%
\frac{qp^{\ast }}{\left[ p^{\ast }\left( q-\delta \right) -q\beta \right] }}%
\right] ^{\frac{p^{\ast }\left( q-\delta \right) -q\beta }{p^{\ast }\left(
q-\delta \right) }}\left( \int\limits_{E_{k,s}^{1}}\left\vert \varphi
\right\vert ^{p^{\ast }}\,dx\right) ^{\frac{\beta q}{p^{\ast }\left(
q-\delta \right) }}\right] ^{\frac{q-\delta }{q}}$ \\ 
$\leq \left\vert E_{k,s}^{1}\right\vert ^{\Theta }\left(
\int\limits_{E_{k,s}^{1}}\left\vert \nabla u\right\vert ^{p}\,dx\right) ^{%
\frac{\delta }{p}}\left\Vert c\right\Vert _{L^{\sigma }\left( \Sigma \right)
}\left( \int\limits_{E_{k,s}^{1}}\left\vert \varphi \right\vert ^{p^{\ast
}}\,dx\right) ^{\frac{\beta }{p^{\ast }}}$ \\ 
\end{tabular}%
\end{equation}%
where $\Theta =\frac{\left( p-q\right) \delta }{pq}+\left( 1-\frac{qp^{\ast }%
}{\sigma \left[ p^{\ast }\left( q-\delta \right) -q\beta \right] }\right)
\left( \frac{p^{\ast }\left( q-\delta \right) -q\beta }{p^{\ast }\left(
q-\delta \right) }\right) \frac{q-\delta }{q}$; moreover, with some
algebraic calculations we obtain $\Theta =1-\left( \frac{\delta }{p}+\frac{%
\beta }{p^{\ast }}+\frac{1}{\sigma }\right) $. Using Sobolev's inequality we
have%
\begin{equation}
\begin{tabular}{l}
$\int\limits_{E_{k,s}^{1}}\left\vert G\left( x,u+\varphi ,\nabla u\right)
-G\left( x,u,\nabla u\right) \right\vert \,dx$ \\ 
$\leq C_{SN}\left\vert E_{k,s}^{1}\right\vert ^{\Theta }\left(
\int\limits_{E_{k,s}^{1}}\left\vert \nabla u\right\vert ^{p}\,dx\right) ^{%
\frac{\delta }{p}}\left\Vert c\right\Vert _{L^{\sigma }\left( \Sigma \right)
}\left( \int\limits_{B_{s}\left( x_{0}\right) }\left\vert \nabla \varphi
\right\vert ^{p}\,dx\right) ^{\frac{\beta }{p}}$ \\ 
$\leq C_{SN}\left\vert E_{k,s}^{1}\right\vert ^{\Theta }\left(
\int\limits_{E_{k,s}^{1}}\left\vert \nabla u\right\vert ^{p}\,dx\right) ^{%
\frac{\delta }{p}}\left\Vert c\right\Vert _{L^{\sigma }\left( \Sigma \right)
}\left( \int\limits_{E_{k,s}^{1}}\eta ^{p}\left\vert \nabla u^{1}\right\vert
^{p}+2^{p}\left( \frac{u^{1}-k}{s-t}\right) ^{p}\,dx\right) ^{\frac{\beta }{p%
}}$%
\end{tabular}%
\end{equation}%
Now, applying Young's inequality, it follows%
\begin{equation}
\begin{tabular}{l}
$\int\limits_{E_{k,s}^{1}}\left\vert G\left( x,u+\varphi ,\nabla u\right)
-G\left( x,u,\nabla u\right) \right\vert \,dx$ \\ 
$\leq \frac{1}{\varepsilon ^{\frac{p^{2}}{\beta \left( p-\beta \right) }}}%
\left( C_{SN}\left\vert E_{k,s}^{1}\right\vert ^{\Theta }\left(
\int\limits_{E_{k,s}^{1}}\left\vert \nabla u\right\vert ^{p}\,dx\right) ^{%
\frac{\delta }{p}}\left\Vert c\right\Vert _{L^{\sigma }\left( \Sigma \right)
}\right) ^{\frac{p}{p-\beta }}$ \\ 
$+\varepsilon \int\limits_{E_{k,s}^{1}}\eta ^{p}\left\vert \nabla
u^{1}\right\vert ^{p}+2^{p}\left( \frac{u^{1}-k}{r-t}\right) ^{p}\,dx$ \\ 
$\leq \frac{C_{\Sigma }}{\varepsilon ^{\frac{p^{2}}{\beta \left( p-\beta
\right) }}}\left\vert E_{k,s}^{1}\right\vert ^{\frac{p\Theta }{p-\beta }%
}+\varepsilon \int\limits_{E_{k,s}^{1}}\eta ^{p}\left\vert \nabla
u^{1}\right\vert ^{p}+2^{p}\left( \frac{u^{1}-k}{r-t}\right) ^{p}\,dx$%
\end{tabular}
\label{81}
\end{equation}%
where%
\begin{equation}
C_{\Sigma }=\left( C_{SN}\left( \int\limits_{\Sigma }\left\vert \nabla
u\right\vert ^{p}\,dx\right) ^{\frac{\delta }{p}}\left\Vert c\right\Vert
_{L^{\sigma }\left( \Sigma \right) }\right) ^{\frac{p}{p-\beta }}
\end{equation}%
Fixed $\varepsilon =\frac{1}{4}$, using (\ref{77}) and (\ref{81}), we get 
\begin{equation}
\begin{tabular}{l}
$\frac{1}{4}\int\limits_{E_{k,s}^{1}}\left\vert \nabla u^{1}\right\vert
^{p}+b_{1}\left( x\right) \left( u^{1}\right) ^{\gamma _{1}}\,dx$ \\ 
$\leq 2c\left( p,\gamma \right) \int\limits_{E_{k,s}^{1}}\left( 1-\eta
^{p}\right) \left( \left\vert \nabla u^{1}\right\vert ^{p}+b_{1}\left(
x\right) \left( u^{1}\right) ^{\gamma _{1}}\right) \,dx$ \\ 
$+\left( 2^{2p-1}p^{p}+2^{p}\right) \int\limits_{E_{k,s}^{1}}\frac{\left(
u^{1}-k\right) ^{p}}{\left( s-t\right) ^{p}}\,dx+2L_{1}\int%
\limits_{E_{k,s}^{1}}a_{1}\left( x\right) \,dx$ \\ 
$+\left\vert k\right\vert ^{p}R^{-\varepsilon N}\left\vert
A_{k,s}^{1}\right\vert ^{1-\frac{p}{N}+\varepsilon }$ \\ 
$+4^{\frac{p^{2}}{\beta \left( p-\beta \right) }}C_{\Sigma }\left\vert
E_{k,s}^{1}\right\vert ^{\frac{p\Theta }{p-\beta }}$ \\ 
$+\int\limits_{E_{k,s}^{1}}\eta ^{p}\left\vert G\left( x,u+\varphi ,A\right)
-G\left( x,u+\varphi ,\nabla u\right) \right\vert \,dx$%
\end{tabular}
\label{4.49}
\end{equation}%
To get our energy estimate we need to estimate the following integral 
\begin{equation*}
\int\limits_{E_{k,s}^{1}}\eta ^{p}\left\vert G\left( x,u+\varphi ,A\right)
-G\left( x,u+\varphi ,\nabla u\right) \right\vert \,dx
\end{equation*}%
from hypothesis H.1.3 and from Remark 1 we have 
\begin{equation}
\begin{tabular}{l}
$\int\limits_{E_{k,s}^{1}}\eta ^{p}\left\vert G\left( x,u+\varphi ,A\right)
-G\left( x,u+\varphi ,\nabla u\right) \right\vert \,dx$ \\ 
$\leq \int\limits_{E_{k,s}^{1}}\eta ^{p}\left\vert G\left( x,u+\varphi
,A\right) \right\vert \,dx+\int\limits_{E_{k,s}^{1}}\eta ^{p}\left\vert
G\left( x,u+\varphi ,\nabla u\right) \right\vert \,dx$ \\ 
$\leq C\int\limits_{E_{k,s}^{1}}\eta ^{p}\left\vert A\right\vert ^{q}+\eta
^{p}c\left( x\right) \left\vert A\right\vert ^{\delta }\left\vert u+\varphi
\right\vert ^{\beta }+\eta ^{p}a\left( x\right) \,dx$ \\ 
$+C\int\limits_{E_{k,s}^{1}}\eta ^{p}\left\vert \nabla u\right\vert
^{q}+\eta ^{p}c\left( x\right) \left\vert \nabla u\right\vert ^{\delta
}\left\vert u+\varphi \right\vert ^{\beta }+\eta ^{p}a\left( x\right) \,dx$
\\ 
$\leq C\int\limits_{E_{k,s}^{1}}\eta ^{p}\left\vert A\right\vert
^{q}\,dx+C\int\limits_{E_{k,s}^{1}}\eta ^{p}\left\vert \nabla u\right\vert
^{q}\,dx$ \\ 
$+C\int\limits_{E_{k,s}^{1}}\eta ^{p}c\left( x\right) \left\vert
A\right\vert ^{\delta }\left\vert u+\varphi \right\vert ^{\beta
}\,dx+C\int\limits_{E_{k,s}^{1}}\eta ^{p}c\left( x\right) \left\vert \nabla
u\right\vert ^{\delta }\left\vert u+\varphi \right\vert ^{\beta }\,dx$ \\ 
$+2C\int\limits_{E_{k,s}^{1}}\eta ^{p}a\left( x\right) \,dx$%
\end{tabular}
\label{4.50}
\end{equation}%
Now we have to estimate the various elements of (\ref{4.50}), let us start
with the following term%
\begin{equation*}
\int\limits_{E_{k,s}^{1}}\eta ^{p}\left\vert A\right\vert ^{q}\,dx
\end{equation*}%
since%
\begin{equation}
\left\vert A\right\vert \leq p\eta ^{-1}\left\vert \nabla \eta \right\vert
\left( u^{1}-k\right) +n\left\vert \nabla u\right\vert
\end{equation}%
on $E_{k,s}^{1}$, then%
\begin{equation}
\begin{tabular}{l}
$\int\limits_{E_{k,s}^{1}}\eta ^{p}\left\vert A\right\vert ^{q}\,dx$ \\ 
$\leq \int\limits_{E_{k,s}^{1}}\eta ^{p}\left\vert p\eta ^{-1}\left\vert
\nabla \eta \right\vert \left( u^{1}-k\right) +n\left\vert \nabla
u\right\vert \right\vert ^{q}\,dx$ \\ 
$\leq 2^{q-1}\int\limits_{E_{k,s}^{1}}\eta ^{p}p^{q}\eta ^{-q}\left\vert
\nabla \eta \right\vert ^{q}\left( u^{1}-k\right) ^{q}+\eta
^{p}n^{q}\left\vert \nabla u\right\vert ^{q}\,dx$ \\ 
$\leq 2^{q-1}p^{q}\int\limits_{E_{k,s}^{1}}\eta ^{p-q}\left\vert \nabla \eta
\right\vert ^{q}\left( u^{1}-k\right)
^{q}\,dx+2^{q-1}n^{q}\int\limits_{E_{k,s}^{1}}\eta ^{p}\left\vert \nabla
u\right\vert ^{q}\,dx$%
\end{tabular}
\label{86}
\end{equation}%
Using Young's inequality we can estimate the term%
\begin{equation*}
\int\limits_{E_{k,s}^{1}}\eta ^{p-q}{}\left\vert \nabla \eta \right\vert
^{q}\left( u^{1}-k\right) ^{q}\,dx
\end{equation*}%
getting 
\begin{equation}
\int\limits_{E_{k,s}^{1}}\eta ^{p-q}\left\vert \nabla \eta \right\vert
^{q}\left( u^{1}-k\right) ^{q}\,dx\leq \frac{q}{p}\int\limits_{E_{k,s}^{1}}%
\left\vert \nabla \eta \right\vert ^{p}\left( u^{1}-k\right)
^{p}\,dx+\int\limits_{E_{k,s}^{1}}\eta ^{p}\,dx  \label{87}
\end{equation}%
furthermore, with the H\"{o}lder inequality we can estimate the term%
\begin{equation*}
\int\limits_{E_{k,s}^{1}}\eta ^{p}\left\vert \nabla u\right\vert ^{q}\,dx
\end{equation*}%
getting%
\begin{equation}
\int\limits_{E_{k,s}^{1}}\eta ^{p}\left\vert \nabla u\right\vert
^{q}\,dx\leq \left[ \int\limits_{E_{k,s}^{1}}\eta ^{p}{}\,dx\right] ^{\frac{%
p-q}{p}}\cdot \left[ \int\limits_{E_{k,s}^{1}}\eta ^{p}\left\vert \nabla
u\right\vert ^{p}\,dx\right] ^{\frac{q}{p}}  \label{88}
\end{equation}%
then, remembering the properties of $\eta $\ and using (\ref{86}), (\ref{87}%
) and (\ref{88}) we have 
\begin{equation}
\begin{tabular}{l}
$\int\limits_{E_{k,s}^{1}}\eta ^{p}\left\vert A\right\vert ^{q}\,dx$ \\ 
$\leq 2^{q-1}p^{q-1}q\int\limits_{E_{k,s}^{1}-E_{k,t}^{1}}\frac{\left(
u^{1}-k\right) ^{p}}{\left( s-t\right) ^{p}}\,dx+2^{q-1}p^{q}\left\vert
E_{k,s}^{1}\right\vert $ \\ 
$+2^{q-1}n^{q}\left\vert E_{k,s}^{1}\right\vert ^{\frac{p-q}{p}}\cdot \left[
\int\limits_{E_{k,s}^{1}}\left\vert \nabla u\right\vert ^{p}\,dx\right] ^{%
\frac{q}{p}}$ \\ 
$\leq 2^{q-1}p^{q-1}q\int\limits_{E_{k,s}^{1}-E_{k,t}^{1}}\frac{\left(
u^{1}-k\right) ^{p}}{\left( s-t\right) ^{p}}\,dx+2^{q-1}p^{q}\left\vert
E_{k,s}^{1}\right\vert $ \\ 
$+2^{q-1}n^{q}\left\vert E_{k,s}^{1}\right\vert ^{\frac{p-q}{p}}\cdot \left[
\int\limits_{\Sigma }\left\vert \nabla u\right\vert ^{p}\,dx\right] ^{\frac{q%
}{p}}$%
\end{tabular}
\label{89}
\end{equation}%
To estimate the following term%
\begin{equation*}
\int\limits_{E_{k,s}^{1}}\eta ^{p}\left\vert \nabla u\right\vert ^{q}\,dx
\end{equation*}%
of (3.50) we can use the properties of $\eta $\ and the H\"{o}lder
inequality getting%
\begin{equation}
\begin{tabular}{l}
$\int\limits_{E_{k,s}^{1}}\eta ^{p}\left\vert \nabla u\right\vert ^{q}\,dx$
\\ 
$\leq \left\vert E_{k,s}^{1}\right\vert ^{\frac{p-q}{p}}\left[
\int\limits_{\Sigma }\left\vert \nabla u\right\vert ^{p}\,dx\right] ^{\frac{q%
}{p}}$%
\end{tabular}
\label{91}
\end{equation}%
We now proceed to estimate the following integral%
\begin{equation*}
\int\limits_{E_{k,s}^{1}}\eta ^{p}c\left( x\right) \left\vert A\right\vert
^{\delta }\left\vert u+\varphi \right\vert ^{\beta }\,dx
\end{equation*}%
of (3.50), using H\"{o}lder's inequality we have%
\begin{equation}
\begin{tabular}{l}
$\int\limits_{E_{k,s}^{1}}\eta ^{p}c\left( x\right) \left\vert A\right\vert
^{\delta }\left\vert u+\varphi \right\vert ^{\beta }\,dx$ \\ 
$\leq \left( \int\limits_{E_{k,s}^{1}}\eta ^{q}\left\vert A\right\vert
^{q}\,dx\right) ^{\frac{\delta }{q}}\left( \int\limits_{E_{k,s}^{1}}\left(
\eta ^{p-\delta }c\left( x\right) \left\vert u+\varphi \right\vert ^{\beta
}\right) ^{\frac{q}{q-\delta }}\,dx\right) ^{\frac{q-\delta }{q}}$ \\ 
$\leq \left[ \left\vert E_{k,s}^{1}\right\vert ^{\frac{p-q}{p}}\left(
\int\limits_{E_{k,s}^{1}}\eta ^{p}\left\vert A\right\vert ^{p}\,dx\right) ^{%
\frac{q}{p}}\right] ^{\frac{\delta }{q}}$ \\ 
$\cdot \left[ \left( \left\vert E_{k,s}^{1}\right\vert ^{1-\frac{qp^{\ast }}{%
\sigma \left[ p^{\ast }\left( q-\delta \right) -q\beta \right] }}\left(
\left\Vert c\right\Vert _{L^{\sigma }\left( E_{k,s}^{1}\right) }\right) ^{%
\frac{qp^{\ast }}{\left[ p^{\ast }\left( q-\delta \right) -q\beta \right] }%
}\right) ^{\frac{p^{\ast }\left( q-\delta \right) -q\beta }{p^{\ast }\left(
q-\delta \right) }}\cdot \left( \int\limits_{E_{k,s}^{1}}\left( \eta ^{\frac{%
p-\delta }{\beta }}\left\vert u+\varphi \right\vert \right) ^{p^{\ast
}}\,dx\right) ^{\frac{\beta q}{p^{\ast }\left( q-\delta \right) }}\right] ^{%
\frac{q-\delta }{q}}$%
\end{tabular}
\label{92}
\end{equation}%
Remembering that $\Theta =1-\left( \frac{\delta }{p}+\frac{\beta }{p^{\ast }}%
+\frac{1}{\sigma }\right) $ we obtain 
\begin{equation}
\begin{tabular}{l}
$\int\limits_{E_{k,s}^{1}}\eta ^{p}c\left( x\right) \left\vert A\right\vert
^{\delta }\left\vert u+\varphi \right\vert ^{\beta }\,dx$ \\ 
$\leq \left\vert E_{k,s}^{1}\right\vert ^{\Theta }\left\Vert c\right\Vert
_{L^{\sigma }\left( \Sigma \right) }\left( \int\limits_{E_{k,s}^{1}}\eta
^{p}\left\vert A\right\vert ^{p}\,dx\right) ^{\frac{\delta }{p}}$ \\ 
$\cdot \left( \int\limits_{E_{k,s}^{1}}\left( \eta ^{\frac{p-\delta }{\beta }%
}\left\vert u+\varphi \right\vert \right) ^{p^{\ast }}\,dx\right) ^{\frac{%
\beta }{p^{\ast }}}$%
\end{tabular}
\label{4.58}
\end{equation}%
Using (3.2) it follows 
\begin{equation}
\int\limits_{E_{k,s}^{1}}\eta ^{p}\left\vert A\right\vert ^{p}\,dx\leq
C_{1,p,q}\int\limits_{E_{k,s}^{1}-E_{k,t}^{1}}\frac{\left( u^{1}-k\right)
^{p}}{\left( s-t\right) ^{p}}\,dx+C_{2,p,q}\left\vert E_{k,s}^{1}\right\vert
+C_{3,p,q}\left\vert E_{k,s}^{1}\right\vert ^{\frac{p-q}{p}}\left[
\int\limits_{\Sigma }\left\vert \nabla u\right\vert ^{p}\,dx\right] ^{\frac{q%
}{p}}  \label{4.59}
\end{equation}%
and by the properties of $\eta $\ we get 
\begin{equation}
\begin{tabular}{ll}
$\left( \int\limits_{E_{k,s}^{1}}\left( \eta ^{\frac{p-\delta }{\beta }%
}\left\vert u+\varphi \right\vert \right) ^{p^{\ast }}\,dx\right) ^{\frac{1}{%
p^{\ast }}}$ & $\leq \left( \int\limits_{E_{k,s}^{1}}\left\vert u+\varphi
\right\vert ^{p^{\ast }}\,dx\right) ^{\frac{1}{p^{\ast }}}$ \\ 
& $\leq \left( \int\limits_{E_{k,s}^{1}}\left\vert u\right\vert ^{p^{\ast
}}\,dx\right) ^{\frac{1}{p^{\ast }}}+\left(
\int\limits_{E_{k,s}^{1}}\left\vert \varphi \right\vert ^{p^{\ast
}}\,dx\right) ^{\frac{1}{p^{\ast }}}$%
\end{tabular}
\label{4.60}
\end{equation}%
by the Sobolev Inequality and the Immersion theorem it follows 
\begin{equation}
\begin{tabular}{l}
$\left( \int\limits_{E_{k,s}^{1}}\left( \eta ^{\frac{p-\delta }{\beta }%
}\left\vert u+\varphi \right\vert \right) ^{p^{\ast }}\,dx\right) ^{\frac{1}{%
p^{\ast }}}$ \\ 
$\leq C\left\Vert u\right\Vert _{W^{1,p}\left( \Sigma ,%
\mathbb{R}
^{n}\right) }+$ \\ 
$+C_{SN}\left( \int\limits_{E_{k,s}^{1}}\left\vert \nabla \varphi
\right\vert ^{p}\,dx\right) ^{\frac{1}{p}}$%
\end{tabular}
\label{4.61}
\end{equation}%
Now, using (\ref{4.58}), (\ref{4.59}), (\ref{4.60}) and (\ref{4.61}) we
obtain 
\begin{equation}
\begin{tabular}{l}
$\int\limits_{E_{k,s}^{1}}\eta ^{p}c\left( x\right) \left\vert A\right\vert
^{\delta }\left\vert u+\varphi \right\vert ^{\beta }\,dx$ \\ 
$\leq \left\vert E_{k,s}^{1}\right\vert ^{\Theta }\left\Vert c\right\Vert
_{L^{\sigma }\left( \Sigma \right) }\left( C\left\Vert u\right\Vert
_{W^{1,p}\left( \Sigma ,%
\mathbb{R}
^{n}\right) }+C_{SN}\left( \int\limits_{E_{k,s}^{1}}\left\vert \nabla
\varphi \right\vert ^{p}\,dx\right) ^{\frac{1}{p}}\right) ^{\beta }$ \\ 
$\cdot \left( C_{1,p,q}\int\limits_{E_{k,s}^{1}-E_{k,t}^{1}}\frac{\left(
u^{1}-k\right) ^{p}}{\left( s-t\right) ^{p}}\,dx+C_{2,p,q}\left\vert
E_{k,s}^{1}\right\vert +C_{3,p,q}\left\vert E_{k,s}^{1}\right\vert ^{\frac{%
p-q}{p}}\left[ \int\limits_{\Sigma }\left\vert \nabla u\right\vert ^{p}\,dx%
\right] ^{\frac{q}{p}}\right) ^{\frac{\delta }{p}}$%
\end{tabular}%
\end{equation}%
Since 
\begin{equation}
\begin{tabular}{l}
$\left( C\left\Vert u\right\Vert _{W^{1,p}\left( \Sigma ,%
\mathbb{R}
^{n}\right) }+C_{SN}\left( \int\limits_{E_{k,s}^{1}}\left\vert \nabla
\varphi \right\vert ^{p}\,dx\right) ^{\frac{1}{p}}\right) ^{\beta }$ \\ 
$\leq \left( C^{\beta }\left\Vert u\right\Vert _{W^{1,p}\left( \Sigma ,%
\mathbb{R}
^{n}\right) }^{\beta }+C_{SN}^{\beta }\left(
\int\limits_{E_{k,s}^{1}}\left\vert \nabla \varphi \right\vert
^{p}\,dx\right) ^{\frac{\beta }{p}}\right) $%
\end{tabular}
\label{4.63}
\end{equation}%
and 
\begin{equation}
\begin{tabular}{l}
$\left( C_{1,p,q}\int\limits_{E_{k,s}^{1}-E_{k,t}^{1}}\frac{\left(
u^{1}-k\right) ^{p}}{\left( s-t\right) ^{p}}\,dx+C_{2,p,q}\left\vert
E_{k,s}^{1}\right\vert +C_{3,p,q}\left\vert E_{k,s}^{1}\right\vert ^{\frac{%
p-q}{p}}\left[ \int\limits_{\Sigma }\left\vert \nabla u\right\vert ^{p}\,dx%
\right] ^{\frac{q}{p}}\right) ^{\frac{\delta }{p}}$ \\ 
$\leq \left( C_{1,p,q}^{\frac{\delta }{p}}\left(
\int\limits_{E_{k,s}^{1}-E_{k,t}^{1}}\frac{\left( u^{1}-k\right) ^{p}}{%
\left( s-t\right) ^{p}}\,dx\right) ^{\frac{\delta }{p}}+B_{p,q,\delta
}\right) $%
\end{tabular}
\label{4.64}
\end{equation}%
where 
\begin{equation*}
B_{p,q,\delta }=C_{2,p,q}^{\frac{\delta }{p}}\left\vert
E_{k,s}^{1}\right\vert ^{\frac{\delta }{p}}+C_{3,p,q}^{\frac{\delta }{p}%
}\left\vert E_{k,s}^{1}\right\vert ^{\frac{\delta }{p}\frac{p-q}{p}}\left[
\left\Vert u\right\Vert _{W^{1,p}\left( \Sigma ,%
\mathbb{R}
^{n}\right) }\right] ^{\frac{\delta q}{p}}
\end{equation*}%
we get 
\begin{equation}
\begin{tabular}{l}
$\int\limits_{E_{k,s}^{1}}\eta ^{p}c\left( x\right) \left\vert A\right\vert
^{\delta }\left\vert u+\varphi \right\vert ^{\beta }\,dx$ \\ 
$\leq \left\vert E_{k,s}^{1}\right\vert ^{\Theta }\left\Vert c\right\Vert
_{L^{\sigma }\left( \Sigma \right) }\left( C^{\beta }\left\Vert u\right\Vert
_{W^{1,p}\left( \Sigma ,%
\mathbb{R}
^{n}\right) }^{\beta }+C_{SN}^{\beta }\left(
\int\limits_{E_{k,s}^{1}}\left\vert \nabla \varphi \right\vert
^{p}\,dx\right) ^{\frac{\beta }{p}}\right) $ \\ 
$\cdot \left( C_{1,p,q}^{\frac{\delta }{p}}\left(
\int\limits_{E_{k,s}^{1}-E_{k,t}^{1}}\frac{\left( u^{1}-k\right) ^{p}}{%
\left( s-t\right) ^{p}}\,dx\right) ^{\frac{\delta }{p}}+B_{p,q,\delta
}\right) $%
\end{tabular}
\label{4.65}
\end{equation}%
Since 
\begin{equation}
\int\limits_{E_{k,s}^{1}}\left\vert \nabla \varphi \right\vert ^{p}\,dx\leq
\int\limits_{E_{k,s}^{1}}\eta ^{p}\left\vert \nabla u^{1}\right\vert
^{p}\,dx+2^{p}\int\limits_{E_{k,s}^{1}}\left( \frac{u^{1}-k}{s-t}\right)
^{p}\,dx  \label{4.66}
\end{equation}%
from (\ref{4.65}) we obtain 
\begin{equation}
\begin{tabular}{l}
$\int\limits_{E_{k,s}^{1}}\eta ^{p}c\left( x\right) \left\vert A\right\vert
^{\delta }\left\vert u+\varphi \right\vert ^{\beta }\,dx$ \\ 
$\leq \left\vert E_{k,s}^{1}\right\vert ^{\Theta }\left\Vert c\right\Vert
_{L^{\sigma }\left( \Sigma \right) }\left( D_{p,\beta }+2^{\beta
}C_{SN}^{\beta }\left( \int\limits_{E_{k,s}^{1}}\left( \frac{u^{1}-k}{s-t}%
\right) ^{p}\,dx\right) ^{\frac{\beta }{p}}\right) $ \\ 
$\cdot \left( C_{1,p,q}^{\frac{\delta }{p}}\left(
\int\limits_{E_{k,s}^{1}-E_{k,t}^{1}}\frac{\left( u^{1}-k\right) ^{p}}{%
\left( s-t\right) ^{p}}\,dx\right) ^{\frac{\delta }{p}}+B_{p,q,\delta
}\right) $%
\end{tabular}
\label{4.67}
\end{equation}%
where%
\begin{equation*}
D_{p,\beta }=\left( C^{\beta }+C_{SN}^{\beta }\right) \left\Vert
u\right\Vert _{W^{1,p}\left( \Sigma ,%
\mathbb{R}
^{n}\right) }^{\beta }
\end{equation*}%
Let us define

\bigskip 
\begin{equation*}
\begin{tabular}{l}
$Z$ \\ 
$=\left\vert E_{k,s}^{1}\right\vert ^{\Theta }\left\Vert c\right\Vert
_{L^{\sigma }\left( \Sigma \right) }\left( D_{p,\beta }+2^{\beta
}C_{SN}^{\beta }\left( \int\limits_{E_{k,s}^{1}}\left( \frac{u^{1}-k}{s-t}%
\right) ^{p}\,dx\right) ^{\frac{\beta }{p}}\right) $ \\ 
$\cdot \left( C_{1,p,q}^{\frac{\delta }{p}}\left(
\int\limits_{E_{k,s}^{1}-E_{k,t}^{1}}\frac{\left( u^{1}-k\right) ^{p}}{%
\left( s-t\right) ^{p}}\,dx\right) ^{\frac{\delta }{p}}+B_{p,q,\delta
}\right) $%
\end{tabular}%
\end{equation*}%
it follows 
\begin{equation*}
\begin{tabular}{l}
$Z$ \\ 
$=\left\vert E_{k,s}^{1}\right\vert ^{\Theta }\left\Vert c\right\Vert
_{L^{\sigma }\left( \Sigma \right) }D_{p,\beta }\,C_{1,p,q}^{\frac{\delta }{p%
}}\left( \int\limits_{E_{k,s}^{1}-E_{k,t}^{1}}\frac{\left( u^{1}-k\right)
^{p}}{\left( s-t\right) ^{p}}\,dx\right) ^{\frac{\delta }{p}}$ \\ 
$+\left\vert E_{k,s}^{1}\right\vert ^{\Theta }\left\Vert c\right\Vert
_{L^{\sigma }\left( \Sigma \right) }D_{p,\beta }\,B_{p,q,\delta }$ \\ 
$+2^{\beta }\left\vert E_{k,s}^{1}\right\vert ^{\Theta }\left\Vert
c\right\Vert _{L^{\sigma }\left( \Sigma \right) }\,\,C_{1,p,q}^{\frac{\delta 
}{p}}C_{SN}^{\beta }\left( \int\limits_{E_{k,s}^{1}}\left( \frac{u^{1}-k}{s-t%
}\right) ^{p}\,dx\right) ^{\frac{\beta }{p}+\frac{\delta }{p}}$ \\ 
$+2^{\beta }\left\vert E_{k,s}^{1}\right\vert ^{\Theta }\left\Vert
c\right\Vert _{L^{\sigma }\left( \Sigma \right) }\,B_{p,q,\delta
}\,C_{SN}^{\beta }\left( \int\limits_{E_{k,s}^{1}}\left( \frac{u^{1}-k}{s-t}%
\right) ^{p}\,dx\right) ^{\frac{\beta }{p}}$%
\end{tabular}%
\end{equation*}%
since $\beta <p-\delta $, using Young inequality, we get

\begin{equation}
\begin{tabular}{l}
$Z$ \\ 
$\leq \left\vert E_{k,s}^{1}\right\vert ^{\frac{p}{p-\delta }\Theta
}\left\Vert c\right\Vert _{L^{\sigma }\left( \Sigma \right) }^{\frac{p}{%
p-\delta }}D_{p,\beta }^{\frac{p}{p-\delta }}\,+C_{1,p,q}\int%
\limits_{E_{k,s}^{1}-E_{k,t}^{1}}\frac{\left( u^{1}-k\right) ^{p}}{\left(
s-t\right) ^{p}}\,dx$ \\ 
$+\left\vert E_{k,s}^{1}\right\vert ^{\Theta }\left\Vert c\right\Vert
_{L^{\sigma }\left( \Sigma \right) }D_{p,\beta }\,B_{p,q,\delta }$ \\ 
$+\left( 2^{\beta }\left\vert E_{k,s}^{1}\right\vert ^{\Theta }\left\Vert
c\right\Vert _{L^{\sigma }\left( \ \Sigma \right) \ }\,\,C_{1,p,q}^{\frac{%
\delta }{p}}C_{SN}^{\beta }\right) ^{\frac{p}{p-\beta -\delta }%
}+\int\limits_{E_{k,s}^{1}}\left( \frac{u^{1}-k}{s-t}\right) ^{p}\,dx$ \\ 
$+\left[ 2^{\beta }\left\vert E_{k,s}^{1}\right\vert ^{\Theta }\left\Vert
c\right\Vert _{L^{\sigma }\left( \Sigma \right) }\,B_{p,q,\delta }\,\right]
^{\frac{p}{p-\beta }}+C_{SN}^{p}\int\limits_{E_{k,s}^{1}}\left( \frac{u^{1}-k%
}{s-t}\right) ^{p}\,dx$%
\end{tabular}
\label{4.68}
\end{equation}

Using (\ref{4.67}) and (\ref{4.68}) we have 
\begin{equation}
\begin{tabular}{l}
$\int\limits_{E_{k,s}^{1}}\eta ^{p}c\left( x\right) \left\vert A\right\vert
^{\delta }\left\vert u+\varphi \right\vert ^{\beta }\,dx$ \\ 
$\leq \left\vert E_{k,s}^{1}\right\vert ^{\frac{p}{p-\delta }\Theta
}\left\Vert c\right\Vert _{L^{\sigma }\left( \Sigma \right) }^{\frac{p}{%
p-\delta }}D_{p,\beta }^{\frac{p}{p-\delta }}+C_{1,p,q}\int%
\limits_{E_{k,s}^{1}-E_{k,t}^{1}}\frac{\left( u^{1}-k\right) ^{p}}{\left(
s-t\right) ^{p}}\,dx$ \\ 
$+\left\vert E_{k,s}^{1}\right\vert ^{\Theta }\left\Vert c\right\Vert
_{L^{\sigma }\left( \Sigma \right) }D_{p,\beta }\,B_{p,q,\delta }$ \\ 
$+\left( 2^{\beta }\left\vert E_{k,s}^{1}\right\vert ^{\Theta }\left\Vert
c\right\Vert _{L^{\sigma }\left( \ \Sigma \right) \ }\,\,C_{1,p,q}^{\frac{%
\delta }{p}}C_{SN}^{\beta }\right) ^{\frac{p}{p-\beta -\delta }%
}+\int\limits_{E_{k,s}^{1}}\left( \frac{u^{1}-k}{s-t}\right) ^{p}\,dx$ \\ 
$+\left[ 2^{\beta }\left\vert E_{k,s}^{1}\right\vert ^{\Theta }\left\Vert
c\right\Vert _{L^{\sigma }\left( \Sigma \right) }\,B_{p,q,\delta }\,\right]
^{\frac{p}{p-\beta }}+C_{SN}^{p}\int\limits_{E_{k,s}^{1}}\left( \frac{u^{1}-k%
}{s-t}\right) ^{p}\,dx$%
\end{tabular}
\label{4.69}
\end{equation}%
Since $\frac{p}{p-\delta },\frac{p}{p-\beta -\delta },\frac{p}{p-\beta }>1$
we can write (\ref{4.69}) in the following simpler form 
\begin{equation}
\begin{tabular}{l}
$\int\limits_{E_{k,s}^{1}}\eta ^{p}c\left( x\right) \left\vert A\right\vert
^{\delta }\left\vert u+\varphi \right\vert ^{\beta }\,dx$ \\ 
$\leq C_{1,\Sigma }\left\vert E_{k,s}^{1}\right\vert ^{\Theta }$ \\ 
$+C_{2,\Sigma }\int\limits_{E_{k,s}^{1}-E_{k,t}^{1}}\frac{\left(
u^{1}-k\right) ^{p}}{\left( s-t\right) ^{p}}\,dx$%
\end{tabular}
\label{4.70}
\end{equation}%
where%
\begin{equation}
\begin{tabular}{l}
$C_{1,\Sigma }=\left\Vert c\right\Vert _{L^{\sigma }\left( \Sigma \right)
}D_{p,\beta }\,B_{p,q,\delta }+\varpi _{N}^{\Theta \frac{\delta }{p-\delta }%
}\left\Vert c\right\Vert _{L^{\sigma }\left( \Sigma \right) }^{\frac{p}{%
p-\delta }}D_{p,\beta }^{\frac{p}{p-\delta }}$ \\ 
$+\varpi _{N}^{\Theta \frac{\beta +\delta }{p-\beta -\delta }}\left(
2^{\beta }\left\Vert c\right\Vert _{L^{\sigma }\left( \ \Sigma \right) \
}\,\,C_{1,p,q}^{\frac{\delta }{p}}C_{SN}^{\beta }\right) ^{\frac{p}{p-\beta
-\delta }}+\varpi _{N}^{\Theta \frac{\beta }{p-\beta }}\left[ 2^{\beta
}\left\Vert c\right\Vert _{L^{\sigma }\left( \Sigma \right) }\,B_{p,q,\delta
}\,\right] ^{\frac{p}{p-\beta }}$%
\end{tabular}%
\end{equation}%
and%
\begin{equation}
C_{2,\Sigma }=\left( C_{1,p,q}+1+C_{1,p,q}^{p}\right)
\end{equation}%
To finish estimating all the terms of (3.50) it remains for us to consider
the term%
\begin{equation*}
\int\limits_{E_{k,s}^{1}}\eta ^{p}c\left( x\right) \left\vert \nabla
u\right\vert ^{\delta }\left\vert u+\varphi \right\vert ^{\beta }\,dx
\end{equation*}%
of (3.50), to estimate this term we will proceed in a similar way to the
previous case; using the H\"{o}lder inequality we have%
\begin{equation}
\begin{tabular}{l}
$\int\limits_{E_{k,s}^{1}}\eta ^{p}c\left( x\right) \left\vert \nabla
u\right\vert ^{\delta }\left\vert u+\varphi \right\vert ^{\beta }\,dx$ \\ 
$\leq \left( \int\limits_{E_{k,s}^{1}}\eta ^{q}\left\vert \nabla
u\right\vert ^{q}\,dx\right) ^{\frac{\delta }{q}}\left(
\int\limits_{E_{k,s}^{1}}\left( \eta ^{p-\delta }c\left( x\right) \left\vert
u+\varphi \right\vert ^{\beta }\right) ^{\frac{q}{q-\delta }}\,dx\right) ^{%
\frac{q-\delta }{q}}$ \\ 
$\leq \left[ \left\vert E_{k,s}^{1}\right\vert ^{\frac{p-q}{p}}\left(
\int\limits_{E_{k,s}^{1}}\eta ^{p}\left\vert \nabla u\right\vert
^{p}\,dx\right) ^{\frac{q}{p}}\right] ^{\frac{\delta }{q}}$ \\ 
$\cdot \left[ \left( \left\vert E_{k,s}^{1}\right\vert ^{1-\frac{qp^{\ast }}{%
\sigma \left[ p^{\ast }\left( q-\delta \right) -q\beta \right] }}\left(
\left\Vert c\right\Vert _{L^{\sigma }\left( \Sigma \right) }\right) ^{\frac{%
qp^{\ast }}{\left[ p^{\ast }\left( q-\delta \right) -q\beta \right] }%
}\right) ^{\frac{p^{\ast }\left( q-\delta \right) -q\beta }{p^{\ast }\left(
q-\delta \right) }}\cdot \left( \int\limits_{E_{k,s}^{1}}\left( \eta ^{\frac{%
p-\delta }{\beta }}\left\vert u+\varphi \right\vert \right) ^{p^{\ast
}}\,dx\right) ^{\frac{\beta q}{p^{\ast }\left( q-\delta \right) }}\right] ^{%
\frac{q-\delta }{q}}$%
\end{tabular}%
\end{equation}%
Remembering that $\Theta =1-\left( \frac{\delta }{p}+\frac{\beta }{p^{\ast }}%
+\frac{1}{\sigma }\right) $ we have 
\begin{equation}
\begin{tabular}{l}
$\int\limits_{E_{k,s}^{1}}\eta ^{p}c\left( x\right) \left\vert \nabla
u\right\vert ^{\delta }\left\vert u+\varphi \right\vert ^{\beta }\,dx$ \\ 
$\leq \left\vert E_{k,s}^{1}\right\vert ^{\Theta }\left\Vert c\right\Vert
_{L^{\sigma }\left( \Sigma \right) }\left( \int\limits_{E_{k,s}^{1}}\eta
^{p}\left\vert \nabla u\right\vert ^{p}\,dx\right) ^{\frac{\delta }{p}}$ \\ 
$\cdot \left( \int\limits_{E_{k,s}^{1}}\left( \eta ^{\frac{p-\delta }{\beta }%
}\left\vert u+\varphi \right\vert \right) ^{p^{\ast }}\,dx\right) ^{\frac{%
\beta }{p^{\ast }}}$%
\end{tabular}%
\end{equation}%
using (\ref{4.60}), (\ref{4.61}), (\ref{4.63}), (\ref{4.64}), (\ref{4.66})
and (\ref{4.68}) we obtain 
\begin{equation}
\begin{tabular}{l}
$\int\limits_{E_{k,s}^{1}}\eta ^{p}c\left( x\right) \left\vert \nabla
u\right\vert ^{\delta }\left\vert u+\varphi \right\vert ^{\beta }\,dx$ \\ 
$\leq C_{1,\Sigma }\left\vert E_{k,s}^{1}\right\vert ^{\Theta }$ \\ 
$+C_{2,\Sigma }\int\limits_{E_{k,s}^{1}-E_{k,t}^{1}}\frac{\left(
u^{1}-k\right) ^{p}}{\left( s-t\right) ^{p}}\,dx$%
\end{tabular}
\label{4.75}
\end{equation}%
Now, using (\ref{4.49}), (3.50), (3.55), (3.56), (\ref{4.70}) and (\ref{4.75}%
)\ we get 
\begin{equation}
\begin{tabular}{l}
$\frac{1}{4}\int\limits_{E_{k,s}^{1}}\left\vert \nabla u^{1}\right\vert
^{p}+b_{1}\left( x\right) \left( u^{1}\right) ^{\gamma _{1}}\,dx$ \\ 
$\leq 2c\left( p,\gamma \right) \int\limits_{E_{k,s}^{1}}\left( 1-\eta
^{p}\right) \left( \left\vert \nabla u^{1}\right\vert ^{p}+b_{1}\left(
x\right) \left( u^{1}\right) ^{\gamma _{1}}\right) \,dx$ \\ 
$+\left( 2^{2p-1}p^{p}+2^{p}\right) \int\limits_{E_{k,s}^{1}}\frac{\left(
u^{1}-k\right) ^{p}}{\left( s-t\right) ^{p}}\,dx+2L_{1}\int%
\limits_{E_{k,s}^{1}}a_{1}\left( x\right) \,dx$ \\ 
$+\left\vert k\right\vert ^{p}R^{-\varepsilon N}\left\vert
A_{k,s}^{1}\right\vert ^{1-\frac{p}{N}+\varepsilon }+\left\vert
E_{k,s}^{1}\right\vert ^{\frac{p-q}{p}}\left[ \int\limits_{\Sigma
}\left\vert \nabla u\right\vert ^{p}\,dx\right] ^{\frac{q}{p}%
}+2^{q-1}n^{q}\left\vert E_{k,s}^{1}\right\vert ^{\frac{p-q}{p}}\cdot \left[
\int\limits_{\Sigma }\left\vert \nabla u\right\vert ^{p}\,dx\right] ^{\frac{q%
}{p}}$ \\ 
$+4^{\frac{p^{2}}{\beta \left( p-\beta \right) }}C_{\Sigma }\left\vert
E_{k,s}^{1}\right\vert ^{\frac{p\Theta }{p-\beta }}+2C_{1,\Sigma }\left\vert
E_{k,s}^{1}\right\vert ^{\Theta }+2C_{2,\Sigma
}\int\limits_{E_{k,s}^{1}-E_{k,t}^{1}}\frac{\left( u^{1}-k\right) ^{p}}{%
\left( s-t\right) ^{p}}\,dx$ \\ 
$+2^{q-1}p^{q-1}q\int\limits_{E_{k,s}^{1}-E_{k,t}^{1}}\frac{\left(
u^{1}-k\right) ^{p}}{\left( s-t\right) ^{p}}\,dx+2^{q-1}p^{q}\left\vert
E_{k,s}^{1}\right\vert $%
\end{tabular}
\label{4.76}
\end{equation}%
Since $\frac{p}{p-\beta }>1$, from hypothesis H.1.4, estimating $%
\int\limits_{E_{k,s}^{1}}a_{1}\left( x\right) \,dx<\left\Vert a\right\Vert
_{L^{\kappa _{1}}\left( \Sigma \right) }\left\vert E_{k,s}^{1}\right\vert
^{1-\frac{1}{\kappa _{1}}}$ with the H\"{o}lder inequality and observing
that $E_{k,t}^{1}=A_{k,t}^{1}$ and $E_{k,s}^{1}\subset A_{k,s}^{1}$ from (%
\ref{4.76}) we obtain 
\begin{equation}
\begin{tabular}{l}
$\frac{1}{4}\int\limits_{A_{k,t}^{1}}\left\vert \nabla u^{1}\right\vert
^{p}+b_{1}\left( x\right) \left( u^{1}\right) ^{\gamma _{1}}\,dx$ \\ 
$\leq 2c\left( p,\gamma \right) \int\limits_{A_{k,s}^{1}\backslash
A_{k,t}^{1}}\left( \left\vert \nabla u^{1}\right\vert ^{p}+b_{1}\left(
x\right) \left( u^{1}\right) ^{\gamma _{1}}\right) \,dx$ \\ 
$+D_{1,\Sigma }\int\limits_{A_{k,s}^{1}}\frac{\left( u^{1}-k\right) ^{p}}{%
\left( s-t\right) ^{p}}\,dx+\left( D_{2,\Sigma }+\left\vert k\right\vert
^{p}R^{-\varepsilon N}\right) \left\vert A_{k,s}^{1}\right\vert ^{1-\frac{p}{%
N}+\varepsilon }$%
\end{tabular}
\label{4.77}
\end{equation}

By adding the quantity 
\begin{equation*}
2c\left( p,\gamma \right) \int\limits_{A_{k,t}^{1}}\left\vert \nabla
u^{1}\right\vert ^{p}+b_{1}\left( x\right) \left( u^{1}\right) ^{\gamma
_{1}}\,dx
\end{equation*}%
to both terms of (\ref{4.77}) and carrying out simple algebraic
calculations, we have

\begin{equation*}
\begin{tabular}{l}
$\int\limits_{A_{k,t}^{1}}\left\vert \nabla u^{1}\right\vert
^{p}+b_{1}\left( x\right) \left( u^{1}\right) ^{\gamma _{1}}\,dx$ \\ 
$\leq \frac{8c\left( p,\gamma \right) }{1+8c\left( p,\gamma \right) }%
\int\limits_{A_{k,s}^{1}}\left( \left\vert \nabla u^{1}\right\vert
^{p}+b_{1}\left( x\right) \left( u^{1}\right) ^{\gamma _{1}}\right) \,dx$ \\ 
$+B_{1}\int\limits_{A_{k,s}^{1}}\frac{\left( u^{1}-k\right) ^{p}}{\left(
s-t\right) ^{p}}\,dx+B_{2}\left( 1+\left\vert k\right\vert
^{p}R^{-\varepsilon N}\right) \left\vert A_{k,s}^{1}\right\vert ^{1-\frac{p}{%
N}+\varepsilon }$%
\end{tabular}%
\end{equation*}%
where $B_{1}=\frac{4D_{1,\Sigma }}{1+8c\left( p,\gamma \right) }$ and $B_{2}=%
\frac{4\left( D_{2,\Sigma }+1\right) }{1+8c\left( p,\gamma \right) }$.

Now, using Lemma \ref{lem3} we get%
\begin{equation}
\begin{tabular}{l}
$\int\limits_{A_{k,\varrho }^{1}}\left\vert \nabla u^{1}\right\vert
^{p}\,dx\leq \frac{C_{C,1}}{\left( R-\varrho \right) ^{p}}%
\int\limits_{A_{k,R}^{1}}\left( u^{1}-k\right) ^{p}\,dx+C_{C,2}\left(
1+\left\vert k\right\vert ^{p}R^{-\varepsilon N}\right) \left\vert
A_{k,s}^{1}\right\vert ^{1-\frac{p}{N}+\varepsilon }$%
\end{tabular}%
\end{equation}

Since $-u$ is a local minimizer of the following integral functional

\bigskip 
\begin{equation*}
\tilde{J}\left( v,\Omega \right) =\int\limits_{\Omega }\sum\limits_{\alpha
=1}^{n}\tilde{f}_{\alpha }\left( x,v^{\alpha }\left( x\right) ,\nabla
v^{\alpha }\left( x\right) \right) +\tilde{G}\left( x,v\left( x\right)
,\nabla v\left( x\right) \right) \,dx
\end{equation*}%
where $\tilde{f}_{\alpha }\left( x,v^{\alpha }\left( x\right) ,\nabla
v^{\alpha }\left( x\right) \right) =f_{\alpha }\left( x,-v^{\alpha }\left(
x\right) ,-\nabla v^{\alpha }\left( x\right) \right) $ and $\tilde{G}\left(
x,v\left( x\right) ,\nabla v\left( x\right) \right) =G\left( x,-v\left(
x\right) ,-\nabla v\left( x\right) \right) $ then we get%
\begin{equation}
\begin{tabular}{l}
$\int\limits_{B_{k,\varrho }^{1}}\left\vert \nabla u^{1}\right\vert
^{p}\,dx\leq \frac{C_{C,1}}{\left( R-\varrho \right) ^{p}}%
\int\limits_{B_{k,R}^{1}}\left( k-u^{1}\right) ^{p}\,dx+C_{C,2}\left(
1+\left\vert k\right\vert ^{p}R^{-\varepsilon N}\right) \left\vert
B_{k,s}^{1}\right\vert ^{1-\frac{p}{N}+\varepsilon }$%
\end{tabular}%
\end{equation}%
Similarly we can proceed for $u^{\alpha }$ with $\alpha =2,...,n$.\bigskip

\subsection{The case of H.1.4 (bis)}

The proof follows as in the previous section 3.1 only with some
modifications in the choice of the parameters of the H\"{o}lder inequalities.

\section{Proof of Theorem 1}

Let $u\in W^{1,p}\left( \Omega ;%
\mathbb{R}
^{n}\right) $ a minimizer of the functional (1.1) then by Theorem 2 it
follows that $u^{\alpha }\in DG\left( \Omega ,p,\lambda ,\lambda _{\ast
},\chi ,\varepsilon ,R_{0},k_{0}\right) $ for $\alpha =1,...,n$ \ then by
Theorem 6 it follows that $u\in L_{loc}^{\infty }\left( \Omega ;%
\mathbb{R}
^{n}\right) $; moreover if $\Sigma $ is a compac subset of $\Omega $ then $%
u\in L^{\infty }\left( \Sigma ;%
\mathbb{R}
^{n}\right) $ and%
\begin{equation}
\left\vert u\right\vert \leq M=\sqrt{\sum\limits_{\alpha =1}^{n}\left(
M^{\alpha }\right) ^{2}}
\end{equation}%
where%
\begin{equation}
M^{\alpha }=\sup_{\Sigma }\left\{ \left\vert u^{\alpha }\right\vert \right\}
\end{equation}%
Since H.1.1 and (1.2) hold then proceeding as in [22] we get%
\begin{equation}
\begin{tabular}{l}
$\int\limits_{A_{k,\varrho }^{\alpha }}\left\vert \nabla u^{\alpha
}\right\vert ^{p}\,dx\leq \frac{\tilde{C}_{C,1}}{\left( R-\varrho \right)
^{p}}\int\limits_{A_{k,R}^{\alpha }}\left( u^{\alpha }-k\right) ^{p}\,dx+%
\tilde{C}_{C,2}\left\vert A_{k,s}^{\alpha }\right\vert ^{1-\frac{p}{N}%
+\varepsilon }$%
\end{tabular}%
\end{equation}%
and%
\begin{equation}
\int\limits_{B_{k,\varrho }^{\alpha }}\left\vert \nabla u^{\alpha
}\right\vert ^{p}\,dx\leq \frac{\tilde{C}_{C,1}}{\left( R-\varrho \right)
^{p}}\int\limits_{B_{k,R}^{\alpha }}\left( k-u^{\alpha }\right) ^{p}\,dx+%
\tilde{C}_{C,2}\left\vert B_{k,s}^{\alpha }\right\vert ^{1-\frac{p}{N}%
+\varepsilon }
\end{equation}%
for every $\alpha =1,...,n$. Theorem 1 follows using (4.3), (4.4) and
Proposition 7.1, Lemma 7.2 and Theorem 7.6 of [24].

\bigskip

\section{Proof of Theorm 5}

\subsection{The case H.2.4}

\bigskip The proof is almost identical to the one presented in section 3
except for some inequalities that we are going to mention.

Let us consider $y\in \Omega $ then we fix $R_{0}=\frac{1}{4}\min \left\{ 
\frac{1}{\sqrt[N]{\varpi _{N}}},dist\left( \partial \Omega ,y\right)
\right\} $, where\ $\varpi _{N}=\left\vert B_{1}\left( 0\right) \right\vert $%
, and we define $\Sigma =\left\{ x\in \Omega :\left\vert x-y\right\vert \leq
R_{0}\right\} $. We fix $x_{0}\in \Sigma $, $R_{1}=\frac{1}{4}dist\left(
\partial \Sigma ,x_{0}\right) $, $R_{0}<\min\limits_{\alpha =1,...,n}\left\{
R_{1},A_{\alpha },B_{\alpha }\right\} $, where

$A_{\alpha }=\frac{1}{\left[ 2\cdot 4^{\gamma _{\alpha }}L_{\alpha
}C_{N,p,\gamma _{\alpha }}D^{\frac{p^{\ast }-\gamma _{\alpha }}{p^{\ast }}%
}\left\Vert b_{\alpha }\right\Vert _{L^{\sigma _{\alpha }}\left( \Sigma
\right) }\left\Vert u^{\alpha }\right\Vert _{W^{1,p}\left( \Sigma \right)
}^{\gamma _{\alpha }-p}\varpi _{N}^{\epsilon _{\alpha }}\right] ^{\frac{1}{%
\epsilon _{\alpha }N}}}$, $B_{\alpha }=\frac{1}{\left[ 2D^{p^{\ast }-\gamma
_{\alpha }}\left\Vert b_{\alpha }\right\Vert _{L^{\sigma _{\alpha }}\left(
\Sigma \right) }\left\Vert u^{\alpha }\right\Vert _{W^{1,p}\left( \Sigma
\right) }^{\gamma _{\alpha }-p}\right] ^{\frac{1}{\epsilon _{\alpha }N}}}$, $%
C_{N,p,\gamma _{\alpha }}$ and $D^{\frac{p^{\ast }-\gamma _{\alpha }}{%
p^{\ast }}}$ are universal positive constants.

We fix $0<\varrho \leq t<s\leq R<R_{0}$, $B_{z}\left( x_{0}\right) =\left\{
x:\left\vert x-x_{0}\right\vert <z\right\} $, $k\in 
\mathbb{R}
$\ and we choose $\eta \in C_{c}^{\infty }\left( B_{s}\left( x_{0}\right)
\right) $\ such that $\eta =1$\ on $B_{t}\left( x_{0}\right) $, $0\leq \eta
\leq 1$\ on $B_{s}\left( x_{0}\right) $\ and $\left\vert \nabla \eta
\right\vert \leq \frac{2}{s-t}$\ on $B_{s}\left( x_{0}\right) $. Let us
define%
\begin{equation*}
\varphi =-\eta ^{p}w
\end{equation*}%
where $w\in W^{1,p}\left( \Sigma ,%
\mathbb{R}
^{n}\right) $ with%
\begin{equation*}
w^{1}=\max \left( u^{1}-k,0\right) ,w^{\alpha }=0,\alpha =2,...,n
\end{equation*}

\bigskip Let us observe that $\varphi =0$ $\mathcal{L}^{N}$-a.e. in $\Omega
\backslash \left( \left\{ \eta >0\right\} \cap \left\{ u^{1}>k\right\}
\right) $ thus%
\begin{equation}
\nabla u+\nabla \varphi =\nabla u
\end{equation}%
$\mathcal{L}^{N}$-a.e. in $\Omega \backslash \left( \left\{ \eta >0\right\}
\cap \left\{ u^{1}>k\right\} \right) $. Let us define 
\begin{equation}
A=\left( 
\begin{tabular}{l}
$p\eta ^{-1}\nabla \eta \left( k-u^{1}\right) $ \\ 
$\nabla u^{2}$ \\ 
$\vdots $ \\ 
$\nabla u^{n}$%
\end{tabular}%
\right)
\end{equation}%
since%
\begin{equation}
\nabla w=\left( 
\begin{tabular}{l}
$\nabla u^{1}$ \\ 
$0$ \\ 
$\vdots $ \\ 
$0$%
\end{tabular}%
\right)
\end{equation}%
$\mathcal{L}^{N}$-a.e. in $\Omega \backslash \left( \left\{ \eta >0\right\}
\cap \left\{ u^{1}>k\right\} \right) $ then we deduce that%
\begin{equation}
\nabla u+\nabla \varphi =\left( 1-\eta ^{p}\right) \nabla u+\eta ^{p}A
\end{equation}%
$\mathcal{L}^{N}$-a.e. in $\Omega \backslash \left( \left\{ \eta >0\right\}
\cap \left\{ u^{1}>k\right\} \right) $. Since $u$ is a local minimizer of
the functional (1.1) then we get%
\begin{equation}
J\left( u,\Sigma \right) \leq J\left( u+\varphi ,\Sigma \right)
\end{equation}

it is%
\begin{equation}
\begin{tabular}{l}
$\int\limits_{\Sigma }\sum\limits_{\alpha =1}^{n}f_{\alpha }\left(
x,u^{\alpha },\nabla u^{\alpha }\right) +\sum\limits_{\alpha
=1}^{n}G_{\alpha }\left( x,u,\left( adj_{n-1}\nabla u\right) ^{\alpha
}\right) \,dx$ \\ 
$\leq \int\limits_{\Sigma }\sum\limits_{\alpha =1}^{n}f_{\alpha }\left(
x,u^{\alpha }+\varphi ^{\alpha },\nabla u^{\alpha }+\nabla \varphi ^{\alpha
}\right) +\sum\limits_{\alpha =1}^{n}G_{\alpha }\left( x,u+\varphi ,\left(
adj_{n-1}\nabla u+\nabla \varphi \right) ^{\alpha }\right) \,dx$%
\end{tabular}%
\end{equation}%
and%
\begin{equation}
\begin{tabular}{l}
$\int\limits_{\Sigma }\sum\limits_{\alpha =2}^{n}f_{\alpha }\left(
x,u^{\alpha },\nabla u^{\alpha }\right) \,dx+\int\limits_{\Sigma
}f_{1}\left( x,u^{1},\nabla u^{1}\right) +\sum\limits_{\alpha
=1}^{n}G_{\alpha }\left( x,u,\left( adj_{n-1}\nabla u\right) ^{\alpha
}\right) \,dx$ \\ 
$\leq \int\limits_{\Sigma }\sum\limits_{\alpha =2}^{n}f_{\alpha }\left(
x,u^{\alpha },\nabla u^{\alpha }\right) \,dx+\int\limits_{\Sigma
}f_{1}\left( x,u^{1}+\varphi ^{1},\nabla u^{1}+\nabla \varphi ^{1}\right)
+\sum\limits_{\alpha =1}^{n}G_{\alpha }\left( x,u+\varphi ,\left(
adj_{n-1}\nabla u+\nabla \varphi \right) ^{\alpha }\right) \,dx$%
\end{tabular}%
\end{equation}%
From (5.7) we deduce%
\begin{equation}
\begin{tabular}{l}
$\int\limits_{\Sigma }f_{1}\left( x,u^{1},\nabla u^{1}\right)
+\sum\limits_{\alpha =1}^{n}G_{\alpha }\left( x,u,\left( adj_{n-1}\nabla
u\right) ^{\alpha }\right) \,dx$ \\ 
$\leq \int\limits_{\Sigma }f_{1}\left( x,u^{1}+\varphi ^{1},\nabla
u^{1}+\nabla \varphi ^{1}\right) +\sum\limits_{\alpha =1}^{n}G_{\alpha
}\left( x,u+\varphi ,\left( adj_{n-1}\nabla u+\nabla \varphi \right)
^{\alpha }\right) \,dx$ \\ 
$=\int\limits_{B_{r}\left( x_{0}\right) }f_{1}\left( x,u^{1}+\varphi
^{1},\nabla u^{1}+\nabla \varphi ^{1}\right) \,dx+\int\limits_{\Sigma
-B_{r}\left( x_{0}\right) \backslash }f_{1}\left( x,u^{1},\nabla
u^{1}\right) \,dx$ \\ 
$+\int\limits_{B_{s}\left( x_{0}\right) }\sum\limits_{\alpha
=1}^{n}G_{\alpha }\left( x,u+\varphi ,\left( adj_{n-1}\nabla u+\nabla
\varphi \right) ^{\alpha }\right) \,dx+\int\limits_{\Sigma -B_{s}\left(
x_{0}\right) \backslash }\sum\limits_{\alpha =1}^{n}G_{\alpha }\left(
x,u,\left( adj_{n-1}\nabla u\right) ^{\alpha }\right) \,dx$%
\end{tabular}%
\end{equation}%
and%
\begin{equation}
\begin{tabular}{l}
$\int\limits_{B_{s}\left( x_{0}\right) }f_{1}\left( x,u^{1},\nabla
u^{1}\right) +\sum\limits_{\alpha =1}^{n}G_{\alpha }\left( x,u,\left(
adj_{n-1}\nabla u\right) ^{\alpha }\right) \,dx$ \\ 
$\leq \int\limits_{B_{s}\left( x_{0}\right) }f_{1}\left( x,u^{1}+\varphi
^{1},\nabla u^{1}+\nabla \varphi ^{1}\right) \,dx+$ \\ 
$+\int\limits_{B_{s}\left( x_{0}\right) }\sum\limits_{\alpha
=1}^{n}G_{\alpha }\left( x,u+\varphi ,\left( adj_{n-1}\nabla u+\nabla
\varphi \right) ^{\alpha }\right) \,dx$%
\end{tabular}%
\end{equation}%
Let us define $E_{k,s}^{1}=\left\{ \eta >0\right\} \cap \left\{
u^{1}>k\right\} \cap B_{s}\left( x_{0}\right) \subset B_{s}\left(
x_{0}\right) $\ then%
\begin{equation}
\begin{tabular}{l}
$\int\limits_{E_{k,s}^{1}}f_{1}\left( x,u^{1},\nabla u^{1}\right)
+\sum\limits_{\alpha =1}^{n}G_{\alpha }\left( x,u,\left( adj_{n-1}\nabla
u\right) ^{\alpha }\right) \,dx$ \\ 
$+\int\limits_{B_{s}\left( x_{0}\right) -E_{k,s}^{1}}f_{1}\left(
x,u^{1},\nabla u^{1}\right) +\sum\limits_{\alpha =1}^{n}G_{\alpha }\left(
x,u,\left( adj_{n-1}\nabla u\right) ^{\alpha }\right) \,dx$ \\ 
$\leq \int\limits_{E_{k,s}^{1}}f_{1}\left( x,u^{1}+\varphi ^{1},\nabla
u^{1}+\nabla \varphi ^{1}\right) \,dx+\int\limits_{B_{s}\left( x_{0}\right)
-E_{k,s}^{1}}f_{1}\left( x,u^{1},\nabla u^{1}\right) \,dx$ \\ 
$+\int\limits_{E_{k,s}^{1}}\sum\limits_{\alpha =1}^{n}G_{\alpha }\left(
x,u+\varphi ,\left( adj_{n-1}\nabla u+\nabla \varphi \right) ^{\alpha
}\right) \,dx$ \\ 
$+\int\limits_{B_{s}\left( x_{0}\right) -E_{k,s}^{1}}\sum\limits_{\alpha
=1}^{n}G_{\alpha }\left( x,u,\left( adj_{n-1}\nabla u\right) ^{\alpha
}\right) \,dx$%
\end{tabular}%
\end{equation}%
and%
\begin{equation}
\begin{tabular}{l}
$\int\limits_{E_{k,s}^{1}}f_{1}\left( x,u^{1},\nabla u^{1}\right)
+\sum\limits_{\alpha =1}^{n}G_{\alpha }\left( x,u,\left( adj_{n-1}\nabla
u\right) ^{\alpha }\right) \,dx$ \\ 
$\leq \int\limits_{E_{k,s}^{1}}f_{1}\left( x,u^{1}+\varphi ^{1},\nabla
u^{1}+\nabla \varphi ^{1}\right) \,dx+$ \\ 
$+\int\limits_{E_{k,s}^{1}}\sum\limits_{\alpha =1}^{n}G_{\alpha }\left(
x,u+\varphi ,\left( adj_{n-1}\nabla u+\nabla \varphi \right) ^{\alpha
}\right) \,dx$%
\end{tabular}%
\end{equation}%
Since $\sum\limits_{\alpha =1}^{n}G_{\alpha }\left( x,s,\left( adj_{n-1}\xi
\right) ^{\alpha }\right) $ is a rank one convex function and $\nabla u-A$
is a rank one matrix then proceeding as in the section 3.1 we have%
\begin{equation*}
\begin{tabular}{l}
$\int\limits_{E_{k,s}^{1}}f_{1}\left( x,u^{1},\nabla u^{1}\right)
+\sum\limits_{\alpha =1}^{n}G_{\alpha }\left( x,u,\left( adj_{n-1}\nabla
u\right) ^{\alpha }\right) \,dx$ \\ 
$\leq \int\limits_{E_{k,s}^{1}}f_{1}\left( x,u^{1}+\varphi ^{1},\nabla
u^{1}+\nabla \varphi ^{1}\right) \,dx+$ \\ 
$+\int\limits_{E_{k,s}^{1}}\left( 1-\eta ^{p}\right) \sum\limits_{\alpha
=1}^{n}G_{\alpha }\left( x,u+\varphi ,\left( adj_{n-1}\nabla u\right)
^{\alpha }\right) +\eta ^{p}\sum\limits_{\alpha =1}^{n}G_{\alpha }\left(
x,u+\varphi ,\left( adj_{n-1}A\right) ^{\alpha }\right) \,dx$%
\end{tabular}%
\end{equation*}%
then using the growth assumptions (1.2) we get%
\begin{equation}
\begin{tabular}{l}
$\int\limits_{E_{k,s}^{1}}\left\vert \nabla u^{1}\right\vert
^{p}-b_{1}\left( x\right) \left\vert u^{1}\right\vert ^{\gamma
_{1}}-a_{1}\left( x\right) \,dx$ \\ 
$\leq L_{1}\int\limits_{E_{k,s}^{1}}\left\vert \nabla u^{1}+\nabla \varphi
^{1}\right\vert ^{p}+b_{1}\left( x\right) \left\vert u^{1}+\varphi
^{1}\right\vert ^{\gamma _{1}}+a_{1}\left( x\right) \,dx$ \\ 
$+\int\limits_{E_{k,s}^{1}}\eta ^{p}\sum\limits_{\alpha =1}^{n}\left[
G_{\alpha }\left( x,u+\varphi ,\left( adj_{n-1}A\right) ^{\alpha }\right)
-G_{\alpha }\left( x,u+\varphi ,\left( adj_{n-1}\nabla u\right) ^{\alpha
}\right) \right] \,dx$ \\ 
$+\int\limits_{E_{k,s}^{1}}\sum\limits_{\alpha =1}^{n}\left[ G_{\alpha
}\left( x,u+\varphi ,\left( adj_{n-1}\nabla u\right) ^{\alpha }\right)
-G_{\alpha }\left( x,u,\left( adj_{n-1}\nabla u\right) ^{\alpha }\right) %
\right] \,dx$%
\end{tabular}%
\end{equation}%
Moreover proceeding in a similar way to the previous \ section 1.3 we have 
\begin{equation*}
\begin{tabular}{l}
$\frac{1}{2}\int\limits_{E_{k,s}^{1}}\left\vert \nabla u^{1}\right\vert
^{p}+b_{1}\left( x\right) \left( u^{1}\right) ^{\gamma _{1}}\,dx$ \\ 
$\leq 2c\left( p,\gamma \right) \int\limits_{E_{k,s}^{1}}\left( 1-\eta
^{p}\right) \left( \left\vert \nabla u^{1}\right\vert ^{p}+b_{1}\left(
x\right) \left( u^{1}\right) ^{\gamma _{1}}\right) \,dx$ \\ 
$+\left( 2^{2p-1}p^{p}+2^{p-1}\right) \int\limits_{E_{k,s}^{1}}\frac{\left(
u^{1}-k\right) ^{p}}{\left( s-t\right) ^{p}}\,dx+2L_{1}\int%
\limits_{E_{k,s}^{1}}a_{1}\left( x\right) \,dx$ \\ 
$+\left\vert k\right\vert ^{p}R^{-\varepsilon N}\left\vert
A_{k,s}^{1}\right\vert ^{1-\frac{p}{N}+\varepsilon }$ \\ 
$+\int\limits_{E_{k,s}^{1}}\eta ^{p}\sum\limits_{\alpha =1}^{n}\left[
G_{\alpha }\left( x,u+\varphi ,\left( adj_{n-1}A\right) ^{\alpha }\right)
-G_{\alpha }\left( x,u+\varphi ,\left( adj_{n-1}\nabla u\right) ^{\alpha
}\right) \right] \,dx$ \\ 
$+\int\limits_{E_{k,s}^{1}}\sum\limits_{\alpha =1}^{n}\left[ G_{\alpha
}\left( x,u+\varphi ,\left( adj_{n-1}\nabla u\right) ^{\alpha }\right)
-G_{\alpha }\left( x,u,\left( adj_{n-1}\nabla u\right) ^{\alpha }\right) %
\right] \,dx$%
\end{tabular}%
\end{equation*}%
Now we have to estimate%
\begin{equation*}
\int\limits_{E_{k,s}^{1}}\sum\limits_{\alpha =1}^{n}\left[ G_{\alpha }\left(
x,u+\varphi ,\left( adj_{n-1}\nabla u\right) ^{\alpha }\right) -G_{\alpha
}\left( x,u,\left( adj_{n-1}\nabla u\right) ^{\alpha }\right) \right] \,dx
\end{equation*}%
remembering that $\varphi ^{1}=-\eta ^{p}w^{1}$ and that $\varphi ^{\alpha
}=0$ for $\alpha =2,...,n$ we deduce that 
\begin{equation*}
\begin{tabular}{l}
$\sum\limits_{\alpha =1}^{n}\left[ G_{\alpha }\left( x,u+\varphi ,\left(
adj_{n-1}\nabla u\right) ^{\alpha }\right) -G_{\alpha }\left( x,u,\left(
adj_{n-1}\nabla u\right) ^{\alpha }\right) \right] $ \\ 
$=G_{1}\left( x,u+\varphi ,\left( adj_{n-1}\nabla u\right) ^{1}\right)
-G_{\alpha }\left( x,u,\left( adj_{n-1}\nabla u\right) ^{1}\right) $ \\ 
$+\sum\limits_{\alpha =2}^{n}\left[ G_{\alpha }\left( x,u,\left(
adj_{n-1}\nabla u\right) ^{\alpha }\right) -G_{\alpha }\left( x,u,\left(
adj_{n-1}\nabla u\right) ^{\alpha }\right) \right] $ \\ 
$=G_{1}\left( x,u+\varphi ,\left( adj_{n-1}\nabla u\right) ^{1}\right)
-G_{1}\left( x,u,\left( adj_{n-1}\nabla u\right) ^{1}\right) $%
\end{tabular}%
\end{equation*}%
Using the hypothesis H.2.4 we have%
\begin{equation*}
\begin{tabular}{l}
$\sum\limits_{\alpha =1}^{n}\left[ G_{\alpha }\left( x,u+\varphi ,\left(
adj_{n-1}\nabla u\right) ^{\alpha }\right) -G_{\alpha }\left( x,u,\left(
adj_{n-1}\nabla u\right) ^{\alpha }\right) \right] $ \\ 
$\leq \left\vert G_{1}\left( x,u+\varphi ,\left( adj_{n-1}\nabla u\right)
^{1}\right) -G_{1}\left( x,u,\left( adj_{n-1}\nabla u\right) ^{1}\right)
\right\vert $ \\ 
$\leq c\left( x\right) \left\vert \left( adj_{n-1}\nabla u\right)
^{1}\right\vert ^{\delta }\left\vert \varphi \right\vert ^{\beta }$%
\end{tabular}%
\end{equation*}%
it follows that%
\begin{equation*}
\begin{tabular}{l}
$\int\limits_{E_{k,s}^{1}}\sum\limits_{\alpha =1}^{n}\left[ G_{\alpha
}\left( x,u+\varphi ,\left( adj_{n-1}\nabla u\right) ^{\alpha }\right)
-G_{\alpha }\left( x,u,\left( adj_{n-1}\nabla u\right) ^{\alpha }\right) %
\right] \,dx$ \\ 
$\leq \int\limits_{E_{k,s}^{1}}c\left( x\right) \left\vert \left(
adj_{n-1}\nabla u\right) ^{1}\right\vert ^{\delta }\left\vert \varphi
\right\vert ^{\beta }\,dx$ \\ 
$\leq \int\limits_{E_{k,s}^{1}}c\left( x\right) \left\vert \nabla
u\right\vert ^{\left( n-1\right) \delta }\left\vert \varphi \right\vert
^{\beta }\,dx$%
\end{tabular}%
\end{equation*}%
Since $0<\delta <\frac{r}{n-1}$, $0<\beta <\frac{p^{\ast }\left(
r-(n-1)\delta \right) }{r}$, $\sigma >\frac{rp^{\ast }}{(r-(n-1)\delta
)p^{\ast }-r\beta }$ then, using the H\"{o}lder Inequality, it follows%
\begin{equation}
\begin{tabular}{l}
$\int\limits_{E_{k,s}^{1}}c\left( x\right) \left\vert \nabla u\right\vert
^{\left( n-1\right) \delta }\left\vert \varphi \right\vert ^{\beta }\,dx$ \\ 
$\leq \left[ \int\limits_{E_{k,s}^{1}}\left\vert \nabla u\right\vert ^{r}\,dx%
\right] ^{\frac{(n-1)\delta }{r}}\left[ \int\limits_{E_{k,s}^{1}}\left(
c\left( x\right) \left\vert \varphi \right\vert ^{\beta }\right) ^{\frac{r}{%
r-(n-1)\delta }}\,dx\right] ^{\frac{r-(n-1)\delta }{r}}$ \\ 
$\leq \left[ \left\vert E_{k,s}^{1}\right\vert ^{1-\frac{r}{p}}\left(
\int\limits_{E_{k,s}^{1}}\left\vert \nabla u\right\vert ^{p}\,dx\right) ^{%
\frac{r}{p}}\right] ^{\frac{(n-1)\delta }{r}}$ \\ 
$\cdot \left[ \left( \int\limits_{E_{k,s}^{1}}\left( c\left( x\right)
\right) ^{\frac{rp^{\ast }}{\left( r-(n-1)\delta \right) p^{\ast }-r\beta }%
}\,dx\right) ^{\frac{\left( r-(n-1)\delta \right) p^{\ast }-r\beta }{\left(
r-(n-1)\delta \right) p^{\ast }}}\left( \int\limits_{E_{k,s}^{1}}\left\vert
\varphi \right\vert ^{p^{\ast }}\,dx\right) ^{\frac{r\beta }{\left(
r-(n-1)\delta \right) p^{\ast }}}\right] ^{\frac{r-(n-1)\delta }{r}}$ \\ 
$\leq \left\vert E_{k,s}^{1}\right\vert ^{\left( 1-\frac{r}{p}\right) \left( 
\frac{(n-1)\delta }{r}\right) }\left\Vert \nabla u\right\Vert _{L^{p}\left(
E_{k,s}^{1}\right) }^{(n-1)\delta }$ \\ 
$\cdot \left[ \left( \left\vert E_{k,s}^{1}\right\vert ^{1-\frac{rp^{\ast }}{%
\left[ \left( r-(n-1)\delta \right) p^{\ast }-r\beta \right] \sigma }%
}\left\Vert c\right\Vert _{L^{\sigma }\left( E_{k,s}^{1}\right) }^{\frac{%
rp^{\ast }}{\left[ \left( r-(n-1)\delta \right) p^{\ast }-r\beta \right] }%
}\right) ^{\frac{\left( r-(n-1)\delta \right) p^{\ast }-r\beta }{\left(
r-(n-1)\delta \right) p^{\ast }}}\left( \int\limits_{E_{k,s}^{1}}\left\vert
\varphi \right\vert ^{p^{\ast }}\,dx\right) ^{\frac{r\beta }{\left(
r-(n-1)\delta \right) p^{\ast }}}\right] ^{\frac{r-(n-1)\delta }{r}}$ \\ 
$\leq \left\vert E_{k,s}^{1}\right\vert ^{\left( 1-\frac{r}{p}\right) \left( 
\frac{(n-1)\delta }{r}\right) }\left\Vert \nabla u\right\Vert _{L^{p}\left(
E_{k,s}^{1}\right) }^{(n-1)\delta }$ \\ 
$\cdot \left[ \left\vert E_{k,s}^{1}\right\vert ^{\frac{\left( r-(n-1)\delta
\right) p^{\ast }-r\beta }{\left( r-(n-1)\delta \right) p^{\ast }}-\frac{%
rp^{\ast }}{\left[ \left( r-(n-1)\delta \right) p^{\ast }\right] \sigma }%
}\left\Vert c\right\Vert _{L^{\sigma }\left( E_{k,s}^{1}\right) }^{\frac{r}{%
\left( r-(n-1)\delta \right) }}\left( \int\limits_{E_{k,s}^{1}}\left\vert
\varphi \right\vert ^{p^{\ast }}\,dx\right) ^{\frac{r\beta }{\left(
r-(n-1)\delta \right) p^{\ast }}}\right] ^{\frac{r-(n-1)\delta }{r}}$ \\ 
$\leq \left\vert E_{k,s}^{1}\right\vert ^{\Lambda }\left\Vert \nabla
u\right\Vert _{L^{p}\left( E_{k,s}^{1}\right) }^{(n-1)\delta }\left\Vert
c\right\Vert _{L^{\sigma }\left( E_{k,s}^{1}\right) }\left(
\int\limits_{E_{k,s}^{1}}\left\vert \varphi \right\vert ^{p^{\ast
}}\,dx\right) ^{\frac{\beta }{p^{\ast }}}$%
\end{tabular}%
\end{equation}%
where $\Lambda =1-\left( \frac{(n-1)\delta }{p}+\frac{\beta }{p^{\ast }}+%
\frac{1}{\sigma }\right) $. Using the Sobolev Inequality we get%
\begin{equation*}
\begin{tabular}{l}
$\int\limits_{E_{k,s}^{1}}c\left( x\right) \left\vert \nabla u\right\vert
^{\left( n-1\right) \delta }\left\vert \varphi \right\vert ^{\beta }\,dx$ \\ 
$\leq C_{SI}^{\beta }\left\vert E_{k,s}^{1}\right\vert ^{\Lambda }\left\Vert
\nabla u\right\Vert _{L^{p}\left( E_{k,s}^{1}\right) }^{(n-1)\delta
}\left\Vert c\right\Vert _{L^{\sigma }\left( E_{k,s}^{1}\right) }\left(
\int\limits_{E_{k,s}^{1}}\left\vert \nabla \varphi \right\vert
^{p}\,dx\right) ^{\frac{\beta }{p}}$%
\end{tabular}%
\end{equation*}%
Since%
\begin{equation*}
\int\limits_{E_{k,s}^{1}}\left\vert \nabla \varphi \right\vert ^{p}\,dx\leq
\int\limits_{E_{k,s}^{1}}\eta ^{p}\left\vert \nabla u^{1}\right\vert
^{p}\,dx+2^{p}\int\limits_{E_{k,s}^{1}}\frac{\left( u^{1}-k\right) ^{p}}{%
\left( s-t\right) ^{p}}\,dx
\end{equation*}%
using Young Inequality it follows%
\begin{equation}
\begin{tabular}{l}
$\int\limits_{E_{k,s}^{1}}c\left( x\right) \left\vert \nabla u\right\vert
^{\left( n-1\right) \delta }\left\vert \varphi \right\vert ^{\beta }\,dx$ \\ 
$\leq \tilde{\varepsilon}\int\limits_{E_{k,s}^{1}}\eta ^{p}\left\vert \nabla
u^{1}\right\vert ^{p}\,dx+\tilde{\varepsilon}2^{p}\int\limits_{E_{k,s}^{1}}%
\frac{\left( u^{1}-k\right) ^{p}}{\left( s-t\right) ^{p}}\,dx$ \\ 
$+\frac{1}{\tilde{\varepsilon}^{\frac{1}{p-\beta }}}\left( C_{SI}^{\beta
}\left\vert E_{k,s}^{1}\right\vert ^{\Lambda }\left\Vert \nabla u\right\Vert
_{L^{p}\left( E_{k,s}^{1}\right) }^{(n-1)\delta }\left\Vert c\right\Vert
_{L^{\sigma }\left( E_{k,s}^{1}\right) }\right) ^{\frac{p}{p-\beta }}$ \\ 
$\leq \tilde{\varepsilon}\int\limits_{E_{k,s}^{1}}\eta ^{p}\left\vert \nabla
u^{1}\right\vert ^{p}\,dx+\tilde{\varepsilon}2^{p}\int\limits_{E_{k,s}^{1}}%
\frac{\left( u^{1}-k\right) ^{p}}{\left( s-t\right) ^{p}}\,dx$ \\ 
$+\frac{C_{\Sigma }}{\tilde{\varepsilon}^{\frac{1}{p-\beta }}}\left\vert
E_{k,s}^{1}\right\vert ^{\Lambda \frac{p}{p-\beta }}$%
\end{tabular}%
\end{equation}%
where%
\begin{equation*}
C_{\Sigma }=\left( C_{SI}^{\beta }\left\Vert \nabla u\right\Vert
_{L^{p}\left( \Sigma \right) }^{(n-1)\delta }\left\Vert c\right\Vert
_{L^{\sigma }\left( \Sigma \right) }\right) ^{\frac{p}{p-\beta }}
\end{equation*}%
Fixed $\tilde{\varepsilon}=\frac{1}{4}$, using (5.13) and (5.14), we get%
\begin{equation}
\begin{tabular}{l}
$\frac{1}{4}\int\limits_{E_{k,s}^{1}}\left\vert \nabla u^{1}\right\vert
^{p}+b_{1}\left( x\right) \left( u^{1}\right) ^{\gamma _{1}}\,dx$ \\ 
$\leq 2c\left( p,\gamma \right) \int\limits_{E_{k,s}^{1}}\left( 1-\eta
^{p}\right) \left( \left\vert \nabla u^{1}\right\vert ^{p}+b_{1}\left(
x\right) \left( u^{1}\right) ^{\gamma _{1}}\right) \,dx$ \\ 
$+\left( 2^{2p-1}p^{p}+2^{p+1}\right) \int\limits_{E_{k,s}^{1}}\frac{\left(
u^{1}-k\right) ^{p}}{\left( s-t\right) ^{p}}\,dx+2L_{1}\int%
\limits_{E_{k,s}^{1}}a_{1}\left( x\right) \,dx$ \\ 
$+\left\vert k\right\vert ^{p}R^{-\varepsilon N}\left\vert
A_{k,s}^{1}\right\vert ^{1-\frac{p}{N}+\varepsilon }$ \\ 
$+\int\limits_{E_{k,s}^{1}}\eta ^{p}\sum\limits_{\alpha =1}^{n}\left[
G_{\alpha }\left( x,u+\varphi ,\left( adj_{n-1}A\right) ^{\alpha }\right)
-G_{\alpha }\left( x,u+\varphi ,\left( adj_{n-1}\nabla u\right) ^{\alpha
}\right) \right] \,dx$ \\ 
$+\frac{C_{\Sigma }}{\left( \frac{1}{4}\right) ^{\frac{1}{p-\beta }}}%
\left\vert E_{k,s}^{1}\right\vert ^{\Lambda \frac{p}{p-\beta }}$%
\end{tabular}%
\end{equation}

Now we will have to give a new estimate of the integral%
\begin{equation*}
\int\limits_{E_{k,s}^{1}}\eta ^{p}\sum\limits_{\alpha =1}^{n}\left[
G_{\alpha }\left( x,u+\varphi ,\left( adj_{n-1}A\right) ^{\alpha }\right)
-G_{\alpha }\left( x,u+\varphi ,\left( adj_{n-1}\nabla u\right) ^{\alpha
}\right) \right] \,dx
\end{equation*}%
remembering that $G_{1}\left( x,u+\varphi ,\left( adj_{n-1}A\right)
^{1}\right) =G_{1}\left( x,u+\varphi ,\left( adj_{n-1}\nabla u\right)
^{1}\right) $ and that $G_{\alpha }\left( x,s,\left( adj_{n-1}\xi \right)
^{\alpha }\right) \geq 0$ for every $x\in \Omega ,s\in 
\mathbb{R}
^{n}$ and $\xi \in 
\mathbb{R}
^{N\times n}$ we get%
\begin{equation*}
\begin{tabular}{l}
$\int\limits_{E_{k,s}^{1}}\eta ^{p}\sum\limits_{\alpha =1}^{n}\left[
G_{\alpha }\left( x,u+\varphi ,\left( adj_{n-1}A\right) ^{\alpha }\right)
-G_{\alpha }\left( x,u+\varphi ,\left( adj_{n-1}\nabla u\right) ^{\alpha
}\right) \right] \,dx$ \\ 
$\leq \int\limits_{E_{k,s}^{1}}\eta ^{p}\sum\limits_{\alpha =2}^{n}G_{\alpha
}\left( x,u+\varphi ,\left( adj_{n-1}A\right) ^{\alpha }\right) \,dx$%
\end{tabular}%
\end{equation*}%
Using from hypothesis H.2.3, it follows that%
\begin{equation*}
\begin{tabular}{l}
$\int\limits_{E_{k,s}^{1}}\eta ^{p}\sum\limits_{\alpha =1}^{n}\left[
G_{\alpha }\left( x,u+\varphi ,\left( adj_{n-1}A\right) ^{\alpha }\right)
-G_{\alpha }\left( x,u+\varphi ,\left( adj_{n-1}\nabla u\right) ^{\alpha
}\right) \right] \,dx$ \\ 
$\leq \int\limits_{E_{k,s}^{1}}\eta ^{p}\sum\limits_{\alpha
=2}^{n}\left\vert \left( adj_{n-1}A\right) ^{\alpha }\right\vert
^{r}+\sum\limits_{\alpha =2}^{n}\eta ^{p}c\left( x\right) \left\vert \left(
adj_{n-1}A\right) ^{\alpha }\right\vert ^{\delta }\left\vert u+\varphi
\right\vert ^{\beta }+\eta ^{p}a\left( x\right) \,dx$%
\end{tabular}%
\end{equation*}%
We use here the following notation:%
\begin{equation*}
\hat{u}=\left( u^{2},...,u^{n}\right)
\end{equation*}%
thus, for all $\alpha =2,...,n$%
\begin{equation*}
\left\vert \left( adj_{n-1}A\right) ^{\alpha }\right\vert \leq \tilde{c}%
\left\vert A^{1}\right\vert \left\vert adj_{n-2}\nabla \hat{u}\right\vert
\end{equation*}%
hence remembering that%
\begin{equation*}
\left\vert adj_{n-2}\nabla \hat{u}\right\vert \leq \tilde{c}_{2}\left\vert
\nabla \hat{u}\right\vert ^{n-2}
\end{equation*}%
we get%
\begin{equation*}
\begin{tabular}{l}
$\int\limits_{E_{k,s}^{1}}\eta ^{p}\sum\limits_{\alpha =1}^{n}\left[
G_{\alpha }\left( x,u+\varphi ,\left( adj_{n-1}A\right) ^{\alpha }\right)
-G_{\alpha }\left( x,u+\varphi ,\left( adj_{n-1}\nabla u\right) ^{\alpha
}\right) \right] \,dx$ \\ 
$\leq c_{3}\int\limits_{E_{k,s}^{1}}\eta ^{p}\left\vert A^{1}\right\vert
^{r}\left\vert \nabla \hat{u}\right\vert ^{\left( n-2\right) r}+\eta
^{p}c\left( x\right) \left\vert A^{1}\right\vert ^{\delta }\left\vert \nabla 
\hat{u}\right\vert ^{\left( n-2\right) \delta }\left\vert u+\varphi
\right\vert ^{\beta }+\eta ^{p}a\left( x\right) \,dx$%
\end{tabular}%
\end{equation*}%
where $c_{3}=\left[ 1+\left( \tilde{c}\tilde{c}_{2}\right) ^{r}+\left( 
\tilde{c}\tilde{c}_{2}\right) ^{\delta }\right] n$. We consider the integral 
$\int\limits_{E_{k,s}^{1}}\eta ^{p}\left\vert A^{1}\right\vert
^{r}\left\vert \nabla \hat{u}\right\vert ^{\left( n-2\right) r}\,dx$\ since $%
r<\frac{p}{n-1}$\ using H\"{o}lder Inequality and Young Inequality it folows%
\begin{equation*}
\begin{tabular}{l}
$\int\limits_{E_{k,s}^{1}}\eta ^{p}\left\vert A^{1}\right\vert
^{r}\left\vert \nabla \hat{u}\right\vert ^{\left( n-2\right) r}\,dx$ \\ 
$\leq \left[ \int\limits_{E_{k,s}^{1}}\eta ^{p}\left\vert A^{1}\right\vert
^{p}\,dx\right] ^{\frac{r}{p}}\left[ \int\limits_{E_{k,s}^{1}}\eta
^{p}\left\vert \nabla \hat{u}\right\vert ^{\frac{\left( n-2\right) rp}{p-r}%
}\,dx\right] ^{\frac{p-r}{p}}$ \\ 
$\leq \left[ \int\limits_{E_{k,s}^{1}}\eta ^{p}\left\vert A^{1}\right\vert
^{p}\,dx\right] ^{\frac{r}{p}}\left[ \left\vert E_{k,s}^{1}\right\vert ^{1-%
\frac{\left( n-2\right) r}{p-r}}\left( \int\limits_{E_{k,s}^{1}}\eta
^{p}\left\vert \nabla u\right\vert ^{p}\,dx\right) ^{\frac{\left( n-2\right)
r}{p-r}}\right] ^{\frac{p-r}{p}}$ \\ 
$\leq \int\limits_{E_{k,s}^{1}}\eta ^{p}\left\vert A^{1}\right\vert
^{p}\,dx+\left\vert E_{k,s}^{1}\right\vert ^{1-\frac{\left( n-2\right) r}{p-r%
}}\left( \int\limits_{E_{k,s}^{1}}\eta ^{p}\left\vert \nabla u\right\vert
^{p}\,dx\right) ^{\frac{\left( n-2\right) r}{p-r}}$ \\ 
$\leq \int\limits_{E_{k,s}^{1}}\eta ^{p}\left\vert A^{1}\right\vert
^{p}\,dx+\left\vert E_{k,s}^{1}\right\vert ^{1-\frac{\left( n-2\right) r}{p-r%
}}\left\Vert \nabla u\right\Vert _{L^{p}\left( \Sigma \right) }^{\frac{%
\left( n-2\right) rp}{p-r}}$%
\end{tabular}%
\end{equation*}

\bigskip Since $A^{1}=p\eta ^{-1}\nabla \eta \left( k-u^{1}\right) $ we get%
\begin{equation*}
\begin{tabular}{l}
$\int\limits_{E_{k,s}^{1}}\eta ^{p}\left\vert A^{1}\right\vert
^{r}\left\vert \nabla \hat{u}\right\vert ^{\left( n-2\right) r}\,dx$ \\ 
$\leq 2^{p}\int\limits_{E_{k,s}^{1}}\frac{\left( u^{1}-k\right) ^{p}}{\left(
s-t\right) ^{p}}\,dx+\left\vert E_{k,s}^{1}\right\vert ^{1-\frac{\left(
n-2\right) r}{p-r}}\left\Vert \nabla u\right\Vert _{L^{p}\left( \Sigma
\right) }^{\frac{\left( n-2\right) rp}{p-r}}$%
\end{tabular}%
\end{equation*}%
Now we need to estimate the integral $\int\limits_{E_{k,s}^{1}}\eta
^{p}c\left( x\right) \left\vert A^{1}\right\vert ^{\delta }\left\vert \nabla 
\hat{u}\right\vert ^{\left( n-2\right) \delta }\left\vert u+\varphi
\right\vert ^{\beta }\,dx$. Using the H\"{o}lder inequality we obtain%
\begin{equation*}
\begin{tabular}{l}
$\int\limits_{E_{k,s}^{1}}\eta ^{p}c\left( x\right) \left\vert
A^{1}\right\vert ^{\delta }\left\vert \nabla \hat{u}\right\vert ^{\left(
n-2\right) \delta }\left\vert u+\varphi \right\vert ^{\beta }\,dx$ \\ 
$\leq \left[ \int\limits_{E_{k,s}^{1}}\eta ^{r}\left\vert A^{1}\right\vert
^{r}\left\vert \nabla \hat{u}\right\vert ^{\left( n-2\right) r}\,dx\right] ^{%
\frac{\delta }{r}}\left[ \int\limits_{E_{k,s}^{1}}\left( \eta ^{p-\delta
}c\left( x\right) \left\vert u+\varphi \right\vert ^{\beta }\right) ^{\frac{r%
}{r-\delta }}\,dx\right] ^{\frac{r-\delta }{r}}$ \\ 
$\leq \left[ \left( \int\limits_{E_{k,s}^{1}}\eta ^{p}\left\vert
A^{1}\right\vert ^{p}\,dx\right) ^{\frac{r}{p}}\left(
\int\limits_{E_{k,s}^{1}}\left\vert \nabla \hat{u}\right\vert ^{\frac{\left(
n-2\right) rp}{p-r}}\,dx\right) ^{\frac{p-r}{p}}\right] ^{\frac{\delta }{r}}$
\\ 
$\cdot \left[ \int\limits_{E_{k,s}^{1}}\left( \eta ^{p-\delta }c\left(
x\right) \right) ^{\frac{r}{r-\delta }}\left\vert u+\varphi \right\vert ^{%
\frac{\beta r}{r-\delta }}\,dx\right] ^{\frac{r-\delta }{r}}$ \\ 
$\leq \left[ \left[ \int\limits_{E_{k,s}^{1}}\eta ^{p}\left\vert
A^{1}\right\vert ^{p}\,dx\right] ^{\frac{r}{p}}\left[ \left\vert
E_{k,s}^{1}\right\vert ^{1-\frac{\left( n-2\right) r}{p-r}}\left(
\int\limits_{E_{k,s}^{1}}\left\vert \nabla u\right\vert ^{p}\,dx\right) ^{%
\frac{\left( n-2\right) r}{p-r}}\right] ^{\frac{p-r}{p}}\right] ^{\frac{%
\delta }{r}}$ \\ 
$\cdot \left[ \left( \int\limits_{E_{k,s}^{1}}\left( \eta ^{p-\delta
}c\left( x\right) \right) ^{\frac{r\left( r-\delta \right) p^{\ast }}{\left(
r-\delta \right) \left[ \left( r-\delta \right) p^{\ast }-\beta r\right] }%
}\,dx\right) ^{\frac{\left[ \left( r-\delta \right) p^{\ast }-\beta r\right] 
}{\left( r-\delta \right) p^{\ast }}}\left(
\int\limits_{E_{k,s}^{1}}\left\vert u+\varphi \right\vert ^{p^{\ast
}}\,dx\right) ^{\frac{\beta r}{\left( r-\delta \right) p^{\ast }}}\right] ^{%
\frac{r-\delta }{r}}$ \\ 
$\leq \left[ \left[ \int\limits_{E_{k,s}^{1}}\eta ^{p}\left\vert
A^{1}\right\vert ^{p}\,dx\right] ^{\frac{r}{p}}\left[ \left\vert
E_{k,s}^{1}\right\vert ^{1-\frac{\left( n-2\right) r}{p-r}}\left(
\int\limits_{E_{k,s}^{1}}\left\vert \nabla u\right\vert ^{p}\,dx\right) ^{%
\frac{\left( n-2\right) r}{p-r}}\right] ^{\frac{p-r}{p}}\right] ^{\frac{%
\delta }{r}}$ \\ 
$\cdot \left[ \left( \int\limits_{E_{k,s}^{1}}\left( \eta ^{p-\delta
}c\left( x\right) \right) ^{\frac{rp^{\ast }}{\left[ \left( r-\delta \right)
p^{\ast }-\beta r\right] }}\,dx\right) ^{\frac{\left[ \left( r-\delta
\right) p^{\ast }-\beta r\right] }{\left( r-\delta \right) p^{\ast }}}\left(
\left\Vert u\right\Vert _{L^{p^{\ast }}\left( B_{R_{0}}\left( x_{0}\right)
\right) }+\left\Vert \varphi \right\Vert _{L^{p^{\ast }}\left( B_{s}\left(
x_{0}\right) \right) }\right) ^{\frac{\beta r}{\left( r-\delta \right) }}%
\right] ^{\frac{r-\delta }{r}}$%
\end{tabular}%
\end{equation*}%
Since $A^{1}=p\eta ^{-1}\nabla \eta \left( k-u^{1}\right) $ and using the
Sobolev Inequaility and The Sobolev Embedding Theorem we get%
\begin{equation*}
\begin{tabular}{l}
$\int\limits_{E_{k,s}^{1}}\eta ^{p}c\left( x\right) \left\vert
A^{1}\right\vert ^{\delta }\left\vert \nabla \hat{u}\right\vert ^{\left(
n-2\right) \delta }\left\vert u+\varphi \right\vert ^{\beta }\,dx$ \\ 
$\leq \left[ \left[ \int\limits_{E_{k,s}^{1}}p^{p}\left\vert \nabla \eta
\right\vert ^{p}\left( u^{1}-k\right) ^{p}\,dx\right] ^{\frac{r}{p}}\left[
\left\vert E_{k,s}^{1}\right\vert ^{1-\frac{\left( n-2\right) r}{p-r}}\left(
\int\limits_{E_{k,s}^{1}}\left\vert \nabla u\right\vert ^{p}\,dx\right) ^{%
\frac{\left( n-2\right) r}{p-r}}\right] ^{\frac{p-r}{p}}\right] ^{\frac{%
\delta }{r}}$ \\ 
$\cdot \left[ \left( \int\limits_{E_{k,s}^{1}}\left( \eta ^{p-\delta
}c\left( x\right) \right) ^{\frac{rp^{\ast }}{\left[ \left( r-\delta \right)
p^{\ast }-\beta r\right] }}\,dx\right) ^{\frac{\left[ \left( r-\delta
\right) p^{\ast }-\beta r\right] }{\left( r-\delta \right) p^{\ast }}}\left(
C_{ES}\left\Vert u\right\Vert _{W^{1,p}\left( \Sigma ,%
\mathbb{R}
^{n}\right) }+C_{IS}\left\Vert \nabla \varphi \right\Vert _{L^{p}\left(
B_{s}\left( x_{0}\right) \right) }\right) ^{\frac{\beta r}{\left( r-\delta
\right) }}\right] ^{\frac{r-\delta }{r}}$%
\end{tabular}%
\end{equation*}%
By H\"{o}lder Inequality it follows%
\begin{equation*}
\begin{tabular}{l}
$\int\limits_{E_{k,s}^{1}}\eta ^{p}c\left( x\right) \left\vert
A^{1}\right\vert ^{\delta }\left\vert \nabla \hat{u}\right\vert ^{\left(
n-2\right) \delta }\left\vert u+\varphi \right\vert ^{\beta }\,dx$ \\ 
$\leq \left[ \left[ \int\limits_{E_{k,s}^{1}}p^{p}\left\vert \nabla \eta
\right\vert ^{p}\left( u^{1}-k\right) ^{p}\,dx\right] ^{\frac{r}{p}}\left[
\left\vert E_{k,s}^{1}\right\vert ^{1-\frac{\left( n-2\right) r}{p-r}}\left(
\int\limits_{E_{k,s}^{1}}\left\vert \nabla u\right\vert ^{p}\,dx\right) ^{%
\frac{\left( n-2\right) r}{p-r}}\right] ^{\frac{p-r}{p}}\right] ^{\frac{%
\delta }{r}}$ \\ 
$\cdot \left[ \left( \left\vert E_{k,s}^{1}\right\vert ^{1-\frac{rp^{\ast }}{%
\left[ \left( r-\delta \right) p^{\ast }-r\beta \right] \sigma }}\left\Vert
c\right\Vert _{L^{\sigma }\left( E_{k,s}^{1}\right) }^{\frac{rp^{\ast }}{%
\left[ \left( r-\delta \right) p^{\ast }-r\beta \right] }}\right) ^{\frac{%
\left[ \left( r-\delta \right) p^{\ast }-\beta r\right] }{\left( r-\delta
\right) p^{\ast }}}\right. $ \\ 
$\left. \cdot \left( C_{ES}\left\Vert u\right\Vert _{W^{1,p}\left( \Sigma ,%
\mathbb{R}
^{n}\right) }+C_{IS}\left\Vert \nabla \varphi \right\Vert _{L^{p}\left(
B_{s}\left( x_{0}\right) \right) }\right) ^{\frac{\beta r}{\left( r-\delta
\right) }}\right] ^{\frac{r-\delta }{r}}$%
\end{tabular}%
\end{equation*}%
moreover 
\begin{equation*}
\int\limits_{E_{k,s}^{1}}\left\vert \nabla \varphi \right\vert ^{p}\,dx\leq
\int\limits_{E_{k,s}^{1}}\eta ^{p}\left\vert \nabla u^{1}\right\vert
^{p}\,dx+2^{p}\int\limits_{E_{k,s}^{1}}\frac{\left( u^{1}-k\right) ^{p}}{%
\left( s-t\right) ^{p}}\,dx
\end{equation*}%
then we have%
\begin{equation*}
\begin{tabular}{l}
$\int\limits_{E_{k,s}^{1}}\eta ^{p}c\left( x\right) \left\vert
A^{1}\right\vert ^{\delta }\left\vert \nabla \hat{u}\right\vert ^{\left(
n-2\right) \delta }\left\vert u+\varphi \right\vert ^{\beta }\,dx$ \\ 
$\leq \left[ \left[ p^{p}2^{p}\int\limits_{E_{k,s}^{1}}\frac{\left(
u^{1}-k\right) ^{p}}{\left( s-t\right) ^{p}}\,dx\right] ^{\frac{r}{p}}\left[
\left\vert E_{k,s}^{1}\right\vert ^{1-\frac{\left( n-2\right) r}{p-r}}\left(
\int\limits_{E_{k,s}^{1}}\left\vert \nabla u\right\vert ^{p}\,dx\right) ^{%
\frac{\left( n-2\right) r}{p-r}}\right] ^{\frac{p-r}{p}}\right] ^{\frac{%
\delta }{r}}$ \\ 
$\cdot \left[ \left( \left\vert E_{k,s}^{1}\right\vert ^{1-\frac{rp^{\ast }}{%
\left[ \left( r-\delta \right) p^{\ast }-r\beta \right] \sigma }}\left\Vert
c\right\Vert _{L^{\sigma }\left( E_{k,s}^{1}\right) }^{\frac{rp^{\ast }}{%
\left[ \left( r-\delta \right) p^{\ast }-r\beta \right] }}\right) ^{\frac{%
\left[ \left( r-\delta \right) p^{\ast }-\beta r\right] }{\left( r-\delta
\right) p^{\ast }}}\right. $ \\ 
$\left. \cdot \left( C_{ES}\left\Vert u\right\Vert _{W^{1,p}\left( \Sigma ,%
\mathbb{R}
^{n}\right) }+C_{IS}\left( \int\limits_{E_{k,s}^{1}}\eta ^{p}\left\vert
\nabla u^{1}\right\vert ^{p}\,dx+2^{p}\int\limits_{E_{k,s}^{1}}\frac{\left(
u^{1}-k\right) ^{p}}{\left( s-t\right) ^{p}}\,dx\right) ^{\frac{1}{p}%
}\right) ^{\frac{\beta r}{\left( r-\delta \right) }}\right] ^{\frac{r-\delta 
}{r}}$%
\end{tabular}%
\end{equation*}%
Using the Young Inequality we obtain%
\begin{equation*}
\begin{tabular}{l}
$\int\limits_{E_{k,s}^{1}}\eta ^{p}c\left( x\right) \left\vert
A^{1}\right\vert ^{\delta }\left\vert \nabla \hat{u}\right\vert ^{\left(
n-2\right) \delta }\left\vert u+\varphi \right\vert ^{\beta }\,dx$ \\ 
$\leq \left\Vert c\right\Vert _{L^{\sigma }\left( E_{k,s}^{1}\right)
}\left\vert E_{k,s}^{1}\right\vert ^{\Lambda }$ \\ 
$\cdot \left[ p^{p}2^{p}\int\limits_{E_{k,s}^{1}}\frac{\left( u^{1}-k\right)
^{p}}{\left( s-t\right) ^{p}}\,dx+C\left\Vert u\right\Vert _{W^{1,p}\left(
\Sigma ,%
\mathbb{R}
^{n}\right) }+1\right] ^{\frac{\delta }{r}}$ \\ 
$\cdot \left[ C\left\Vert u\right\Vert _{W^{1,p}\left( \Sigma ,%
\mathbb{R}
^{n}\right) }+1+2^{p}\int\limits_{E_{k,s}^{1}}\frac{\left( u^{1}-k\right)
^{p}}{\left( s-t\right) ^{p}}\,dx\right] ^{\frac{r-\delta }{r}}$ \\ 
$\leq 2\left\Vert c\right\Vert _{L^{\sigma }\left( \Sigma \right)
}\left\vert E_{k,s}^{1}\right\vert ^{\Lambda }$ \\ 
$\cdot \left[ p^{p}2^{p}\int\limits_{E_{k,s}^{1}}\frac{\left( u^{1}-k\right)
^{p}}{\left( s-t\right) ^{p}}\,dx+C\left\Vert u\right\Vert _{W^{1,p}\left(
\Sigma ,%
\mathbb{R}
^{n}\right) }+1\right] $ \\ 
$\leq p^{p}2^{p+1}\left\Vert c\right\Vert _{L^{\sigma }\left( \Sigma \right)
}\int\limits_{E_{k,s}^{1}}\frac{\left( u^{1}-k\right) ^{p}}{\left(
s-t\right) ^{p}}\,dx+C\left( \Sigma \right) \left\vert
E_{k,s}^{1}\right\vert ^{\Lambda }$%
\end{tabular}%
\end{equation*}%
where $\Lambda =1-\left( \frac{(n-1)\delta }{p}+\frac{\beta }{p^{\ast }}+%
\frac{1}{\sigma }\right) $. Proceeding as already done in the previous
section 3.1 it follows 
\begin{equation}
\begin{tabular}{l}
$\frac{1}{4}\int\limits_{E_{k,s}^{1}}\left\vert \nabla u^{1}\right\vert
^{p}+b_{1}\left( x\right) \left( u^{1}\right) ^{\gamma _{1}}\,dx$ \\ 
$\leq 2c\left( p,\gamma \right) \int\limits_{E_{k,s}^{1}}\left( 1-\eta
^{p}\right) \left( \left\vert \nabla u^{1}\right\vert ^{p}+b_{1}\left(
x\right) \left( u^{1}\right) ^{\gamma _{1}}\right) \,dx$ \\ 
$+\left( 2^{2p-1}p^{p}+2^{p}\right) \int\limits_{E_{k,s}^{1}}\frac{\left(
u^{1}-k\right) ^{p}}{\left( s-t\right) ^{p}}\,dx+2L_{1}\int%
\limits_{E_{k,s}^{1}}a_{1}\left( x\right) \,dx$ \\ 
$+\left\vert k\right\vert ^{p}R^{-\varepsilon N}\left\vert
A_{k,s}^{1}\right\vert ^{1-\frac{p}{N}+\varepsilon }+C_{1,\Sigma }\left\vert
E_{k,s}^{1}\right\vert ^{\Lambda }+2C_{2,\Sigma
}\int\limits_{E_{k,s}^{1}-E_{k,t}^{1}}\frac{\left( u^{1}-k\right) ^{p}}{%
\left( s-t\right) ^{p}}\,dx$%
\end{tabular}%
\end{equation}%
From this last inequality, (5.16), proceeding as already done in the
previous section 3.1 by applying Lemma 3 we obtain the inequalities of
Cacioppoli (1.14) and (1.15).\bigskip

\subsection{\protect\bigskip The case H.2.4 (bis)}

\bigskip The proof follows as in the previous section 4.1 only with some
modifications in the choice of the parameters of the H\"{o}lder inequalities.

\section{Proof of Theorem 4}

Let $u\in W^{1,p}\left( \Omega ;%
\mathbb{R}
^{n}\right) $ a minimizer of the functional (1.1) then by Theorem 2 it
follows that $u^{\alpha }\in DG\left( \Omega ,p,\lambda ,\lambda _{\ast
},\chi ,\varepsilon ,R_{0},k_{0}\right) $ for $\alpha =1,...,n$ \ then by
Theorem 6 it follows that $u\in L_{loc}^{\infty }\left( \Omega ;%
\mathbb{R}
^{n}\right) $; moreover if $\Sigma $ is a compac subset of $\Omega $ then $%
u\in L^{\infty }\left( \Sigma ;%
\mathbb{R}
^{n}\right) $ and%
\begin{equation}
\left\vert u\right\vert \leq M=\sqrt{\sum\limits_{\alpha =1}^{n}\left(
M^{\alpha }\right) ^{2}}
\end{equation}%
where%
\begin{equation}
M^{\alpha }=\sup_{\Sigma }\left\{ \left\vert u^{\alpha }\right\vert \right\}
\end{equation}%
Since H.1.1 and (1.2) hold then proceeding as in [22] we get%
\begin{equation}
\begin{tabular}{l}
$\int\limits_{A_{k,\varrho }^{\alpha }}\left\vert \nabla u^{\alpha
}\right\vert ^{p}\,dx\leq \frac{\tilde{C}_{C,1}}{\left( R-\varrho \right)
^{p}}\int\limits_{A_{k,R}^{\alpha }}\left( u^{\alpha }-k\right) ^{p}\,dx+%
\tilde{C}_{C,2}\left\vert A_{k,s}^{\alpha }\right\vert ^{1-\frac{p}{N}%
+\varepsilon }$%
\end{tabular}%
\end{equation}%
and%
\begin{equation}
\int\limits_{B_{k,\varrho }^{\alpha }}\left\vert \nabla u^{\alpha
}\right\vert ^{p}\,dx\leq \frac{\tilde{C}_{C,1}}{\left( R-\varrho \right)
^{p}}\int\limits_{B_{k,R}^{\alpha }}\left( k-u^{\alpha }\right) ^{p}\,dx+%
\tilde{C}_{C,2}\left\vert B_{k,s}^{\alpha }\right\vert ^{1-\frac{p}{N}%
+\varepsilon }
\end{equation}%
for every $\alpha =1,...,n$. Theorem 1 follows using (6.3), (6.4) and
Proposition 7.1, Lemma 7.2 and Theorem 7.6 of [24].


\bigskip

\begin{thebibliography}{10}
\bibitem[1]{1} Acerbi E., Fusco N., Regularity for minimizers of
non-quadratic functionals:\ the case $1<p<2$, \textit{J. Math. Anal. Appl.},
140, 1989, 115-134

\bibitem[2]{2} Acerbi E., Fusco N., Partial regularity under anisotropic $%
(p,q)$ growth conditions, \textit{J. Diff. Equ.}, 107, no 1, 1994, 46-67

\bibitem[3]{3} Aubin T., Probl\`{e}mes isop\'{e}rim\'{e}triques et espaces
de Sobolev, \textit{Journal of Differential Geometry}, 11, 1976, 573-598

\bibitem[4]{4} Bildhauer M., Fuchs M., Partial regularity for variational
integrals with $(s,\mu ,q)$-growth, \textit{Calc. Var.}, 13, 2001, 537-560

\bibitem[5]{5} Bildhauer M., Fuchs M., Mingione G., A priori gradient bounds
and local $C^{1,\alpha }$-estimates for (doble) obstacle problems under
nonstandard growth conditions, \textit{Z. Anal. Anw.}, 20, no 4, 2001,
959-985

\bibitem[6]{6} Bildhauer M., \textit{Convex Variational Problems, Linear,
Nearly Linear and Anisotropic Growth Conditions}, Springer, Berlin, 2003

\bibitem[7]{7} Breit D., Stroffolini B., Verde A., A general regularity
theorem for functionals with $\varphi $-growth, \textit{J. Math. Anal. Appl.}%
, 383 (2011), 226-233. https://doi.org/10.1016/j.jmaa.2011.05.012

\bibitem[8]{8} Cupini G., Focardi M., Leonetti E., Mascolo F., On the Holder
continuity for a class of vectorial problems, \textit{Advances in Nonlinear
Analysis}, 9 (2020), no. 1, 1008-1025.https://doi.org/10.1515/anona-2020-0039

\bibitem[9]{9} Cupini G., Leonetti F., Mascolo E., Local boundedness for
minimizers of some polyconvex integrals, \textit{Arch. Rational Mech. Anal.}%
, 224 (2017), no.1, 269-289.

\bibitem[10]{10} Cupini G., Marcellini P., Mascolo E., Local boundedness of
solutions to some anisotropic elliptic systems, \textit{Contemporary
Mathematics}, 595 (2013), 169-186.

\bibitem[11]{11} Cupini G., Marcellini P., Mascolo E., Local boundedness of
solutions to quasilinear elliptic systems, \textit{Manuscripta Mathematica},
137, 287-315, 2012

\bibitem[12]{12} Dacorogna B., \textit{Direct Methods in the Calculus of
Variations}, Springer-Verlag, 1989

\bibitem[13]{13} De Giorgi E., Sulla differenziabilit`a e l'analicit`a delle
estremali degli integrali multipli regolari, \textit{Mem. Accad. Sci. Torino
(Classe di Sci. mat. fis. e nat.)}, 3 (3) (1957), 25-43.

\bibitem[14]{14} De Giorgi E., Un esempio di estremali discontinue per un
problema variazionale di tipo ellittico, \textit{Boll. U.M.I.}, 4 (1968),
135-137.

\bibitem[15]{15} Diening L., Stroffolini B., Verde A., Everywhere regularity
of functional with $\varphi $-growth, \textit{Manus. Math.}, 129 (2009),
440-481.

\bibitem[16]{16} Evans L. C., Quasiconvexity and partial regularity in the
calculus of variations, \textit{Arch. Raz. Mech. Anal.}, 95, 1986, 227-252

\bibitem[17]{17} Esposito L., Mingione G., Some remarks on the regulariy of
weak solutions of degenerate elliptic systems, \textit{Rev. Mat. Complu.},
11, no 1, 1998, 203-219

\bibitem[18]{18} Esposito L., Mingione G., Partial regularity for minimizers
of convex integrals with $L\,log\,L$ -growth, \textit{Nonlinear Diff. Equ.
Appl.}, 7, 2000, 107-125

\bibitem[19]{19} Frehse J., A discontinuous solution of a mildly nonlinear
system, \textit{Math. Z.}, 124 (1973), 229-230.
https://doi.org/10.1007/bf01214096

\bibitem[20]{20} Fuchs M., Serengin G., A regularity theory for variational
integrals with $L\,log\,L$ -growth, \textit{Cal. Var.}, 6, 1998, 171-187\ 

\bibitem[21]{21} Fuchs M., Local Lipschitz regularity of vector valued local
minimizers of variational integrals with densities depending on the modulus
of the gradient, \textit{Math. Nachr.}, 284 (2011), 266-272.

\bibitem[22]{22} Giaquinta M., Giusti E., On the regularity of minima of
variational integrals, \textit{Acta Mathematica}, 148 (1983), 285-298.

\bibitem[23]{23} Giusti E., Miranda M., Un esempio di soluzioni discontinue
per un problema di minimo relativo ad un integrale regolare del calcolo
delle variazioni, \textit{Boll. U.M.I.}, 2 (1968), 1-8

\bibitem[24]{24} Giusti E., \textit{Direct methods in the Calculus of
variations}, World Scientific, 2003

\bibitem[25]{25} Granucci T., Randolfi M., Regularity for local minima of a
special class of vectorial problems with fully anisotropic growth,
manuscripta math. (2022). https://doi.org/10.1007/s00229-021-01360-0

\bibitem[26]{26} Granucci T., Randolfi M.: Local boundedness of
quasi-minimizers of fully anisotropic scalar variational problems. Manuscr.
Math. 160, 99--152 (2019)

\bibitem[27]{27} Granucci T., An Harnack inequality for quasi-minima of
scalar integral functionals with general growth conditions. Manuscr. Math.
152, 345--380 (2017)

\bibitem[28]{28} Granucci T., L$^{\Phi }$-L$^{\infty }$ inequalities and new
remarks on the Ho%
\"{}%
lder continuity of the quasi-minima of scalar integral functionals with
general growths. Bol. Soc. Mat. Mex. 22, 165--212 (2016)

\bibitem[29]{29} Granucci T., An Example of Everywhere Regularity for Minima
of Vectorial Integral Fuctionals of the Calculus of Variation, International
Journal of Mathematical Analysis, Vol. 15, 2021, no. 7, 291 - 302

\bibitem[30]{30} Granucci T.,\ An Example of Everywhere H\"{o}lder
Continuity for a Weak Solution of a Particular Vectorial Problem,
International Journal of Mathematical Analysis, Vol. 14, 2020, no. 7, 305 -
313

\bibitem[31]{31} Granucci T., On the everywhere h\"{o}lder continuity of the
minima of a class of vectorial integral functionals of the calculus of
variation, submitted

\bibitem[32]{32} Granucci T., Everywhere h\"{o}lder continuity of vectorial
local minimizers of special classes of integral functionals with rank one
integrands, Monatsh Math (2022). https://doi.org/10.1007/s00605-022-01763-5

\bibitem[33]{33} Marcellini P., Everywhere regularity for a class of
elliptic systems without growth conditions, \textit{Ann. Sc. Norm. Super.
Pisa}, 23 (1996), 1-25.

\bibitem[34]{34} Mingione G., Singularities of minima: a walk on the wild
side of the calculus of variations, \textit{Journal of global optimization},
40, (2008) 209-223

\bibitem[35]{35} Mingione G., Regularity of minima: an invitation to the
dark side of the calculus of variations, \textit{Applications of mathematics}%
, 51, (2006) 355-426

\bibitem[36]{36} Morrey C. B., Partial regularity results for nonlinear
elliptic systems, \textit{J. Math. Mech.}, 17, 1968, 649-670

\bibitem[37]{37} Moser J., A new proof of De Giorgi's theorem concerning the
regularity problem for elliptic differential equations, \textit{Comm. Pure
Appl. Math.}, 14 (1961), 457-468. https://doi.org/10.1002/cpa.3160130308

\bibitem[38]{38} Nash J., Continuity of solution of parabolic and elliptic
equations, \textit{Amer. J. of Math.}, 80 (1958), 931-954.
https://doi.org/10.2307/2372841

\bibitem[39]{39} Sobolev S. L., Sur un th\'{e}or\`{e}me d'analyse
fonctionnelle, \textit{Math. Sb. (N.S)}, 46, 1938, 471-496

\bibitem[40]{40} Talenti G., Best constants in Sobolev inequality, \textit{%
Ann. di Matem. Pura ed Apll.}, 110, 1976, 353-372

\bibitem[41]{41} Tolksdorf P., A new proof of a regularity theorem, \textit{%
Invent. Math.}, 71 (1983), no.1, 43-49.

\bibitem[42]{42} Tolksdorf P., Regularity for a More General Class of
Quasilinear Elliptic Equations, \textit{J. of Differ. Equ.}, 51 (1984),
126-150.

\bibitem[43]{43} Uhlenbeck K., Regularity \ for a class of nonlinear
elliptic systems, \textit{Acta Math.}, 138 (1977), 219-240.
https://doi.org/10.1007/bf02392316
\end{thebibliography}
\end{document}